\documentclass[10pt]{amsart} \usepackage{amssymb, verbatim, epic, eepic}

\def\Re{{\hbox{Re}}}
\def\diam{{\hbox{diam}}} 

\def\C{{\hbox{\bf C}}}
\def\R{{{\mathbf{R}}}}

\def\eps{\varepsilon}
\def \endprf{\qed}
\def\emph#1{{\it #1}} \def\textbf#1{{\bf #1}}
 
\newcommand{\nabb}{\mbox{$\nabla \mkern-13mu /$\,}}
\def\Ad{\mathcal{A}_{\delta}}
\def\CId{C^\infty_\delta}
\def\half{\frac1{2}}

\def\Ginit{G_{\operatorname{init}}}
\def\Cgamma{\mathcal{C}_{\gamma}}

\def\RR{\mathbf{R}}

\def\ang{{\operatorname{ang}}}

\def\Tscstar{{}^{\sc} T^*} 
\def\sc{\operatorname{sc}}
\def\Psisc{\Psi_{\sc}} 
\def\Id{\operatorname{Id}}
\def\bface{\operatorname{bf}}
\def\bfc{\operatorname{bf}}
\def\ang#1{\langle #1 \rangle}
\def\MMb{\overline{M}^2_b}

\newcommand{\pa}{\partial}
\newcommand{\Mbar}{\overline{M}}
\newcommand{\abs}[1]{{\left\lvert{#1}\right\rvert}}

\newcommand{\ep}{{\epsilon}}
\newcommand{\CI}{{\mathcal{C}^\infty}}

\newcommand\lb{\operatorname{lb}}
\newcommand\rb{\operatorname{rb}}

\newcommand\conic{{\operatorname{conic}}}

\DeclareMathOperator{\Op}{{Op}}
\DeclareMathOperator{\supp}{{supp}}
\DeclareMathOperator{\graph}{{graph}}
\DeclareMathOperator{\WF}{{WF}}
\DeclareMathOperator{\RC}{{RC}}

\theoremstyle{plain} \newtheorem{theorem}[subsection]{Theorem}
  
  \newtheorem{proposition}[subsection]{Proposition}
  \newtheorem{lemma}[subsection]{Lemma}
  \newtheorem{corollary}[subsection]{Corollary}

\theoremstyle{remark} 
  \newtheorem{remark}[subsection]{Remark}
  \newtheorem{remarks}[subsection]{Remarks}
  \newtheorem{example}[subsection]{Example}

\theoremstyle{definition} \newtheorem{definition}[subsection]{Definition}
  
\numberwithin{equation}{section}

\include{psfig}

\begin{document}

\title[Strichartz on manifolds]{Sharp Strichartz estimates on non-trapping asymptotically conic manifolds}

\author{Andrew Hassell} \thanks{A.H.\ is supported in part by an Australian
Research Council Fellowship.}  \address{Department of Mathematics, ANU,
Canberra, ACT 0200, AUSTRALIA} \email{hassell@maths.anu.edu.au}

\author{Terence Tao} \thanks{T.T.\ is a Clay Prize Fellow and is supported
in part by a grant from the Packard Foundation.}  \address{Department of
  Mathematics, UCLA, Los Angeles California 90095, USA} \email{tao@math.ucla.edu}

\author{Jared Wunsch} \thanks{J.W.\ is supported in part by NSF grants
DMS-0323021 and DMS-0401323.}  \address{Department of Mathematics, Northwestern University,
  Evanston IL 60208, USA}
\email{jwunsch@math.northwestern.edu}



\vspace{-0.3in}
\begin{abstract}
We obtain the Strichartz inequalities
$$ \| u \|_{L^q_t L^r_x([0,1] \times M)} \leq C \| u(0) \|_{L^2(M)}$$
for any smooth $n$-dimensional Riemannian manifold $M$ which is
asymptotically conic at infinity (with either short-range or long-range metric perturbation)
and non-trapping, where $u$ is a solution
to the Schr\"odinger equation $iu_t + \frac{1}{2} \Delta_M u = 0$, and $2 < q, r \leq \infty$
are admissible Strichartz exponents ($\frac{2}{q} + \frac{n}{r} = \frac{n}{2}$).  
This corresponds with the estimates available for Euclidean space (except for the endpoint 
$(q,r) = (2, \frac{2n}{n-2})$ when $n > 2$). These estimates
imply existence theorems for semi-linear Schr\"odinger equations on $M$, by adapting arguments from Cazenave and Weissler \cite{cwI} and Kato \cite{kato}. 

This result improves
on our previous result in \cite{HTW}, which was an $L^4_{t,x}$ Strichartz estimate in three dimensions. It is closely related to
the results in \cite{st}, \cite{burq}, \cite{tataru}, \cite{rz}, which
consider the case of asymptotically flat manifolds.
\end{abstract}

\maketitle

\tableofcontents

\section{Introduction}\label{sec:introduction}

The purpose of this paper is to establish the full range (except for the endpoint)
of (local-in-time) Strichartz inequalities on a class of non-Euclidean spaces, namely smooth
asymptotically conic Riemannian manifolds
$(M,g)$ which obey a non-trapping condition.  As in our earlier paper \cite{HTW}, by an \emph{asymptotically conic manifold} we mean a smooth complete noncompact Riemannian manifold $M = M^n$ of dimension $\dim M = n \geq 2$, with metric\footnote{We use $z \in \RR^n$ as the spatial co-ordinate on $M$, in order to reserve the letter $x$ for the scattering co-ordinate $x := 1/r$.  We neglect the case $n=1$ as such manifolds are isometrically equivalent to a finite number of copies of 
the Euclidean real line $\R$.} 
$g:=g_{jk}(z) dz^j dz^k$ such that there is a compact set $K_0 \subset M$ and an $(n-1)$-dimensional smooth compact
Riemannian manifold $(\partial M, h)$, with metric $h := h_{jk}(y) dy^j dy^k$ such that the \emph{scattering region} or \emph{asymptotic region} $M \backslash
K_0$ can be parametrized as the collar neighbourhood
$$ M \backslash K_0 \equiv (0, \epsilon_0) \times \partial M := \{ (x, y):
0 < x < \epsilon_0, y \in \partial M \}$$ 
for some $\epsilon_0 > 0$, with
metric which can be brought to the form
\begin{equation} g = \frac{dx^2}{x^4} + \frac{h_{jk}(x,y) dy^jdy^k}{x^2} = dr^2 + r^2 h_{jk}(\frac{1}{r}, y) dy^j dy^k.
\label{ac-metric}
\end{equation}
Here $r := \frac{1}{x} \in (\frac{1}{\epsilon_0},+\infty)$ is the ``radial'' variable, $y^j$ are the $n-1$ ``angular''variables, the ``scattering co-ordinate'' $x = \frac{1}{r}$ is the defining function for $\partial M$, and $h_{jk}(x,y)$ is a smooth function on $[0,\epsilon_0) \times \partial M$ such that $h_{jk}(0,y) =h_{jk}(y)$. If $h_{ij}$ is independent of $x$ (so that $h_{ij}(x,y)= h_{ij}(y)$) we say that $M$ is \emph{perfectly conic near infinity} and write $g_\conic$ instead of $g$. It is then natural to compactify $M$ to $\overline M := M \cup \partial M$ by identifying $\partial M$ with $\{0\} \times \partial M$ in this co-ordinate chart.  In the terminology of Melrose \cite{melrose}, the metric $g$ on $\Mbar$ is thus a \emph{short-range scattering metric}. 
As is customary, we say that $M$ is \emph{non-trapping} if every geodesic $\gamma: \R \to M$ reaches $\partial M$ at $\pm \infty$.

We shall also consider long range metrics in this paper. We say that $g$ is a \emph{long range scattering metric} on $M$ if it can be brought to the form \eqref{ac-metric} in the scattering region, but where the function $h_{jk}(x,y)$ is no longer smoothup to the boundary $x=0$, but instead obeys the expansion\begin{equation}\label{ac-lr-metric}h_{jk}(x,y) = h_{jk}(y) + x^\delta e_{jk}(x,y)\end{equation}for some $0 < \delta \leq 1$, where $e_{jk}$ are symbols of order 0 (see Definition \ref{symbol-def}).

\begin{remark}\label{general-lr-metrics} By a result of Joshi and S\'a Barreto, every short range scattering metric can be brought to the form \eqref{ac-metric}. In section~\ref{normal-form} we prove a generalization: any metric of the form 
\begin{equation}\label{ac-lr-metric2}
\frac{dx^2}{x^4} + \frac{h(y, dy)}{x^2} + x^\delta \Big( k_{00} \frac{dx^2}{x^4} + k_{0j} \frac{dx}{x^2} \frac{dy^j}{x} + k_{ij} \frac{dy^i}{x} \frac{dy^j}{x} \Big)\end{equation}
where the $k_{ij}$ are symbols of order zero, can be brought to the form \eqref{ac-metric}, with $h_{jk}$ as in \eqref{ac-lr-metric}. If $M = \RR^n$ and $h$ is the standard metric on $S^{n-1}$, then \eqref{ac-lr-metric2} may be expressed as follows: the metric is of the form
$$
g_{ij} = \delta_{ij} + g'_{ij},
$$
where the $g'_{ij}$ are symbols of order $-\delta$. Thus our class of metrics include the
long range metrics considered by Burq in \cite{burq}, and is also somewhat more general than the definition given in \cite{Vasy}. 
\end{remark}

It will be convenient to select a somewhat artificial, but globally defined radial positive weight $\langle z \rangle$ on $M$, chosen so that $\langle z \rangle$ is equal to $r$ in the scattering region $(0,\epsilon_0) \times \partial M$,
is comparable to $1/\epsilon_0$ in the compact interior region $K_0$, and is smooth throughout.

\begin{example}
The most important example  of an asympotically conic manifold is Euclidean
space $M = \R^n$ itself, which is perfectly conic and nontrapping (with 
$\partial M$ being the unit sphere $S^{n-1}$ with the standard metric, and $(r,y)$ being
the usual polar co-ordinates). More generally, any 
compact perturbation of Euclidean space is also asymptotically conic (and perfectly conic at infinity), although 
it may or may not be nontrapping. Any compactly supported perturbation of Euclidean space which is sufficiently small in $C^2$ will be nontrapping. 
\end{example}

\

Let $(M, g)$ be an asymptotically conic manifold. We let $dg(z) := \sqrt{g}\ dz$ be the measure induced by the 
metric $g$, and let $L^2(M)$ be the complex Hilbert space given by the inner product
$$ \langle f_1, f_2 \rangle_{L^2(M)} := \int_M f_1(z) \overline{f_2(z)}\
dg(z).$$ 
Let $H = \frac1{2} \nabla^* \nabla$ be minus one-half the Laplace-Beltrami operator on $M$. As is well known,
$H$ is self-adjoint and positive-definite on $L^2(M)$ and has a functional calculus on its
spectrum $[0,\infty)$, and one can define Sobolev spaces $H^s(M)$ for all
$s \in \R$ by taking the closure of test functions under the norm $\| u
\|_{H^s(M)} := \| (1 + H)^{s/2} u \|_{L^2(M)}$.

We now consider Schwartz\footnote{This is in order to avoid needless technicalities, but our estimates
will not depend on any of the Schwartz semi-norms of the solution and so can be extended to rougher solutions.}
 solutions $u: \R \times M \to \C$ to the (time-dependent) Schr\"odinger equation
\begin{equation}\label{nls}
 u_t = -iH u;
\end{equation}
this is the natural quantization of the geodesic flow equation.
It is well known (see e.g.\ \cite{cks}) that if $u(0)$ is Schwartz then there is a unique global Schwartz solution to 
\eqref{nls}.  

For fixed $n$, let us call a pair of exponents $(q,r)$ \emph{admissible} if we have the conditions
\begin{equation}\label{cond}
 \frac{2}{q} + \frac{n}{r} = \frac{n}{2}; \quad 2 \leq q,r \leq \infty; \quad (q,r,n) \neq (2,\infty,2).
\end{equation}
In Euclidean space $(M, g) = (\R^n, \delta)$ we have the global-in-time \emph{Strichartz estimate} 
(see e.g.\ \cite{tao:keel} and the references therein)
\begin{equation}\label{strich}
\| u \|_{L^q(\R; L^r(\R^n))} \leq C_{n,q,r} 
\| u(0) \|_{L^2(\R^n)}
\end{equation}
for all admissible exponents $(q,r)$.  The admissibility conditions are best
possible (see \cite{tao:keel} for a discussion, and \cite{montgomery},
\cite{tao} for the failure of the estimates at the endpoint $(q,r,n) =
(2,\infty,2)$).  One can also generalize these estimates to other Sobolev
spaces via Sobolev embedding (exploiting the fact that the Schr\"odinger
flow commutes with fractional differentiation operators $(1+H)^{s/2}$).
The purpose of this paper is to show that these Strichartz estimates extend to
more general non-trapping manifolds, at least locally in time and excluding
the endpoint $q = 2$.

\begin{theorem}[Strichartz estimate]\label{main}  For any $n$-dimensional asymptotically conic non-trapping manifold $M$, 
with either short range or long range metric $g$, 
and any (Schwartz) solution $u$ to \eqref{nls}, we have
\begin{equation}\label{strichartz}
\| u \|_{L^q([0,1]; L^r(M))} \leq C_{n,q,r,M} \| u(0) \|_{L^2(M)} 
\end{equation}
for any admissible pair $(q,r)$ with $q > 2$.
\end{theorem}

\begin{remark} We have recently learnt that similar estimates have been independently announced by Tataru \cite{tataru} and by Robbiano-Zuily \cite{rz}. These authors treat asymptotically Euclidean, rather than asymptotically conic metrics, but in other respects their results may be more general than ours (for example, Tataru treats $C^2$ metrics). 
\end{remark}

\begin{remark}There are some standard extensions of this estimate, as well as
applications to low regularity local existence theory of semilinear wave
equations, in the spirit of \cite{cwI}; we discuss these in Section
\ref{remarks:sec}.  The extension to long range metrics is routine but somewhat technical, requiring a re-proof
of a key lemma (Lemma \ref{metric-smooth}) which we defer to Appendix \S \ref{longrange}.
\end{remark}

We contrast Theorem \ref{main} with previously known results.  In the rest
of the paper we allow all our constants $C$ to depend on $n, q, r, M$ and
will not mention this dependence explicitly again.  Without the
non-trapping condition, one does not expect to obtain the estimate
\eqref{strich} without conceding some derivatives, because of
quasimode-type solutions to \eqref{nls} (possibly allowing for a small
forcing term on the right-hand side).  However, it is still possible to
obtain a (local in time) Strichartz estimate with a loss of derivatives.
For instance, for admissible $(q,r)$, Burq, Gerard and Tzvetkov \cite{bgt}
showed that the estimate
\begin{equation}\label{strich-manifold}
\| u \|_{L^q([0,1]; L^r(M))} \leq C \| u(0)
 \|_{H^{\frac{1}{q}}(M)}
\end{equation}
holds for arbitrary smooth manifolds (trapping or non-trapping) in general
dimension; a key tool here is the use of parametrices for the propagator
$e^{-itH}$ for frequency-localized solutions and for extremely short time
intervals.  This should be compared to what one can obtain purely from
Sobolev embedding, which loses $2/q$ derivatives.  This loss of $1/q$
derivatives, as well as the localization in time, was shown in \cite{bgt}
to be sharp in the case of the sphere (which is in some sense maximally
trapping).

In the non-trapping case one expects better estimates, as local smoothing estimates become available.  For our
purposes we need a somewhat refined ``half-angular local smoothing estimate'' introduced in \cite{HTW} (see also \cite{Sug} for a Euclidean space variant).  We need a rather technical Banach space $X$, defined formally in
Section \ref{sec:doi}, but which is heuristically of the form
$$ \| f \|_X \approx \| f \|_{L^2(M)} + \| \langle z \rangle^{-1/2-\eps} |\nabla|^{1/2} f \|_{L^2(M)}
+ \| \langle z \rangle^{-1/2} |\nabb|^{1/2} f \|_{L^2(M)}$$
where $\nabb := x \nabla_y = \frac{1}{r} \nabla_y$ is the angular portion of the derivative\footnote{We have to insert a smooth cutoff to the scattering region $r \gg 1$ in order to ensure that $\nabb$ is well-defined and has smooth coefficients, and to justify taking the square root of $|\nabb|$ we will need the scattering pseudo-differential calculus.  We will do all this rigorously in Section \ref{sec:doi}.}.  

\begin{lemma}[Half-angular local smoothing effect \cite{HTW}]\label{ls-lemma}  Let $M$ be a smooth non-trapping 
asymptotically conic manifold,
and let $u$ be a Schwartz solution to \eqref{nls}.  Then we have
\begin{equation}\label{ls-scat}
\|u\|_{L^2([0,1];X)} := \left(\int_0^1 \| u(t) \|_X^2\ dt\right)^{1/2} \leq C \| u(0) \|_{L^2(M)}.
\end{equation}
\end{lemma}

\begin{remarks}
This Lemma shall follow immediately from \cite[Lemmas 10.4-10.5]{HTW} once
we define the space $X$ properly in Section \ref{sec:doi}.  The point is
that if one wants to gain a half-derivative smoothing effect in all
directions then one must lose an additional factor of $\langle z
\rangle^{-\eps}$, but if one only wishes for a half-derivative smoothing
effect in just the \emph{angular} directions then one does not lose this
$\eps$ of decay.  This is crucial in order not to lose even an epsilon of
derivatives in \eqref{strichartz}.  The well known Morawetz estimate (based
on commuting the flow with the vector field $\partial_r$) allows one to
obtain a slightly weaker ``full-angular local smoothing effect'' in which
involves one full angular derivative and $-1/2$ of a general derivative,
but this estimate turns out to be too weak for our arguments here.  See
\cite{HTW} for some further comparisons among these various local
smoothing estimates.  For our argument we shall also need to prove
\eqref{ls-scat} not just for the actual solution $u = e^{-itH} u(0)$ of the
Schr\"odinger equation, but also for various localized parametrices to this
solution, but it turns out that this will also follow by a modification of
the positive commutator method used in \cite{HTW}.
\end{remarks}

The non-trapping hypothesis was exploited by Staffilani and Tataru
\cite{st}, who were able to remove the $1/q$ loss of derivatives in
\eqref{strich-manifold} completely (i.e.\ to prove \eqref{strichartz}) for
metrics on $\R^n$ which are non-trapping and Euclidean outside a compact
set (requiring only $C^{1,1}$ regularity of the metric); a key tool was the
use of local smoothing estimates such as \eqref{ls-scat} to control the
error term arising from localization to a compact region of space.  They
were also able to attain the $q=2$ endpoint with some additional arguments.
After the work of Staffilani and Tataru, Burq \cite{burq} gave an
alternative proof of the same result (using the local parametrices from
\cite{bgt}), and also considered metrics which were asymptotically flat
rather than flat outside of a compact set.  In this more general setting
Burq was able to recover almost the same Strichartz estimate as Staffilani
and Tataru (for $q>2$ at least), but with an epsilon loss of derivatives:
$$ \| u \|_{L^q_t L^r_x([0,1] \times M)} \leq C \| u(0)
 \|_{H^{\eps}(M)}.
 $$ 
 The idea was to divide $M$ into dyadic shells $R
 \leq |z| \leq 2R$ and apply a variant of the Staffilani-Tataru argument on
 each shell, relying on the non-angular component of the local smoothing estimate
\eqref{ls-scat} to control the error terms.  The presence of the $\eps$ loss of decay in the non-angular
local smoothing estimate eventually contributes to the epsilon
loss of derivatives in Burq's result\footnote{We have learnt that this epsilon loss has very recently been removed
by independent work of Tataru \cite{tataru} and Robbiano-Zuily \cite{rz}.}.

In \cite{HTW}, we were able to remove this $\eps$ loss, but only in the
special context of the specific Strichartz estimate
$$ \int_0^1 \int_M |u(t,z)|^4\ dg(z) dt \leq C_M \| u(0)
\|_{H^{1/4}(M)}^4$$ in the three-dimensional case $n=3$. The technique here
relied on an interaction Morawetz (or positive commutator) method rather
than a parametrix-based approach, and as such was restricted to the
$L^4_{t,x}$ set of exponents.  However, the proof relied geometrically on
some refined estimates on the metric function $d_M(z,z')$---in particular,
that it is smooth and obeys product symbol estimates on thin cones---and we
will use these estimates again to prove Theorem \ref{main}.  Another key
feature of the argument in \cite{HTW} was a crucial reliance on the
half-angular component of the local smoothing estimate \eqref{ls-scat}.
This estimate avoids the epsilon loss problem, but on the other hand only
controls angular derivatives.  This forced us in \cite{HTW} to rely
primarily on cutoff functions which could vary in the angular variable $y$
but not in the radial variable $x$, thus localizing our solutions to thin
cones rather than dyadic annuli.  For similar
reasons\footnote{Semi-classically, the reason why the localization to
dyadic annuli costs an epsilon of derivatives while the localization to
cones do not is that each geodesic in $M$ can (and will) pass through an
infinite number of dyadic annuli, but can only cross a finite number of
cones (and cannot stay in the compact region indefinitely, by the
non-trapping hypothesis).} we shall adopt a similar localization to cones
in this paper.

Our argument here extends to more general exponents than the $L^4$ estimate in \cite{HTW} (but not the endpoint
estimates).  We sketch the proof as follows.  We cover the compact manifold $\overline{M}$ by a finite number of small open sets; back on the original manifold $M$, this corresponds to a covering by a finite number of small balls and a finite number of ``thin cones.''  We then apply cutoff functions to localize to each of these sets, and then construct an approximate ``local''
parametrix for the propagator $e^{-itH}$.  Indeed, our parametrix is rather simple, and is based on the parametrix $U(t)$ from \cite{Hassell-Wunsch1}, whose kernel $U(t,z,z')$ had the form 
$$
U(t, z, z') := (2\pi i t)^{-n/2} e^{id_M(z,z')^2/2t} a(z, z', t)
$$
when $(z,z')$ are sufficiently close, so that there is a unique minimizing geodesic between them. For flat Euclidean space, of course, the kernel has this form globally, with $a \equiv 1$; in \cite{Hassell-Wunsch1} $a$ was smooth in $t$ with the coefficients of its Taylor series  at $t=0$ given by a solution of certain transport equations. Here,
we simply take this as an ansatz, restrict  $a(z, z, t)$ to $t=0$, and insert a cutoff function which is supported in a region, containing the diagonal, where $z'$ and $z$ are joined by unique minimizing geodesics. More precisely, and to fix some notation, we take our local parametrix $U(t)$ for the propagator $e^{-itH}$ to be
\begin{equation}\label{parametrix}
 U(t) f(z) := \frac{1}{(2\pi i t)^{n/2}} \int_{M} e^{i \Phi(z,z')/t} \chi(z,z') a(z,z') f(z')\ dg(z')
\end{equation}
where $\Phi$ is the phase function
\begin{equation}\label{Phi-def}
\Phi(z,z') := d_M(z,z')^2/2,
\end{equation}
and $d_M$ is the metric function on $M$, and $a(z,z') dg(z)^{1/2} dg(z')^{1/2}$ is the half-density defined (for $z$ and $z'$ sufficiently close) 
as\footnote{Note that $-\nabla_z \nabla_{z'} \Phi$ is most naturally a section of the tensor bundle $T_z M \otimes T_{z'} M$, and so its determinant naturally is a section of the tensor bundle $(\bigwedge^n T_z M) \otimes (\bigwedge^n T_{z'} M)$, and is thus a scalar multiple of $dg(z) dg(z')$.  Thus the equation \eqref{a-def} is tensorial, or equivalently is invariant under changes of co-ordinates.  Indeed, the choices of $\Phi$ and $a$ are not arbitrary, but are dictated by the geodesic flow on phase space.}
\begin{equation}\label{a-def}
 a(z,z') dg(z)^{1/2} dg(z')^{1/2} := \det( -\nabla_{z} \nabla_{z'} \Phi(z,z') )^{1/2},
\end{equation}
and $\chi$ is a  sum of cutoff functions, each of which localizes $z$, $z'$ to either a small ball or a thin cone. 
When $z$ and $z'$ are widely separated, this
parametrix cannot be expected to work since the metric function can develop
singularities\footnote{Microlocally, one can think of this parametrix as an
approximation to the initial portion of the flow, before the geodesic flow
develops bad behaviour such as caustics or conjugate points.  One can
obtain a more accurate parametrix for the full evolution by composing this
local parametrix with itself a finite number of times using the Duhamel
formula, but we will not do so in this paper.}.  However, when $z, z'$ are restricted to a small ball or to a thin cone then this
parametrix is reasonably accurate,
essentially because $\Phi$ and $a$
solve the eikonal and Hamilton-Jacobi equations associated to the Schr\"odinger flow \eqref{nls}.  
The parametrix here has two key advantages over that in \cite{Hassell-Wunsch1}. First it is specified when both variables $z$ and $z'$ approach infinity (provided the angular separation between them remains small), while in \cite{Hassell-Wunsch1} the first variable $z$ was restricted to lie in a compact set in $M$, although it was completely global in the other variable $z'$. Second, it satisfies the  $O(|t|^{-n/2})$ $L^1 \to L^\infty$ estimate, while the parametrix of \cite{Hassell-Wunsch1}  in general does not satisfy such an estimate.  In fact, it follows from \cite{Hassell-Wunsch1} that \emph{the propagator itself fails to satisfy such an $L^1 \to L^\infty$ estimate} at any point $(z, z')$ such that there is a geodesic emanating from $z$ with a conjugate point at $z'$, so the cutoff in \eqref{parametrix} is essential for the method used here. 

The fact that the parametrix 
obeys $L^1(M) \to
L^\infty(M)$ type dispersive estimates with an $O(|t|^{-n/2})$ bound allows us\footnote{More precisely, we
also need dispersive inequalities on the $TT^*$ version of this parametrix,
but this can be obtained by stationary phase arguments.  Also we need to
prove energy and local smoothing estimates for the parametrix.}  to apply
the abstract Strichartz estimate in \cite{tao:keel} to obtain the required
estimates for each of the localized parametrices. 
The proof of the Strichartz estimates for the actual propagators $e^{-itH}$
then hinges on dealing with the error
terms caused by the inaccuracies of the parametrices and by the cutoff
functions.  For the terms localized to small balls these terms can be
controlled by the standard local smoothing estimates (as in the arguments
of Staffilani-Tataru \cite{st} and Burq \cite{burq}).  For the terms
localized to thin cones one uses half-angular local smoothing estimates
such as those in Lemma \ref{ls-lemma} to control the error terms; a key
point here is that because the cutoff depends mainly on the angular
variable $y$ and not on the radial variable $x$, the terms arising from
differentiating the cutoff will mainly involve \emph{angular} derivatives
and are thus compatible with the half-angular $X$ spaces in Lemma
\ref{ls-lemma} (a precise formulation of this is in Lemma \ref{X-split}).
The non-endpoint hypothesis $q \neq 2$ is exploited (as in \cite{st},
\cite{burq}) via the Christ-Kiselev lemma \cite{ck}, which allows one to
decouple the Strichartz and the local smoothing aspects of our estimates.

As in \cite{HTW}, a key component of the argument is a detailed analysis of
the metric function $d_M(z,z')$ and its derivatives for various values of
$z, z'$, including cases in which $z$ and $z'$ both lie in a thin cone but where $\ang{z}/\ang{z'}$ 
or $\ang{z'}/\ang{z}$ 
could be arbitrarily large. However, the geometric information we need
on the metric is slightly different from that in \cite{HTW}.  There, the
main objective was to obtain a certain geodesic (weak) convexity of the
metric function $d_M(z,z')$, modulo lower order errors.  Here, by
constrast, it is more important to obtain non-degeneracy estimates on the
various Hessians of this metric function (in order to apply stationary
phase).

This paper is organized as follows.  In Section \ref{G-sec} we state the
necessary geometric estimates we need on asymptotically conic manifolds;
the proof of these estimates are technical and will be deferred to Appendix
\S \ref{sec:geometry}, with one particular long-range estimate deferred to
Appendix \S \ref{longrange}.  Then in Section \ref{sec:LSIO} we introduce
the concept of a \emph{local Schr\"odinger integral operator} (LSIO), which
are basically Fourier integral operators associated to the Schr\"odinger
flow, which are localized with respect to the compactified manifold
$\overline{M}$ but not the original manifold $M$.  We give the basic theory
for these operators in that section, but defer the proofs (which are mostly
standard but technical applications of the $TT^*$ method and the principle
of stationary phase) to Appendix \S \ref{FIO-sec}.  Then, in Section
\ref{sec:doi} we review from \cite{HTW} and \cite{melrose} the calculus of
scattering pseudo-differential operators, which we need to rigourously
define the space $X$ alluded to earlier, and to establish local smoothing
estimates for both the Schr\"odinger propagator and the localized
parametrix (which we construct in Section \ref{parametrix-sec}).  Finally,
in Section \ref{sec:proof} we assemble all these ingredients together,
using the Duhamel formula and the Christ-Kiselev lemma (which we review in
Appendix \S \ref{sec:lemmas}), to quickly prove Theorem \ref{main}.  In
Section \ref{remarks:sec} we discuss some extensions and applications of
Theorem \ref{main}.

We thank Ben Andrews and Nicolas Burq for illuminating conversations.


\section{The geometry of asymptotically conic manifolds}\label{G-sec}

In this section we set out our notation and describe some of the geometry
of asymptotically conic manifolds.  While the statements here are
intuitively plausible, some of the proofs are a little technical and will
be deferred to an Appendix \S \ref{sec:geometry}.  For this discussion we shall consider
both short-range and long-range metrics.

The Riemannian metric $g$ induces a distance function $d_M$ on $M$, and
similarly $h(0,\cdot)$ induces a distance function $d_{\partial M}$ on $\partial M$.
It will be convenient to also select a smooth
Riemannian metric $\overline{g}$ on compactified manifold
$\overline{M}$; the exact choice of this metric is not particularly important since
$\overline{M}$ is compact (and hence all Riemannian metrics are
equivalent), however for sake of concreteness we require $\overline{g}$ to be given
in the asymptotic region $0 \leq x < \epsilon_0$ by the formula\footnote{One can interpret this metric as
the natural metric obtained after using a ``radial compactification'' to identify $M$ with a ``hemisphere'' with
boundary $\partial M$; see the proof of Lemma \ref{glc} in Appendix \S \ref{sec:geometry}.}
\begin{equation}\label{ovg}
 \overline{g} = \frac{dx^2}{(1 + x^2)^2} + \frac{h_{ij}(x,y) y^i y^j}{1+x^2}
\end{equation}
(compare with \eqref{ac-metric}); the choice of metric for $\overline{g}$ in the near region $K_0$ is not relevant
for our argument so long as it is smooth.  The advantage of using \eqref{ovg} for the compactified metric $\overline{g}$ is that
in the perfectly conic case $h(x,y) = h(0,y)$, the geodesics corresponding to $g$ and to
$\overline{g}$ (near infinity at least); see Lemma \ref{glc}.  

Let $d_{\overline{M}}$ be the distance function on $\overline{M}$
associated to $\overline{g}$.  Observe that $d_{\overline{M}}(z,z')$ will
be comparable to $d_M(z,z')$ when $z, z'$ are restricted to a fixed compact
set (with the comparability constants depending of course on this set), and
when $z = (x,y), z' = (x',y')$ are in the asymptotic region then
$d_{\overline{M}}(z,z')$ will be comparable to $|x - x'| + d_{\partial
M}(y, y')$.  Thus in particular, small balls in the $d_{\overline{M}}$
metric will either correspond to small balls in the near region of $M$, or
thin cones in the asymptotic region of $M$.  For any $\epsilon > 0$, we
define the neighbourhoods of the diagonal $\Delta_\epsilon := \{ (z,z') \in
M^2: d_{\overline M}(z,z') \leq \epsilon \}$, together with their
compactifications $\overline{\Delta_\epsilon} := \{ (z,z') \in
\overline{M}^2: d_{\overline M}(z,z') \leq \epsilon \}$.

The basic philosophy of this section is that for sufficiently small
$\epsilon > 0$, the geodesic structure of the original metric $d_M$ (as
well as the structure of the closely related functions $\Phi$ and $a$) on
$\Delta_\epsilon$ behaves very much like its Euclidean counterpart.  The
crucial point (which was also key in \cite{HTW}) is that we will still
control all the geometry in the case when the radial variables $r$, $r'$
are widely separated, so long as we have the angular closeness condition
$d_{\partial M}(y,y') \ll 1$; this is essential in order for us not to lose
an epsilon of derivatives or decay in our final Strichartz estimates,
essentially because the Schr\"odinger flow tends to propagate along
geodesics, and (thanks to the non-trapping condition) any given geodesic
eventually ends up on a thin cone where the angular variable $y$ is in a
small set, but the radial variable $r$ goes to infinity.  When
$d_{\overline M}(z,z')$ is large, we do not attempt to control the geometry
at all, and our parametrix will be highly inaccurate in this region;
however we can still recover estimates for the genuine propagator
$e^{-itH}$ from that of the parametrix by means of Duhamel's formula and
local smoothing estimates, which are genuinely global estimates relying on
the non-trapping property of the manifold.

We treat $\epsilon > 0$ as a small parameter, and adopt the following convenient notation: we use
$o(1)$ to denote any (possibly tensor-valued) quantity which is bounded in magnitude by $c(\epsilon)$ for some quantity $c(\epsilon)$
depending only on $\epsilon$ (and on the manifold $M$) which goes to zero as $\epsilon \to 0$.  This quantity $c(\epsilon)$ will vary from line to line.  More generally we abbreviate 
$o(1) X$ as $o(X)$ for any non-negative quantity $X$.  In
what follows we shall always assume that $\epsilon$ is sufficiently small (depending only on $M$).
We also adopt the convention that if a statement involved a point $z \in M$ in the asymptotic region,
then $r = r(z), y = y(z), x = x(z)$ are assumed to be the radial, angular, and scattering co-ordinates of $z$
respectively; similarly for $z'$ and $r', y', x'$, etc.

We begin with a basic lemma on the geodesic structure of $M$ near the diagonal.

\begin{lemma}[Geodesic local connectedness of $\overline{M}$]\label{glc}
If $z, z' \in M$ is such that $d_{\overline M}(z,z') = o(1)$, there
is a unique length-minimizing geodesic\footnote{All our geodesics will be parameterized in the usual constant-speed (or energy-minimizing) manner, i.e. they will be projections to physical space of the geodesic flow on the cotangent bundle $T^* M$.} $\gamma_{z \to z'}: [0,1] \to M$ which goes from $z$ to $z'$, and furthermore we have 
$\diam_{\overline{M}}(\gamma_{z \to z'}) = o(1)$ (i.e.\ $\gamma_{z \to z'}$ lies inside a ball of radius $o(1)$ in 
the $d_{\overline M}$ metric).  In particular, if $z$ and $z'$ are in the asymptotic region
$x, x' = o(1)$, then so is the geodesic $\gamma_{z \to z'}$.  Furthermore, in this case we have the additional estimate
\begin{equation}\label{cycle}
r(\gamma_{z \to z'}(\theta)) = (1 + o(1)) ( (1-\theta) r + \theta r' )
\end{equation}
for all $0 \leq \theta \leq 1$, where $r(\gamma_{z \to z'}(\theta))$ is the radial co-ordinate of
$\gamma_{z \to z'}(\theta)$.  If $g = g_\conic$ is a perfectly conic metric in this asymptotic region, 
then we also have the cosine rule
\begin{equation}\label{sharp-cosine}
2\Phi_{\conic}(z,z') = d_\conic(z,z')^2 =  r^2 + (r')^2 - 2r r' \cos d_{\partial M}(y,y')
\end{equation}
and furthermore the geodesics associated to $g_\conic$ and to the compactified metric $\overline{g}$ coincide
in this region.
\end{lemma}

\begin{proof}  For short-range metrics, most of these claims are essentially 
in \cite[Lemma 9.1, Proposition 9.2]{HTW}, but we shall give a complete proof (covering the long-range case also)
in Appendix \S \ref{sec:geometry}.  
\end{proof}

\begin{remark} From this lemma we can uniquely define the geodesic $\gamma_{z \to z'}$ for any pair 
of points $z, z'$ which are sufficiently close in the compactified metric $d_{\overline M}$.  One can think of the map $z \mapsto \gamma'_{z_0 \to z}(0)$ as normal co-ordinates around $z_0$, thus Lemma \ref{glc} implies that these normal co-ordinates are smooth as long as $d_{\overline M}(z_0,z)$ is sufficiently small.
\end{remark}

\begin{example}
In Euclidean space $M=\R^n$, we can define $\gamma_{z \to z'}$ for all pairs $z, z'$, indeed we have $\gamma_{z \to z'}(\theta) = (1-\theta) z + \theta z'$.  Also we have $\Phi(z,z') = |z-z'|^2/2$ and $a(z,z') = 1$ in this case.  
\end{example}

Now we study the phase function $\Phi$ defined in \eqref{Phi-def}.  This function is clearly real-valued and symmetric.
From normal co-ordinates around $z$ and Lemma \ref{glc} we see that $\Phi$ is smooth on $\Delta_\epsilon$ for sufficiently small $\epsilon$.  Also, we recall \emph{Gauss's Lemma}, which in our notation asserts that
\begin{equation}\label{gauss-lemma}
\nabla^{z'} \Phi(z,z') := \gamma'_{z \to z'}(1) = - \gamma'_{z' \to z}(0);
\quad \nabla^z \Phi(z,z') := -\gamma'_{z \to z'}(0) = \gamma'_{z \to z'}(1),
\end{equation}
valid when $d_{\overline M}(z,z')$ is sufficiently small.  Also, since geodesics have constant speed, we have
\begin{equation}\label{constant-speed}
|\gamma'_{z \to z'}(\theta)|_{g(\gamma_{z \to z'}(\theta))} = d_M(z,z')
\end{equation}
for all $0 \leq \theta \leq 1$ and all $z, z'$ sufficiently close in $d_{\overline M}$.  
Combining these two equations we obtain the \emph{eikonal equation}
\begin{equation}\label{eikonal}
\Phi(z,z') = \frac{1}{2} \big| \nabla^{z} \Phi(z,z') \big|_{g(z)}^2 = \frac{1}{2} \big| \nabla^{z'} 
\Phi(z,z') \big|_{g(z')}^2
\end{equation}
again valid when $d_{\overline M}(z, z')$ is sufficiently small.  The first equation in \eqref{eikonal}
can be rewritten as
\begin{equation}\label{eikonal-variant}
(i t^2 \partial_t - t^2 H_z) e^{i \Phi(z,z') / t} = O_{z,z'}(t)
\end{equation}
where we use $O_{z,z'}(t)$ to denote a quantity which vanishes to first order at $t=0$ for any fixed $z,z'$,
and $H_z = \frac{1}{2} \nabla_z^* \nabla_z$ is the operator $H$ applied to the $z$ variable.
From \eqref{eikonal-variant} we naturally expect the phase $e^{i\Phi(z,z')/t}$ to arise any parametrix for the Schr\"odinger operators $e^{-itH}$, at least in the region where $z,z'$ are geometrically close; much of this 
paper is devoted to the justification of such a heuristic.

Another consequence of \eqref{gauss-lemma} is the identities
\begin{align*}
\nabla^{z'} \Phi(\gamma_{z_0 \to z_1}(\theta), \gamma_{z_0
\to z_1}(\theta')) &= (\theta'-\theta) \gamma'_{z_0 \to z_1}(\theta')\\
\nabla^{z} \Phi(\gamma_{z_0 \to z_1}(\theta), \gamma_{z_0
\to z_1}(\theta')) &= -(\theta'-\theta) \gamma'_{z_0 \to z_1}(\theta)
\end{align*}
for any $0 \leq \theta, \theta' \leq 1$ and $z_0, z_1$ sufficiently close in $d_{\overline M}$.  In particular we have
\begin{equation}\label{id}
 \nabla_{z''} \big( (1-\theta)\Phi(z'',z) + \theta \Phi(z'',z') \big) |_{z'' = \gamma_{z \to z'}(\theta)} = 0,
\end{equation}
whenever $d_{\overline M}(z,z')$ is sufficiently small and $0 \leq \theta\leq 1$.  In fact this identity can be used 
to define $\gamma_{z \to z'}(\theta)$, locally at least; see \eqref{gradw} below.

Now we show that $\Phi$ is well-behaved whenever the compactified distance $d_{\overline{M}}(z,z')$ is 
small. First we give

\begin{definition}\label{symbol-def}  
We say that a function $b: M \to \C$ is a \emph{symbol of order $j$ on $M$} if we have the bounds
$$
\big| \nabla_y^\alpha \nabla_r^j b \big|  \leq C_{\alpha} \langle z \rangle^j
$$
for all $\alpha, j \geq 0$ and $z \in M$. This is equivalent to requiring that
$$ \big|(\langle z \rangle \nabla_z)^\alpha b(z) \big| \leq C_{\alpha} \langle z \rangle^j$$
for all $\alpha \geq 0$ and $z \in M$, where we use the metric to measure the length of the tensor $\nabla^\alpha b$. 

Similarly, a smooth function $b: \Delta_{o(1)} \to \C$ defined on a small diagonal neighbourhood
$\Delta_{o(1)} \subset M^2$ of $M^2$ is said to be a \emph{local product symbol} if we have the bounds 
$$ |(\langle z \rangle \nabla_z)^\alpha (\langle z' \rangle \nabla_{z'})^\beta b(z,z')|
\leq C_{\alpha,\beta} \langle z \rangle^j \langle z' \rangle^{j'}$$
for all $\alpha, \beta \geq 0$ and all $(z,z') \in \Delta_{o(1)}$.   Thus for instance 
any function which extends smoothly to $\overline{M}$ is a symbol of order 0 on $M$, while $\langle z \rangle^j$ is a symbol of order $j$, and if $b_j(z)$ and $b_{j'}(z)$ are symbols of order $j$ and $j'$ respectively then 
$b_j(z) b_{j'}(z')$ is a local product symbol of order $(j,j')$.
\end{definition}

For a conic metric, the cosine rule \eqref{sharp-cosine} implies in particular
$\Phi_\conic(z,z') - (\langle z \rangle^2 + \langle z' \rangle^2)/2$ is 
a product symbol of order $(1,1)$.  For asymptotically conic metrics, there is no exact cosine rule, nevertheless
one can approximate $\Phi$ by $\Phi_\conic$ up to a lower order error:

\begin{lemma}[Phase regularity estimates \cite{HTW}]\label{metric-smooth}  Let $g_\conic$ be the perfectly conic metric
associated to $g$, and let $\Phi_\conic$ be the associated phase function.  Then for $d_{\overline M}(z,z') = o(1)$,
we can write $\Phi(z,z') = \Phi_\conic(z,z') + e(z,z')$, where we have the error estimates
\begin{align}
e(z,z') &= o( \Phi(z,z') ) \label{e-zeroderiv}\\
\nabla_z e(z,z'), \nabla_{z'} e(z,z') &= o( \langle z \rangle + \langle z' \rangle ) \label{e-onederiv} \\
\nabla_z \nabla_{z'} e(z,z') &= o( 1 ) \label{e-mixedderiv} \\
\nabla_z \nabla_{z} e(z,z') &= o( 1 + \frac{ \langle z' \rangle}{\langle z \rangle} ) \label{e-zzderiv} \\
\nabla_{z'} \nabla_{z'} e(z,z') &= o( 1 + \frac{ \langle z \rangle}{\langle z' \rangle} ). \label{e-zzpderiv} 
\end{align}
Furthermore, the function $\Phi(z,z') - (\langle z \rangle^2 + \langle z' \rangle^2)/2$ is a local product symbol of order $(1,1)$.
\end{lemma}

\begin{proof}  For short-range metrics, the claims follow from \cite[Proposition 9.4]{HTW} (in fact, a more precise
estimates are proven there).  The arguments extend to long-range metrics also but are a little more technical; we give the details in Appendix \ref{longrange}.
\end{proof}

Next, we describe certain non-degeneracy and regularity estimates on the phase $\Phi(z,z')$.  

\begin{lemma}[Phase nondegeneracy estimates]\label{grad1} 
If $z, z', z''\in M$ are within $o(1)$ of each other in the $d_{\overline M}$ metric, then for any $0 \leq \theta
\leq 1$ we have the estimates
\begin{align}
\det(-\nabla_{z''} \nabla_{z} \Phi(z,z'')) &= 1 + o(1)
\label{det-est}\\
\frac{\Big| \nabla^{z''} \big( \Phi(z'',z) - \Phi(z'',z')
\big) \Big|_{g(z'')}}{d_M(z,z')} &= 1 + o(1) \label{grad}\\ 
\frac{\Big| \nabla^{z''} \big( (1-\theta)\Phi(z'',z) + \theta
\Phi(z'',z') \big) \Big|_{g(z'')}}{d_M(z'', \gamma_{z \to z'}(\theta))}
&\geq 1 - o(1). \label{gradw}
\end{align}
We emphasize that the error terms $o(1)$ go to zero as $\epsilon \to 0$ uniformly in the choice of
$z, z', z'', \theta$.
\end{lemma}

As an easy corollary of the above estimates we can obtain some further regularity bounds on
$\gamma_{z \to z'}(\theta)$:

\begin{lemma}[Geodesic regularity estimates]\label{geodesic-lemma}  We have the estimate
\begin{equation}\label{comparable}
 C^{-1} \leq \frac{\langle \gamma_{z \to z'}(\theta) \rangle}{\theta \langle z' \rangle + (1-\theta) \langle z \rangle}
\leq C
\end{equation}
whenever $0 \leq \theta \leq 1$ and $d_{\overline M}(z,z')$ is sufficiently small.  
In fact in this region we have the additional symbol-type estimates
\begin{equation}\label{symbol} 
|(\langle z \rangle \nabla_z)^\alpha (\langle z' \rangle \nabla_{z'})^\beta \gamma_{z \to z'}(\theta)|
\leq C_{\alpha,\beta} (\theta \langle z' \rangle + (1-\theta) \langle z \rangle)
\end{equation}
for all $\alpha, \beta \geq 0$ with $\alpha + \beta > 0$.  Furthermore we have the estimates
\begin{equation}\label{lipschitz}
 C^{-1} d_M(z', z'') \leq |\gamma'_{z_0 \to z'}(0) - \gamma'_{z_0 \to z''}(0)|_{g(z_0)} \leq C d_M(z',z'')
\end{equation}
whenever $z_0, z', z'' \in M$ are sufficiently close in the compactified metric $d_{\overline M}$.
\end{lemma}

The proof of these results will be deferred to Appendix \S\ref{sec:geometry}.

\begin{remarks} The estimate \eqref{det-est} can be interpreted as the
  statement as the ``mixed'' Hessian matrix $-\nabla_{z''} \nabla_{z} \Phi(z,z'')$,
  when measured in an orthonormal frame around $T_z M$ and $T_{z''} M$, is approxiamtely unimodular.
  Note that for Euclidean space, the
  quantities in the middle of \eqref{det-est}, \eqref{grad}, and \eqref{gradw}
  are precisely equal to $1$.  These estimates are
  important for establishing (via stationary phase arguments) basic Fourier
  integral operator type estimates for oscillatory integrals whose phase is
  given by $\Phi$.  The estimate \eqref{lipschitz} asserts
that normal co-ordinates are bilipschitz with respect to the metric $d_M$ as long as one works on a small 
ball in $d_{\overline M}$.
\end{remarks}

We now turn to the amplitude $a(z,z')$,
defined in \eqref{a-def}.  From \eqref{det-est} we see that $a$ is well-defined and positive on $\Delta_\epsilon$
for $\epsilon$ sufficiently small. We in fact have a number of additional estimates:

\begin{lemma}[Amplitude estimates and identities]\label{a-smooth}  The functions $a$, $a^{-1}$ are symmetric and are local product symbols of order $(0,0)$.  We also have the diagonal identity
\begin{equation}\label{diagonal}
 a(z,z) = 1 \hbox{ for all } z \in M
\end{equation}
and the transport equations
\begin{equation}\label{right-transport}
 (\nabla^{z'} \Phi\cdot \nabla_{z'} + \frac{1}{2} \Delta_{z'} \Phi - \frac{n}{2}) a(z,z') = 
 (\nabla^{z} \Phi \cdot \nabla_{z} + \frac{1}{2} \Delta_z \Phi - \frac{n}{2}) a(z,z') = 0
\end{equation}
for $d_{\overline{M}}(z,z')$ sufficiently small. Here $\Delta_z$ and $\Delta_{z'}$ denotes the Laplace-Beltrami operator
in the $z$ and $z'$ variables, and we use the metric $g$ to raise the derivatives $\nabla_z$, $\nabla_{z'}$ in the
 usual manner.
\end{lemma}

\begin{remark}\label{unitary}  The equations \eqref{right-transport} are the natural Hamilton-Jacobi equations associated
to the Schr\"odinger flow \eqref{nls} and the phase $\Phi$, which justifies the use of $a(z,z')$ as an amplitude
for a parametrix \eqref{parametrix} for $e^{-itH}$; one can also use \eqref{diagonal}, \eqref{right-transport} instead of \eqref{a-def} as one's definition for $a$.  The choice \eqref{a-def} can also be justified on the grounds that it will make our parametrix \eqref{parametrix} ``locally unitary''. Namely, if we regard \eqref{parametrix} as a Fourier integral operator, then choosing $a$ as in \eqref{a-def} means that the symbol of the operator is equal to $1$, which will ensure that the parametrix is an approximate semigroup (see Theorem~\ref{sfio-concatenate}). 
\end{remark}

\begin{proof} The claim \eqref{diagonal} follows by working in normal
  co-ordinates around $z$, and observing that $\Phi$ agrees with the
  Euclidean counterpart $\frac{1}{2} |z-z'|^2$ in those co-ordinates up to
  errors of cubic or higher order.  The smoothness follows from Lemma
  \ref{metric-smooth} and \eqref{det-est}, since the mixed Hessian of
  $(\langle z \rangle^2 + \langle z' \rangle^2)/2$ vanishes.  The symmetry
  of $a$ is clear from definition.  Now we prove \eqref{right-transport}.
By symmetry, it suffices to establish the second of the equalities in \eqref{right-transport}.  
We rewrite this identity in terms of the half-density $a(z,z') dg(z)^{1/2} \, dg(z')^{1/2}$
as
\begin{equation}\label{lie-identity}
 {\mathcal L}_{\nabla^z \Phi \cdot \nabla_z} (a(z,z') dg(z)^{1/2} \, dg(z')^{1/2}) - \frac{n}{2}
a(z,z') dg(z)^{1/2} \, dg(z')^{1/2} = 0,
\end{equation}
where ${\mathcal L}$ denotes the Lie derivative.

To prove \eqref{lie-identity}, we shall work on the product manifold 
$$N := \RR_+ \times \Delta_{o(1)}
= \{ (t,z,z'): t > 0, d_{\overline M}(z,z') < o(1) \}.$$
We parameterize the cotangent bundle $T^* N$ using ``scattering co-ordinates''
by $(t,z,z',\tau,\zeta,\zeta')$, where $\tau$, $\zeta$, $\zeta'$ are the dual variables to
$-\partial_t/t^2$, $\nabla_z/t$, and $\nabla_{z'}/t$; in other words, the canonical one-form is
$$ -\tau \frac{dt}{t^2} + \zeta \cdot\frac{dz}{t} + \zeta' \cdot \frac{dz'}{t}.$$
We now view rewrite the eikonal equation \eqref{eikonal-variant} semi-classically, noting that the symbol of
$i t^2 \partial_t - t^2 H$ is $\sigma := \tau - \frac{1}{2} |\zeta|_{g(z)}^2$, and that the graph
$L \subset T^* N$ of the derivative of the phase function $\Phi/t$ is given by the $2n+1$-dimensional manifold
\begin{align}
L &:= \graph d\left(\frac\Phi t\right) \nonumber\\
&= \{
(t,z,z',\Phi(z,z'),\nabla_z \Phi(z,z'), \nabla_{z'} \Phi(z,z')): (t,z,z') \in N \};\label{L-form}
\end{align}
note that our use of scattering co-ordinates allows us to extend this manifold smoothly to the boundary $t=0$.
The eikonal equation \eqref{eikonal} (or \eqref{eikonal-variant}) shows that the symbol $\sigma$ vanishes on $L$.
On the other hand, since $L$ is the graph of an exact form, it is a Lagrangian submanifold of $T^* N$.  
Thus if we let $X$ denote the Hamilton vector field associated to $\sigma$, then $X$ preserves $L$.
This vector field $X$ can also be computed explicitly (cf. \cite[(3.1)]{Hassell-Wunsch1}) as
\begin{equation}\label{x-form}
X = t(t \pa_t + \zeta \cdot \pa_\zeta + \zeta' \cdot \pa_{\zeta'} - \abs{\zeta}^2 \pa_\tau + H_{z,\zeta})
\end{equation}
where $H_{z,\zeta}$ generates geodesic flow on $M$.  Thus $t^{-1} X$ extends smoothly (in scattering co-ordinates)
to the boundary $\{ t=0\}$ of $T^* N$, and preserves that boundary, and thus also preserves 
the $n$-dimensional manifold $L|_{t=0}$.

We now restrict to the boundary $t=0$ of $T^* N$, and consider the canonical 
half-density $\alpha := (dz d\zeta)^{1/2}$.  Since $H_{z,\zeta}$ preserves the symplectic
form $d\zeta\wedge dz,$ we easily compute
$$
\mathcal{L}_{t^{-1} X} (\alpha) -\frac n2 \alpha = 0.
$$
We then restrict this identity to $L|_{t=0}$ (which as observed before, was preserved by $t^{-1} X$)
and then push forward from $L|_{t=0}$ down to $\Delta_{o(1)}$ by the projection map
$$(0,z,z',\Phi(z,z'),\nabla_z \Phi(z,z'), \nabla_{z'} \Phi(z,z')) \mapsto (z,z').$$
When one does so, the vector field $t^{-1} X$ becomes $\zeta \cdot \nabla^z = \nabla^z \Phi \cdot \nabla_z$
by \eqref{x-form}, \eqref{L-form}
and the half-density $\alpha_L$ becomes $a(z,z') dg(z)^{1/2} \, dg(z')^{1/2}$ by \eqref{L-form}, \eqref{a-def}.  Thus 
we obtain \eqref{lie-identity} as desired.
\end{proof}


\section{Local Schr\"odinger integral operators (LSIOs)}\label{sec:LSIO}

In our arguments we shall frequently be dealing with a number of
oscillatory integral operators of the form \eqref{parametrix}, as
approximate parametrices for the Schr\"odinger propagators $e^{-itH}$ and
associated operators.  These operators are almost of the standard
oscillatory integral type, except that the support of the amplitude
function is on a thin cone rather than a compact set, which unfortunately
means that we cannot quite cite a reference for the (standard) results for
these operators we will need here. Nevertheless, the control on the
geometry of thin cones obtained in the previous section is good enough that
it turns out that the basic calculus of these operators is almost identical
to those of standard oscillatory integral operators on compact sets.

We now define the family of oscillatory integral operators we shall use.

\begin{definition}  A \emph{local Schr\"odinger integral operator (or LSIO)
 } is a family of operators  $S(t) = S_b(t)$ parametrized by $t \in
  [-1,1]$ given by the formula
\begin{equation}\label{S-def}
 S_b(t) f(z) := \frac{1}{(2\pi i t)^{n/2}} \int_M e^{i \Phi(z,z')/t} a(z,z') b(z,z',t) f(z')\ dg(z')
\end{equation}
when $t \neq 0$ (compare with \eqref{parametrix}) and
\begin{equation}\label{S0-def}
S_b(0) f(z) := b(z,z,0) f(z)
\end{equation}
when $t=0$, where $b(\cdot, \cdot,t)$ is an $L^\infty(t)$ family of local
product symbols of order (0,0) (which we refer to as the \emph{symbol of
$S(t)$}).  Here we take the standard branch of $z^{1/2}$ (with branch cut
at the negative real axis) to define $(2\pi i t)^{n/2}$.
\end{definition}

\begin{remarks}  One can heuristically think of an LSIO $S_b(t)$ as a Fourier integral operator (FIO) of order 0,
whose canonical relation is given by geodesic flow for time $t$; such
operators can track the Schr\"odinger flow (and operators related to that
flow) accurately as long as one does not travel too far in the compactified
metric $d_{\overline M}$, but become useless once the distances in
$d_{\overline M}$ become large.  Unfortunately as our spatial parameters
can extend to the boundary $\partial M$ of our compactified manifold, it
is not easy to use the standard theory of FIOs directly to handle LSIOs,
and so we must redevelop much of this theory in our setting.
From a geometrical viewpoint it turns out that it is more natural to view the amplitude
function a half-density $a(z,z') dg(z)^{1/2} dg(z')^{1/2}$ rather than
as a scalar function $a(z,z')$ (cf. \eqref{a-def}), and similarly to view
$S_b(t)$ as acting on half-densities $f(z') dg(z')^{1/2}$ rather than on
scalar functions $f(z')$ (which is of course extremely natural for the
$L^2$ theory, although less so for the $L^p$ theory); we already saw this in the proof of \eqref{right-transport}. 
 Thus for instance
\begin{equation}\label{invariant} 
\begin{split}
&\int_M f dg^{1/2} S_b(t)(f' dg^{1/2}) =\\
&\frac{1}{(2\pi i t)^{n/2}} \int_{M^2}  e^{i \Phi(z,z')/t} 
b(z,z',t) (a(z,z') dg(z)^{\frac1{2}} dg(z')^{\frac1{2}})
f(z)\ dg(z)^{\frac1{2}} f'(z')\ dg(z')^{\frac1{2}}.
\end{split}
\end{equation}
In this formulation $S_b(t)$ becomes manifestly invariant under co-ordinate changes (with $b$ being a scalar
function).  We shall thus adopt this invariant viewpoint when convenient (e.g. if we want to work in normal co-ordinates).
\end{remarks}

\begin{remark}
Naively, it would natural to approximate the Schr\"odinger flow $e^{-itH}$
by the global parametrix $S_1(t)$.  Unfortunately, neither the phase $\Phi$ from
\eqref{Phi-def} nor the amplitude $a$ from \eqref{a-def} will be smooth  away from $\Delta_\ep$
in general\footnote{In the case
  when $M$ has non-positive curvature everywhere, then the phase function
  $\Phi$ is globally smooth, and it may be possible to use $S_1(t)$ as a
  global parametrix without recourse to any cutoff.  We have not pursued
  any simplifications in this case however.}.  Thus we must apply an
additional cutoff $\chi$, which on the one hand makes our parametrix into a
``local'' operator, with phase supported in $\overline{\Delta_{\ep}}$ with
$\ep\ll 1,$ but which creates additional error terms that must be dealt
with.  It turns out that the resulting errors can be controlled by local
smoothing estimates, at least for the purposes of proving non-endpoint
Strichartz estimates.
\end{remark}

\begin{remark}
While our parametrix will have a time-independent symbol, error terms will
crop up involving time-dependent symbols.  It is for this reason that we
allow time-dependence in $b.$  The more accurate parametrix constructed in
\cite{Hassell-Wunsch1}, which satisfies the Schr\"odinger equation up to
$O(t^\infty),$ also has a time-dependent symbol.
\end{remark}

We now give the basic properties of LSIOs.  

\begin{lemma}\label{sfio-basic} Let $b$ be a local product symbol of order (0,0) and let $|t| \leq 1$.  If $\tilde b$
is the local product symbol of order (0,0) defined by $\tilde b(z,z',t) :=
 \overline{b(z',z,-t)}$ then $S_b(t)^* = S_{\tilde b}(-t)$.  In particular
 if $b$ is real and symmetric in $z,z'$ and constant in $t$ then $S_b(-t) =
 S_b(t)^*$.  Also, if $b(t)$ is continuous in $t$ at $t=0,$ we have the
 continuity property
\begin{equation}\label{pointwise}
 \lim_{t \to 0} S_b(t) f(z) = S_b(0) f(z)
\end{equation}
for any test function $f$, and for both the $t \to 0^+$ and $t \to 0^-$ approach, where the limit can be taken in either the pointwise or the distributional sense\footnote{In fact one has convergence in much smoother topologies, such as the Schwartz topology, but we will not need to use this here.}.
\end{lemma}

\begin{proof}
The first set of claims follow from \eqref{invariant} since $\Phi$ and $a$ is real and symmetric.  The 
claim \eqref{pointwise} can be proven pointwise for a single $z$ by 
working in normal co-ordinates around $z$ and using stationary phase (see e.g. \cite{stein:large}),
observing from Gauss's Lemma \eqref{gauss-lemma} that $\Phi(z,z')$ is only stationary at $z'=z$, and in normal co-ordinates is equal to
$\frac{1}{2} |z-z'|^2$ plus terms of cubic or higher order.  The distributional convergence then follows by using 
\eqref{diagonal} and the hypothesis that $f$ is a test function and quantitative stationary phase estimates to place uniform bounds on 
$S_b(t) f$; we omit the details.
\end{proof}

The basic $L^2$ estimate (analogous to the H\"ormander-Eskin theorem \cite{hormander}, \cite{eskin} for FIOs) 
is given by the following theorem, which we prove in Appendix \S \ref{FIO-sec}.

\begin{theorem}[$L^2$ estimate for LSIOs]\label{sfio-l2}  If $S_b(t)$ is a LSIO then we have the estimate
$$ \| S_b(t) f \|_{L^2(M)} \leq C_b \| f \|_{L^2(M)}$$ for all test
functions $f$ and all times $|t| \leq 1$.  The bound $C_b$ depends only on
$M$ and on finitely many of the constants in the product symbol bounds for
the symbol $b.$
\end{theorem}

From \eqref{S-def} it is clear that we have the $L^1 \to L^\infty$ estimate
\begin{equation}\label{easy-dispersive}
 \| S_b(t) f \|_{L^\infty} \leq C_b |t|^{-n/2} \| f \|_{L^1}
\end{equation}
for all test functions $f$ and all $0 < |t| < 1$.  We can generalize this estimate as follows.

\begin{theorem}[Dispersive estimate for SFIOs]\label{sfio-infty}  If $b$, $b'$ are local product symbols of order (0,0), then we have the estimate
\begin{equation}\label{estimate-dispersive}
 \| S_b(t)^* S_{b'}(s) f \|_{L^\infty(M)} \leq C_{b,b'} |t-s|^{-n/2} \| f \|_{L^1(M)}
\end{equation}
for all test functions $f$ and all times $-1 \leq s,t \leq 1$ with $s \neq t$.
The bound $C_{b,b'}$ depends only on finitely many of the constants in the product symbol 
bounds for $b$, $b'$.
\end{theorem}

We prove this theorem in Appendix \S \ref{FIO-sec}.  From these two theorems and the abstract Strichartz
estimate in \cite[Theorem 1.2]{tao:keel} we thus have the Strichartz estimate
\begin{equation}\label{S-endpoint}
 \| S_{b}(t) f \|_{L^q_t L^r_x([-1,1] \times M)} \leq C_b \| f \|_{L^2(M)}
\end{equation}
for all admissible exponents $(q,r)$ and all test functions $f$, where for each $t$, $b(t)$ is a local product symbol
of order (0,0), 
with bounds uniform in $t$. Indeed we even obtain the endpoint
case $q=2$ this way (as long as $n \geq 3$).  Unfortunately we will
not be able to exploit this endpoint as we also use the Christ-Kiselev lemma (see Lemma \ref{ck-lemma}), which
requires $q > 2$.  Note that we do not require any regularity of $b(t)$ in time in order to establish \eqref{S-endpoint},
only that each $b(t)$ is a local product symbol of order (0,0) uniformly in $t$.

As discussed in the introduction, our parametrix \eqref{parametrix} can now be interpreted as an LSIO of the form
$S_{\chi a}(t)$ for some suitable cutoff $\chi$, where $a$ is the amplitude \eqref{a-def}.  The motivation for this
parametrix is justified by the following identity, which shows that $S_{\chi a}(t)$ locally solves the Schr\"odinger equation \eqref{nls} in the top two orders in $t^{-1}$:

\begin{lemma}[LSIOs as parametrices]\label{LSIO-operator}  Let $\psi_1, \psi_2: M \to \C$ be two local symbols of order 0, and let $\chi:M^2 \to \C$ be a local product symbol of order (0,0) such that $\chi = 1$ on $\supp(\psi_1) \times \supp(\psi_2)$.
Then for all $|t| \leq 1$, $\psi_1 (\frac{d}{dt} S_{\chi}(t) + i S_{\chi}(t) H) \psi_2$ is an LSIO; in fact we have
the explicit identities
\begin{align*}
\psi_1 (\frac{d}{dt} S_{\chi}(t) + i S_{\chi}(t) H) \psi_2
&= S_{\psi_1(z)  \psi_2(z') (H_{z'} a) / a}(t)\\
\psi_1 (\frac{d}{dt} S_{\chi a}(t) + i H S_{\chi a}(t)) \psi_2
&= S_{\psi_1(z) \psi_2(z') (H_{z} a) / a}(t).
\end{align*}
\end{lemma}

\begin{proof}  We prove just the former claim, as the latter is similar (and can also be deduced from the first
by duality and Lemma \ref{sfio-basic}).  From \eqref{S-def} and hypothesis on $\chi$, we have
$$
\psi_1 (\frac{d}{dt} S_{\chi}(t) + i S_{\chi}(t) H) \psi_2 f(z) =
\int_{M} (\partial_t + i H_{z'})\Big( \frac{e^{i \Phi/t} a}{(2\pi i t)^{n/2}} \Big) 
\psi_1(z) \psi_2(z') f(z')\ dg(z')
$$
By the Leibniz rule and chain rule, we can expand
\begin{equation}\label{t-split}
\begin{split}
(\partial_t + i H_{z'})\Big(\frac{e^{i\Phi/t} a}{(2\pi i t)^{n/2}} 
 \Big)  &= 
\frac{e^{i\Phi/t}}{(2\pi i t)^{n/2}} 
\bigg[ t^{-2} (-i \Phi + \frac{i}{2} |\nabla_{z'} \Phi|^2) a \\
 &+ t^{-1} (-\frac{n}{2} + \frac{1}{2} (\Delta_{z'} \Phi) + \nabla^{z'} \Phi \cdot \nabla_{z'} ) a 
 - i H_z a \bigg].
\end{split}
\end{equation}
But the $t^{-2}$ terms vanish thanks to the eikonal equation \eqref{eikonal} (see also
\eqref{eikonal-variant}), and the
$t^{-1}$ terms vanish thanks to the transport equation \eqref{right-transport}.
\end{proof}

Finally, we require a composition law for LSIOs.  It turns out that one has a good law as long as one composes two
operators moving in the same direction in time.

\begin{theorem}[Composition law for LSIOs]\label{sfio-concatenate}  Let $0 < t, t'$ be such that $t + t' \leq 1$, and let
$b$, $b'$ be local product symbols of order (0,0).  Then we have\footnote{Strictly speaking, for this theorem to work, $b$ and $b'$ have to be supported on $\Delta_{o(1)/2}$ rather than $\Delta_{o(1)}$.  This will however cause no difficulties as we shall only use this theorem a finite number of times and can take $o(1)$ to be as small as we wish.} the composition law
$$ S_b(t) S_{b'}(t') = S_{c}(t+t') + S_{\tilde e_{t,t'}}(t+t')$$
where $c$ is the local product symbol of order (0,0)
\begin{equation}\label{b-weight}
c(z,z') := b(z,w,t) b'(w,z',t') 
\end{equation}
where $w := \gamma_{z \to z'}(\frac{t}{t+t'})$
and $\tilde e_{t,t'} = \frac{t+t'}{\langle z \rangle + \langle z' \rangle} e_{t,t'}$
where $e_{t,t'}$ is another local product symbol of order (0,0) (uniformly in $t, t'$).
\end{theorem}

\begin{remarks}
Note that the point $w$ arises naturally from geometric optics heuristics:
if a particle is at $z'$ at time zero and at $z$ at time $t+t'$, we expect
it to be at $w = \gamma_{z \to z'}(\frac{t}{t+t'}) = \gamma_{z' \to
  z}(\frac{t'}{t+t})$ at time $t'$.  The presence of the amplitude factors
$a$ in \eqref{b-weight} reflects the fact that the operators $e^{-itH}$
form a semigroup, and that the operators $U(t)$ with amplitude $a$ behave
like $e^{-itH}$.  We prove this theorem in Appendix \S \ref{FIO-sec}.
We remark that by using Lemma \ref{sfio-basic} we also have an analogue of
this theorem for negative times $t, t'$, but we will not need to use that
analogue here.  However, the lemma fails when $t, t'$ have opposite sign as
one loses control of the geometry in that case\footnote{The point is that if a geodesic is at $z$ at time $0$
and at $z'$ at time $t+t'$, then when $t, t'$ have opposite sign, we do not necessarily have any good control on
where the geodesic is at time $t$, even if $d_{\overline{M}}(z,z') = o(1)$.}.  The factor
$\frac{t+t'}{\langle z \rangle + \langle z' \rangle}$ is not best possible
(indeed one can obtain a more accurate asymptotic expansion for the error term $\tilde e_{t,t'}$), 
but will suffice for our arguments.  There are also additional regularity estimates for $\tilde e_{t,t'}$
in the time variables $t, t'$, but again we will not need these estimates (because
\eqref{S-endpoint} does not require any regularity of the symbol in time).
\end{remarks}

It seems of interest to develop a more systematic calculus of ``scattering Fourier integral operators'' which would include
these LSIOs as a special case, but can also extend to more non-local operators in which the canonical relation contains
folds (and would therefore not be representable in an oscillatory integral form such as \eqref{S-def}).  We will however not pursue such a theory here.

\section{Scattering symbol calculus and the space $X$}\label{sec:doi}

In this section we set up some standard notation for scattering pseudo-differential operators, as described by
Melrose \cite{melrose}; we will
need this to properly define the $X$ space used in Lemma \ref{ls-lemma}, and also to prove local smoothing estimates
for a certain parametrix for the Schr\"odinger equation which we will encounter later in this paper.  Our notation
is repeated from \cite{HTW}, and so we shall be somewhat brief in our descriptions, referring to \cite{melrose}
and \cite{HTW} for more detail.

The first step is to choose co-ordinates for the cotangent bundle $\Tscstar
\Mbar$ in the scattering region $r > r_0$.  We shall use the co-ordinates
$(x,y^j,\nu,\mu_j)$, where $0 < x < \epsilon_0$, and $y \in \partial M$,
while the co-ordinates $\nu$ and $\mu_j$ are dual to the vector fields
$-x^2 \pa_x = \partial_r$ and $x \partial_{y^j}$ respectively, and the
canonical one-form is thus $-{\nu dx}/{x^2} + {\mu^j dy_j}/{x}$.  Following
Melrose \cite{melrose}, we define the \emph{scattering cotangent bundle}
over the compact manifold $\overline{M}$ as the bundle $\Tscstar \Mbar$
whose sections are locally spanned over $\mathcal{C}^\infty(\overline{M})$
by $dx/x^2$ and $dy/x$ (hence can be paired with the \emph{scattering
  vector fields}, spanned by $-x^2 \pa_x$ and $x \partial_y$).  On any
compact set away from the boundary $\partial M$, we shall parameterize the
cotangent bundle $\Tscstar \Mbar$ in the usual manner by a spatial variable
$z$ and a frequency variable $\zeta$.

\begin{definition}[Scattering symbol classes]\label{def:symb} Let $m,l$ be real numbers and $0 \leq \rho < 1$. 
A smooth function $a: \Tscstar \Mbar \to \C$ is said to be in the \emph{scattering symbol class}
 $S^{m,l}_{1, \rho}(\Mbar)$ provided that one has the usual Kohn-Nirenberg symbol estimates 
$$ | \nabla_z^\alpha \nabla_\zeta^\delta a(z,\zeta) | \leq C_{\alpha,\delta,K} (1 + |\zeta|)^{m-|\delta|}$$
for $z$ ranging in any compact subset $K$ of $M$, and all $\alpha, \beta \geq 0$, as well as the
scattering region estimates
$$
\abs{\partial_x^\alpha \nabla_y^\beta \partial_\nu^\gamma \nabla_\mu^\delta
a(x,y,\nu,\mu)}\leq C_{\alpha,\beta,\gamma,\delta}
\, (1 + |\nu| + |\mu|)^{m-\gamma-\abs\delta}
x^{l-\alpha - \rho(\gamma + \abs\delta)} 
$$ in the scattering region $(x,y) \in (0,\epsilon_0) \times \partial M$.
Every symbol $a$ in $S^{m,l}_{1,\rho}(\Mbar)$ can be quantized to give a
\emph{scattering pseudo-differential operator} $\Op(a)$ in the class
$\Psisc^{m,l; \rho}(M)$; the exact means of quantization is not
particularly important for our purposes, but we can for instance use the
Kohn-Nirenberg quantization (on any $n$-dimensional asymptotically conic
manifold)
\begin{equation}\label{op-def}
\Op(a) u(z') := (2\pi)^{-n} \int_{T^* M} \chi(z',z'')e^{-i\langle \gamma'_{z' \to z''}(0), \zeta \rangle} a(z',\zeta) u(z'')\, d\zeta dz'',
\end{equation}
to define the kernel of $\Op(a)$ near $\partial M$; here, $\chi$ is a
cutoff to $\Delta_\eps$ for some suitably small $\eps > 0$, and $d\zeta dz''$ is the canonical Liouville measure on the cotangent bundle $T^* M$.
\end{definition}

We now briefly review the scattering calculus for these symbols (see \cite{HTW} and especially \cite{melrose}
for more details).  If $A$ is an operator in $\Psisc^{m,l; \rho}(\Mbar)$, then its (principal) symbol $\sigma(A)$ is well defined in $S^{m,l}_{1,\rho}(\Mbar)/S^{m-1,l+1-\rho}_{1,\rho}(\Mbar)$, and its 
adjoint $A^*$ (with respect to the $L^2(M)$ inner product) is also
in this class with symbol
\begin{equation}\label{symbol-adjoint}
\sigma(A^*) = \overline{\sigma(A)}
\end{equation}
If $A \in \Psisc^{m,l; \rho}(\Mbar)$ and $B \in \Psisc^{m',l'; \rho}(\Mbar)$ then
\begin{equation}\label{scattering-calculus}
AB \in \Psisc^{m+m',l+l'; \rho}(\Mbar); \quad i[A,B] \in \Psisc^{m+m'-1,l+l'+1 - \rho; \rho}(\Mbar).
\end{equation}
Moreover
\begin{equation}\label{scattering-calculus-symbol}
\begin{split}
\sigma(AB) &= \sigma(A) \sigma(B)\\
\sigma(i[A,B]) &= \{ \sigma(A), \sigma(B) \}
\end{split}
\end{equation}
where $\{ , \}$ is the Poisson bracket on the scattering cotangent bundle
$\Tscstar \Mbar$.  In particular $\{\sigma(H), \cdot\}$ is the Hamilton vector field in the direction of
geodesic flow.

We introduce the weighted Sobolev spaces $H^{m,l}(M)$ as
$$ H^{m,l}(M) := \{ u: \langle z \rangle^l u \in H^m(M) \}.$$
From \eqref{scattering-calculus}, \eqref{scattering-calculus-symbol} and $L^2$-boundedness of $\Psisc^{0,0; \rho}(\Mbar)$ it is easy to verify that operators in $\Psisc^{m,l;\rho}(\Mbar)$ map $H^{m',l'}(M)$ to $H^{m'-m,l'+l}(M)$.

We now define a half-angular operator $B$.

\begin{definition}[Half-angular derivative]\label{B-def}
Let $\phi: \R \to \R$ be a smooth non-decreasing function such that $\phi =
0$ on $(-\infty,1]$ and $\phi=1$ on $[2,+\infty)$.  For any $\rho \in
  (0,1),$ we let $B \in \Psisc^{1/2,1/2;\rho}(\Mbar)$ be a quantization of
  the $S^{1/2,1/2}_{1,\rho}(\Mbar)$ symbol
\begin{equation}
\sigma(B) := \frac1{2} \phi(\frac{\epsilon_0}{x})
\phi(|\mu|^2 + |\nu|^2) 
\phi\big(\frac{\abs{\mu}}{|\nu| x^\rho}\big) x^{1/2} 
  |\mu|^{1/2}  
\label{b-defn}\end{equation}
where $|\mu|^2 = |\mu|^2_{h(x)} := h^{jk}(x,y) \mu_j \mu_k$.  By adding a lower order term to $B$ if necessary we may
assume that $B$ is self-adjoint.  We then define the $X$ norm on test functions $f$ as
\begin{equation}\label{X-def}
\| f \|_X := \| f \|_{L^2(M)} + \| f \|_{H^{1/2,-1/2-\rho/2}(M)} + \| B f \|_{L^2(M)}.
\end{equation}
\end{definition}

\begin{remark} Ignoring all the cutoffs, the symbol $\sigma(B)$ is essentially $\frac{1}{2} x^{1/2} |\mu|^{1/2}$;
thus $B$ should be thought of heuristically as $B \approx \langle z \rangle^{-1/2} |\nabb|^{1/2}$.  The cutoffs
are designed to keep $B$ close to $\partial_M$, away from the frequency origin $(\mu,\nu) = 0$, and
away from the radial directions $\mu = 0$, in order to avoid a number of singularities in the symbol.
\end{remark}

The local smoothing estimate in Lemma \ref{ls-lemma} then follows directly from \cite[Lemmas 10.4-10.5]{HTW}.  Later on we will also need to modify that argument so that it gives
a local smoothing estimate for our parametrix.

A key fact about the $X$ norm is its ability to absorb ``half'' of a commutator $i[H,\psi]$; the main point is that
that commutator consists primarily of an angular derivative, and each $X$ norm can absorb roughly half of such
an angular derivative.  A more precise statement is the following.

\begin{lemma}\label{X-split}  Let $\psi: \overline{M} \to \R$ be a smooth function.  Then we have the estimate
$$ |\langle f, i[H,\psi] g \rangle_{L^2(M)}| \leq C_\psi \|f\|_X \|g\|_X$$
for all test functions $f, g$.  
\end{lemma}

\begin{proof}  If $f$ or $g$ is supported on the near region $r \leq 2r_0$ then the claim follows simply by
observing from \eqref{X-def} that the $X$ norm controls the $H^{1/2}(M)$
norm in this region, and that $i[H,\psi]$ is a first-order operator in this
region with bounded coefficients and thus maps $H^{1/2}(M)$ to
$H^{-1/2}(M)$.  Now suppose that $f, g$ are supported in the scattering
region $r > r_0$.  Ignoring the zeroth order component of $i[H,\psi]$
(which can be dealt with by the $L^2$ component of the $X$ norm in
\eqref{X-def}), we see from the Liebniz rule that we can write $i[H,\psi]
= a x^4 D_x + b \cdot x^2 D_y$ for some symbols $a(z)$, $b(z)$ of order 0;
the point is that the angular derivatives of $\psi$ decay like $O(\langle z
\rangle^{-1})$ but the radial derivatives of $\psi$ will decay like
$O(\langle z \rangle^{-2})$ since $\psi$ is smooth on the boundary
$\overline{M}$.  Using the scattering calculus
\eqref{scattering-calculus-symbol} and the fact that $\mu x$ is bounded by
a multiple of $\nu x^{1+\rho}$ on $\supp (1-\phi(\frac{\abs\mu}{\abs\nu x^\rho})),$ we can thus write
$$ i[H,\psi] \in B^* \Psisc^{0,0;\rho}(M) B + \Psisc^{1,1+\rho;\rho}(M).$$
Since elements of $ \Psisc^{1,1+\rho;\rho}(M)$ map
$H^{1/2,-1/2-\rho/2}(M)$ to $H^{-1/2,1/2+\rho/2}(M)$ and those of
$\Psisc^{0,0;\rho}(M)$ map $L^2(M)$ to $L^2(M)$, the corresponding terms
are controlled by the last two terms in \eqref{X-def}.
\end{proof}

\begin{remark} The commutators $i[H,\psi]$ arise naturally when trying to approximate the flow $e^{-itH}$ by certain
approximate flows localized to the support of $\psi$---see \S\ref{sec:proof}.  Lemma \ref{X-split}, combined with the local smoothing estimate in Lemma \ref{ls-lemma} (and its analogue for the parametrix in Proposition~\ref{parametrix-smoothing}) will be essential in allowing us to pass from Strichartz 
control of the localized flows to Strichartz control on the global flow.
\end{remark}


\section{Construction of the parametrix}\label{parametrix-sec}

In this section we define precisely our parametrix for the Schr\"odinger propagator $e^{-itH}$, and prove a smoothing estimate for the parametrix which is the exact analogue of Lemma~\ref{ls-lemma}. We start by constructing a family of partitions of unity and associated cutoff functions on $M$. 

\begin{lemma}\label{psilemma}  Let $\epsilon > 0$ be a small number.  Then
  there exist local symbols $\psi^{(j)}_\alpha: M \to \R$ for $j=0,1,2,3,4$ of order 0, each supported
on a ball of radius $\epsilon$ in the $d_{\overline M}$ metric, which enjoy the partition of unity property 
\begin{equation}\label{partition}
 1 = \sum_{\alpha \in A} \psi^{(0)}_\alpha,
\end{equation}
the nesting property
\begin{equation}\label{nesting}
\psi^{(j)}_\alpha \psi^{(j')}_\alpha = \psi^{(j)}_\alpha \hbox{ for all } 0
\leq j < j' \leq 4,\ \alpha \in A,
\end{equation}
and the square root property
\begin{equation}\label{phi-sq}
\psi^{(j)}_\alpha,\ \nabla \psi^{(j)}_\alpha \ \text{are squares of smooth functions.}
\end{equation}
Also, we can strengthen the nesting property \eqref{nesting} in two different ways:
\begin{itemize}
\item[(i)] For any $0 \leq j \leq 3$, $\alpha \in A$ and any two points $z, z' \in \supp(\psi^{(j)}_\alpha)$, the geodesic $\gamma_{z \to z'}([0,1])$ is contained in the region where $\psi^{(j+1)}_\alpha = 1$.
\item[(ii)] For any $0 \leq j \leq 2$, $\alpha \in A$
there exists $c > 0$ such that, if $z' \in \supp \psi^{(j)}_\alpha$ and $z \in \supp \nabla \psi^{(j+1)}_\alpha$, then 
\begin{equation}\label{zetasign}
\nabla_z \psi^{(j+1)}_\alpha \cdot \gamma'_{z \to z'}(0) \geq c |\gamma'_{z
  \to z'}(0)|_{g(z)} | \nabla_z \psi^{(j+1)}_\alpha |_{g(z)^{-1}}.
\end{equation}
\end{itemize}
\end{lemma}

\begin{proof}
Let $z_\alpha$ be an arbitrary point in $M$, and let $\psi^{(4)}_\alpha$ be
a local symbol of order 0 which equals 1 on a small ball
$B_{\overline{M}}(z_\alpha,r_4)$ in the compactified metric $d_{\overline{M}}$
centered at $z_\alpha$.  By Lemma \ref{glc} we can then find an even
smaller ball $B_{\overline{M}}(z_\alpha,r_3)$ centered at $z_\alpha$ with
the property that any minimal length geodesic connecting two points in that
neighbourhood lies in the region where $\psi^{(4)}_\alpha = 1$. We now
claim that there is a local symbol $\psi^{(3)}_\alpha$ of order 0 supported in
$B_{\overline{M}}(z_\alpha,r_3)$ which equals one on a smaller ball
centered at $z_\alpha$, and with the property that $\psi^{(3)}_\alpha$ is
monotone decreasing on any outward ray from $z_\alpha$.  This is clear for
$z_\alpha \in K_0$, since we can take $\psi^{(3)}_\alpha$ to be a function
of the distance $d_M(z_\alpha, \cdot)$, in which case \eqref{zetasign}
holds with $c=1$. Moreover, by continuity we can find a smaller ball
$B_{\overline{M}}(z_\alpha,r_2)$ such that \eqref{zetasign} holds, with
$c=1/2$, say, for all $z' \in B_{\overline{M}}(z_\alpha,r_2)$ and $z' \in
\supp \nabla \psi^{(3)}_\alpha$. Now consider $z_\alpha$ in the scattering
region $M \setminus K_0$. Suppose first that $M$ is perfectly conic.  Then we use the fact
from Lemma \ref{glc} that the geodesic curves given by
$g_\conic$ are also the geodesic curves given by the compactified metric $\overline{g}$.
Then we can take $\psi^{(3)}_\alpha$ to be a smooth
monotone decreasing function of distance with respect to the compactified
metric $\overline{g}$; this is decreasing along every geodesic emanating from $z' =
z_\alpha$, and hence, by continuity, decreasing along every geodesic
emanating from $z'$ in a ball $B_{\overline{M}}(z_\alpha,r_2)$, where $r_2$
may be taken uniform over the scattering region. Therefore \eqref{zetasign}
holds, with the positivity of $c$ following from compactness.  Since in the
general case the metric $g$ is very close to the perfectly conic metric $g_\conic$ 
for $x$ sufficiently small,
$$
\big( 1-\frac{c}{2} \big) g_{\conic} \leq g \leq  \big(1+ \frac c 2 \big) g_{\conic} \text{ for $x$ small},
$$
\eqref{zetasign} holds also for $g$ (with constant $c/2$, say) for $x$ sufficiently small. 
This concludes the proof of \eqref{zetasign}.

We can then find a local
symbol $\psi^{(2)}_\alpha$ of order 0 supported on $B_{\overline{M}}(z_\alpha,r_2)$
which equals 1 in an even smaller ball $B_{\overline{M}}(z_\alpha,r_1)$,
such that any minimal length geodesic connecting two points in $\supp
\psi^{(2)}_\alpha$ lies in the set where $\psi^{(3)}_\alpha = 1$.  We can
likewise find a local symbol $\psi^{(1)}_\alpha$ of order 0 supported inside
$B_{\overline{M}}(z_\alpha,r_1)$ and equal to 1 in an even smaller ball
around $z_\alpha$ satisfying (i) and (ii) with $j=1,$ and then construct
$\psi^{(0)}_\alpha$ similarly.  Note that none of these radii of balls used
in this construction need to depend on $z_\alpha$.

As $z_\alpha$ ranges over $M$, the balls $B_{\overline{M}}(z_\alpha,r_0)$
cover $\overline{M}$.  By compactness and a standard partition of unity
construction we can thus find a finite index set $A$ and a local symbol
$\psi^{(0)}_\alpha$ of order 0 for each $\alpha \in A$ supported on
$B_{\overline{M}}(z_\alpha,r_0)$ obeying \eqref{partition}.  The claims of
the lemma are then easily verified.  (To ensure the property \eqref{phi-sq}
we simply choose all cutoff functions above to have smooth square roots and
derivatives.)
\end{proof}

We now formally define the local parametrices that we shall use:
\begin{definition}
For $\alpha \in A$, let $U_\alpha(t)$ be the LSIO with symbol 
$b_\alpha(z,z') := \psi^{(2)}_\alpha(z) \psi^{(2)}_\alpha(z')$;
this is a product symbol of order $(0,0)$ by Lemma \ref{a-smooth}.  
\end{definition}
From Lemma \ref{sfio-basic} we have the symmetry condition
\begin{equation}\label{symmetry} U^*_\alpha(t) = U_\alpha(-t)
\end{equation}
and the condition
\begin{equation}\label{identity} \operatorname{w-lim}_{t \to 0} U_\alpha(t) = U_\alpha(0) = \psi^{(2)}(z)^2,
\end{equation}
where the limit is meant in the distribution sense and the right hand side is regarded as a multiplication operator. 

Recall from Section~\ref{sec:LSIO} that $U_\alpha(t)$ is bounded on $L^2(M)$ uniformly in time (Theorem~\ref{sfio-l2}), satisfies a composition law (Theorem~\ref{sfio-concatenate}), and an $L^1 \to L^\infty$ estimate (Theorem~\ref{sfio-infty}), and thus also obeys a Strichartz estimate \eqref{S-endpoint}. We shall now round out this collection of estimates
by also demonstrating that $U_\alpha(t)$ satisfies a smoothing estimate which is the analogue of Lemma~\ref{ls-lemma}. 

\begin{proposition}[Half-angular local smoothing effect for parametrix]\label{parametrix-smoothing}
The para\-metrix $U_\alpha(t)$ satisfies the smoothing estimate
$$ 
\| U_\alpha(t) \tilde f \|_{L^2([-1,0]; X)} \leq C \| \tilde f \|_{L^2(M)}$$
(cf. Lemma \ref{ls-lemma}).  
\end{proposition}

\begin{proof}
The $L^2(M)$ portion of \eqref{X-def} of this estimate follows from 
Theorem \ref{sfio-l2}, so it will suffice to prove the estimates
\begin{equation}\label{ls-usual}
 \| U_\alpha(t) \tilde f \|_{L^2_t H^{1/2,-1/2-\rho/2}([-1,0] \times M)} + 
\| B U_\alpha(t) \tilde f \|_{L^2_t L^2([-1,0] \times M)} \leq C \| \tilde f \|_{L^2(M)}.
\end{equation}
For technical reasons it will be convenient to replace $U_\alpha$ by the
slightly larger operator $\psi^{(3)}_\alpha \tilde U_\alpha$, where $\tilde
U_\alpha(t)$ is defined the same as $U_\alpha(t)$ but with amplitude
$\psi^{(4)}_\alpha(z) a(z,z') \psi^{(2)}_\alpha(z')$ instead of
$\psi^{(2)}_\alpha(z) a(z,z') \psi^{(2)}_\alpha(z')$.  This replacement is
justified since the space $H^{1/2,-1/2-\rho/2}$ is stable under
multiplication by cutoff functions such as $\psi^{(j)}_\alpha$ (which are
special cases of operators in $\Psisc^{0,0;\rho}(\Mbar)$), and because the
commutator of $B$ with $\psi^{(3)}_\alpha$ is of zeroth order and hence
gives an acceptable contribution to \eqref{ls-usual} by Theorem
\ref{sfio-l2}.

We therefore aim to prove

\begin{equation}\label{ls-variant}
 \|\psi^{(3)}_\alpha \tilde U_\alpha(t) \tilde f \|_{L^2_t H^{1/2,-1/2-\rho/2}([-1,0] \times M)} + 
\| B \psi^{(3)}_\alpha \tilde U_\alpha(t) \tilde f \|_{L^2_t L^2([-1,0] \times M)} \leq C \| \tilde f \|_{L^2(M)}.
\end{equation}

We use the positive commutator method.  Let $A \in \Psisc^{0,0;\rho}(\Mbar)$
be a self-adjoint scattering pseudo-differential operator to be chosen later; the constants below are allowed
to depend on $A$.  Since $A$ is bounded on $L^2(M)$, we see from Theorem \ref{sfio-l2} that
$$ |\langle A \psi^{(3)}_\alpha \tilde U_\alpha(t) \tilde f, \psi^{(3)}_\alpha \tilde U_\alpha(t) \tilde f \rangle|
\leq C \| \tilde f \|_{L^2(M)}^2$$
for all $-1 \leq t \leq 1$.  Similarly, from Lemma \ref{LSIO-operator} we see that $\psi^{(3)}_\alpha (\partial_t + iH) \tilde U_\alpha(t)$ is an LSIO, and hence by Lemma \ref{sfio-l2} again
$$ |\langle A \psi^{(3)}_\alpha (\partial_t + iH) \tilde U_\alpha(t) \tilde f, \psi^{(3)}_\alpha \tilde U_\alpha(t) \tilde f \rangle|
\leq C \| \tilde f \|_{L^2(M)}^2$$
Thus from the identity
\begin{align*}
\frac{d}{dt} \langle A \psi^{(3)}_\alpha \tilde U_\alpha(t) \tilde f, \psi^{(3)}_\alpha \tilde U_\alpha(t) \tilde f \rangle 
=& \Re \langle i[H,A] \psi^{(3)}_\alpha \tilde U_\alpha(t) \tilde f, \psi^{(3)}_\alpha \tilde U_\alpha(t) \tilde f \rangle \\
&+ 2 \Re \langle A i[H,\psi^{(3)}_\alpha] \tilde U_\alpha(t) \tilde f, \psi^{(3)}_\alpha \tilde U_\alpha(t) \tilde f \rangle \\
&+ 2 \Re \langle A \psi^{(3)}_\alpha (\partial_t + iH) \tilde U_\alpha(t) \tilde f, \psi^{(3)}_\alpha \tilde U_\alpha(t) \tilde f \rangle 
\end{align*}
and the fundamental theorem of calculus we see that
\begin{multline}\label{heisenberg}
 |\int_{-1}^0 \Re \langle i[H,A] \psi^{(3)}_\alpha \tilde U_\alpha(t) \tilde f, \psi^{(3)}_\alpha \tilde U_\alpha(t) \tilde f \rangle \\
+ 2 \Re \langle A i[H,\psi^{(3)}_\alpha] \tilde U_\alpha(t) \tilde f, \psi^{(3)}_\alpha \tilde U_\alpha(t) \tilde f \rangle\ dt|
\leq C \|\tilde f\|_{L^2(M)}^2.
\end{multline}
To exploit this we will choose $A$ so that both terms on the left-hand side are mostly positive.  In view of \eqref{zetasign} and microlocal analysis heuristics, one expects to be able to do this by taking $A$ and $i[H,A]$ to have positive symbols.  We now interrupt the proof of Proposition~\ref{parametrix-smoothing} to give two lemmas formalizing this intuition.

\begin{lemma}[Positive commutator inequality for parametrix]\label{tech}
Suppose that $A \in \Psisc^{0,0;\rho}(\Mbar)$ and $\sigma(A)$ is non-negative. Then we have
\begin{equation}\label{posit}
\Re \langle A i[H, \psi^{(3)}_\alpha] \tilde U_\alpha(t) \tilde f, \psi^{(3)}_\alpha \tilde U_\alpha(t) \tilde f \rangle \geq -C_A \| \tilde f \|_{L^2(M)}^2,
\end{equation}
uniformly for $-1 \leq t < 0$ and all test functions $\tilde f$. 
\end{lemma}

\begin{proof}
We begin by noting that if $A\in \Psisc^{-\infty,0}(M)$, then $A i[H,\psi^{(3)}]$ is zeroth order,
hence bounded on $L^2(M)$, and the result
follows from Theorem~\ref{sfio-l2}.  Thus we may without loss
of generality take the kernel of $A$ to have support in $\Delta_{\ep}$ with
any desired $\ep>0.$

We now write the identity operator as the sum of two (approximate)
pseudodifferential projections. By \eqref{zetasign}, we
have \begin{equation} \begin{split} &\zeta = \nabla_z d(z, z')^2,
    \ (z,z') \in \supp \nabla \psi_\alpha^{(1)} \times \supp \psi_\alpha^{(0)}
    \\ &\Longrightarrow \sigma(i[H, \psi_\alpha^{(1)}])(z, \zeta) \equiv
    \nabla_{\zeta} \psi_\alpha^{(1)}\leq -\epsilon |\zeta| |\nabla
    \psi_\alpha^{(1)} |,
\end{split}
\label{angle-cond}
\end{equation}
Therefore it is possible to construct a pseudodifferential operator $P = P^\alpha \in \Psisc^{0,0}(\Mbar)$ such that 
\begin{equation}\begin{gathered}
\WF' P^* P \subset\{ \sigma(i[H, \psi_\alpha^{(1)}]) < -\frac{\epsilon}{2} |\zeta| |\nabla \psi_\alpha^{(1)} |\}, \text{ and } \\
\WF' (\Id - P^* P) \cap \{ (z, \zeta) \mid 
\zeta = \nabla_z d(z, z')^2 ,  (z,z')  \in \supp \nabla \psi_\alpha^{(1)}
\times \supp \psi_\alpha^{(0)} \} = \emptyset;
\end{gathered}\end{equation} 
we may further assume that $P^*P$ and $\Id-P^*P$ have kernels supported in
$\Delta_\ep$ for any desired $\ep>0.$

Now we write
\begin{multline}
\langle A i[H, \psi_\alpha^{(1)}] \tilde U(t) \psi_\alpha^{(0)} f, U^\alpha(t) f
\rangle = \langle \psi_\alpha^{(1)} A i[H, \psi_\alpha^{(1)}] P^* P \tilde U(t)
\psi_\alpha^{(0)} f, \tilde U^\alpha(t) \psi_\alpha^{(0)} f \rangle \\ + \langle A
i[H, \psi_\alpha^{(1)}] (\Id - P^* P) \tilde U(t) \psi_\alpha^{(0)} f, U^\alpha(t)
f \rangle.
\label{expansion}\end{multline}
In the first of the terms on the right hand side, $ \psi_\alpha^{(1)} A
  i[H, \psi_\alpha^{(1)}] P^* P$ has a positive symbol. In fact, since the
  symbol of $A$ is negative, we have $A = - F^* F$ up to an operator of
  order $(-1, 1)$. Similarly, $\sigma(-i[H, \psi_\alpha^{(1)}]) >
  (\epsilon/2) |\zeta||\nabla \psi_\alpha^{(1)} |$ on the microsupport of
  $P$, so there is an operator $F_1$ of order $(1/2, 0)$ such that
  $\sigma(F_1) = \sqrt{ \sigma(-i[H, \psi_\alpha^{(1)}])|\nabla
  \psi_\alpha^{(1)} |^{-1}}$ on the microsupport of $P$.  Thus up to
  operators of order zero, $\psi_\alpha^{(1)} A i[H, \psi_\alpha^{(1)}] P^*
  P$ may be written in the form $G^* G$, where $G := \sqrt{\psi_\alpha^{(1)}
  |\nabla \psi_\alpha^{(1)} |} FF_1P$, while the order zero error is
  bounded by $C \| f \|_2^2$ by Theorem~\ref{sfio-l2}. Using
  \eqref{phi-sq}
 we see that $\sqrt{\psi_\alpha^{(1)} |\nabla
  \psi_\alpha^{(1)} |}$ is smooth, so $G$ is an order $(1/2, 1/2)$
  pseudodifferential operator.  The contribution of $G^* G$ to \eqref{posit} is non-negative
and can thus be discarded.

We now consider the $\Id-P^* P$ term.  Heuristically, the contribution of this term should
be extremely small because $i[H,\psi_\alpha^{(3)}] \tilde U_\alpha(t)$ should be
microlocally supported in ``outgoing'' directions (starting from the support of $\psi_\alpha^{(2)}$
and emanating towards the support of $\nabla \psi_\alpha^{(3)}$), whereas $\Id-P^* P$ is microlocally supported
away from the outgoing directions.  However, to make this intuition precise
requires a certain amount of stationary phase arguments\footnote{A slightly different alternative to the arguments given here would be to dyadically decompose the region of integration and apply rescaled versions of the localization principle (see e.g. \cite{stein:large}) to each dyadic component.}.  We begin by
observing that we may replace $A i[H,\psi_\alpha^{(1)}]$
by $i[H,\psi_\alpha^{(1)}] A,$ as the difference is in $\Psisc^{0,2}(M),$ hence
its composition with $\tilde U(t)$ is $L^2$-bounded by
Theorem~\ref{sfio-l2}.  Now we may suppose that the symbol $q$ of $A (\Id -
P^*P)$ is written in right-reduced form. Then the kernel of $ i[H,
  \psi_\alpha^{(1)}] A (\Id - P^* P) \tilde U(t) \psi_\alpha^{(0)}$ is
\begin{equation}
t^{-n/2} \int -i \chi(z,z'') (\nabla \psi_\alpha^{(1)}(z) \cdot \zeta)
e^{i(z-z'') \cdot \zeta} q(z'', \zeta) e^{i\Phi(z'',z')/t} a(z'',z')
\psi_\alpha^{(0)}(z') \, d\zeta \, dg_{z''}; 
\label{ker}\end{equation}
here we have explicitly inserted a cutoff $\chi,$ supported in $\Delta_{2 \ep}$, as
permitted by our support assumptions on the kernels of $A$ and $\Id-P^*P.$
By construction, $q(z'', \zeta)$ vanishes when $\zeta = \nabla_z \Phi(z'',z')/t$,  $(z'',z')  \in \supp \nabla \psi_\alpha^{(1)} \times \supp
\psi_\alpha^{(0)}$.  
Thus we may integrate by parts in $z''$. By Lemma~\ref{psilemma}, the angle
between $\zeta$ and $\nabla_{z''} \Phi(z'',z')/t$ is bounded away from zero,
when $(z'',z') \in \supp \nabla \psi_\alpha^{(1)} \times \supp \psi_\alpha^{(0)}$
and $(z'', \zeta)$ is in the cone support of $q$. Therefore, there is an
$\epsilon' > 0$ such that 
\begin{equation}
\big| t\zeta - \nabla_{z''} \Phi(z'',z') \big| \geq \epsilon' \Big(t |\zeta| + \big| 
\nabla_{z''} \Phi(z'',z') \big| \Big) = \epsilon' \Big( t|\zeta| + d_M(z'',z') \Big) \geq C,
\label{elliptic}\end{equation}
when $z'', z', \zeta$ are as above, since $d_{\overline M}(z'',z')$ (and hence $d_M(z'',z')$)
is bounded away from zero. Let 
\begin{equation*}
v_j(z'', z', \zeta, t) =   \frac{t\zeta_j - \nabla_{z''_j} \Phi(z'',z') }{
\big| t\zeta - \nabla_{z''} \Phi(z'',z') \big|_{g(z'')}^2}.
\end{equation*}
Then for all $k,$ the operator $(i t v^j \nabla_{z''_j})^k$ leaves $e^{i(z-z'') \cdot \zeta}
e^{i\Phi(z'',z')/t}$ invariant. Integrating by parts gives the kernel
\begin{multline}\label{kernel}
t^{-n/2} \int -i \chi(z,z'') (\nabla \psi_\alpha^{(1)}(z) \cdot \zeta)
e^{i(z-z'') \cdot \zeta} e^{i\Phi(z'',z')/t} \\ 
\times ( -it \nabla_{z''_j} v^j)^k \Big( q(z'', \zeta) a(z'',z') \Big)  
 \psi_\alpha^{(0)}(z') \, d\zeta \,
dg_{z''}.
\end{multline}

Note that for $(z',z'')$ in the support of the symbol, we have
$d_{\overline M}(z',z'') >\delta>0$ for some $\delta,$ hence
$$
d(z',z'') \geq C \min(\ang{z'},\ang{z''})
$$
for some $C.$  Now this implies that
$$
\abs{t v^j}^{-1} \leq  C \min(\abs{\zeta}^{-1}, t (\ang{z'} + \ang{z''})^{-1}),
$$
hence, for any $N,$ taking $k$ sufficiently large (and using symbol
estimates on $q$ and $a$) we may estimate those terms in \eqref{kernel} in
which no derivative hits a $v^j$ by a multiple of
$$
t^N \int \chi(z,z'') \Big(\ang{z'}+\ang{z''}\Big)^{-N} \, dg_{z''};
$$ the $\zeta$ integral has disappeared as we have taken $k$ sufficiently
large that the integrand is $L^1$ in $\zeta$ with bounds in $z',z''$ as
above.  This term is $L^2$-bounded, uniformly in $t,$ as it is
Hilbert-Schmidt for $N$ large.

To estimate terms in which derivatives fall upon $v_i$'s, note that
by
Lemma \ref{metric-smooth}, we may estimate, for $k\geq 1,$
$$|\nabla_{z''}^k (t\zeta-\Phi(z'',z'))| \leq C ((r'')^{2-k} + r'
(r'')^{1-k});$$
on the support of $a,$ this is bounded by $C ((r')^{1-k} +
(r'')^{1-k}).$  This allows us to conclude that
$$
\abs{(t \nabla_{z''})^k v_j} \leq C \min (\abs{\zeta}^{-k}, t^k (r'+r'')^{-k}),
$$
and the estimates on terms involving derivatives follow as above.

\end{proof}

\begin{lemma}[Construction of commutant]\label{A-ok}  There exists a scattering pseudo-differ\-ential operator $A \in \Psisc^{0,0;\rho}(\Mbar)$ such that $\sigma(A) \geq 0$ and
\begin{equation}\label{aha}
 i[H,A] \geq (1+H)^{1/4} \langle z \rangle^{-1+\rho} (1 + H)^{1/4} + B^* B
 - O(\Psisc^{0,0;\rho}(\Mbar)),
\end{equation}
where we use $A \geq B$ to denote the assertion that $A-B$ is positive
semi-definite and $O(\Psisc^{0,0;\rho}(\Mbar))$ denotes some term in $\Psisc^{0,0;\rho}(\Mbar).$
\end{lemma}

\begin{proof}  This argument is essentially contained in \cite[Lemmas 10.4, 10.5]{HTW}, and we give only a sketch here.
We first observe that the claim $\sigma(A) \geq 0$ is automatic, since if
$\sigma(A)$ is not non-negative we may simply add a large constant to $A$
(taking advantage of the symbol bounds $A \in \Psisc^{0,0;\rho}(\Mbar)$;
note that this does not affect the commutator $i[H,A]$).  Thus we can ignore
the constraint $\sigma(A) \geq 0$.

We still have to find an operator $A \in \Psisc^{0,0;\rho}(\Mbar)$ which obeys the inequality
\eqref{aha}.  In practice we accomplish this in several separate pieces,
which we will add together at the end.

To start with, we exploit the non-trapping hypothesis on the manifold $M$
in the usual manner (see \cite{cks}, \cite{doi}) to construct, for any
given compactly supported bump function $\varphi$, an operator $A_0 \in
\Psisc^{-1,0;\rho}(\Mbar)$ such that $i[H,A_0] \geq \varphi^2 +
O(\Psisc^{-1,0;\rho}(\Mbar)$.  Multiplying this on the left and right by
$(1+H)^{1/4}$ (which is self-adjoint and commutes with $H$) we obtain a
symbol $(1+H)^{1/4} A_0 (1+H)^{1/4} \in \Psisc^{0,0;\rho}(\Mbar)$ such that
$$ i[H,(1+H)^{1/4} A_0 (1+H)^{1/4}] \geq (1+H)^{1/4} \varphi^2 (1+H)^{1/4}
+ O(\Psisc^{0,0;\rho}(\Mbar)).$$ This gives a commutant with the desired
properties over any compact subset of $\Mbar$ (note that $B$ is of order
$1/2$).  Thus it will suffice to obtain the inequality \eqref{aha} in the
asymptotic region $\langle z \rangle \gg 1$, since the claim will then
follow by adding a sufficiently large multiple of $(1+H)^{1/4} A_0
(1+H)^{1/4}$ to $A$.

We now test the radial vector field $A_1 := \eta^2 \partial_r \in
\Psisc^{1,0;\rho}(\Mbar)$, where $r$ is the radial co-ordinate in the
scattering region and $\eta$ is a smooth cutoff to the asymptotic region
$\langle z \rangle \gg 1$.  A computation (see \cite[Lemma 10.4]{HTW})
shows that
$$ i[H,A_1] \geq c (\langle z \rangle^{-1/2} \nabb_j \eta)^* h^{jk} (\langle z \rangle^{-1/2} \nabb_k \eta) - \varphi O(\Psisc^{2,0;\rho}(\Mbar)) \varphi
- O(\Psisc^{1,0;\rho}(\Mbar))$$
for some compactly supported bump function $\varphi$ and some $c>0$.  A similar computation (again see \cite[Lemma 10.4]{HTW})
shows that if we let $A_2 := - \eta^2 \langle z \rangle^{-2\eps} \partial_r \in \Psisc^{1,0;\rho}(\Mbar)$ for some small $\eps > 0$ then
\begin{multline*} i[H,A_2] \geq C_\eps (\langle z \rangle^{-1/2-\eps} \nabla)^* (\langle z \rangle^{-1/2-\eps} \nabla)
- (x^{1/2} \nabb_j \eta)^* O(\Psisc^{0,0;\rho}(\Mbar)) (x^{1/2} \nabb_k \eta) \\ - 
\varphi O(\Psisc^{2,0;\rho}(\Mbar)) \varphi - O(\Psisc^{1,0;\rho}(\Mbar)
\end{multline*}
for some $C_\rho > 0$.  Adding a suitable combination of the two, we can find $A_3 \in \Psisc^{1,0;\rho}(\Mbar)$
such that
$$ i[H,A_3] \geq (\langle z \rangle^{-1/2-\eps} \nabla)^* (\langle z \rangle^{-1/2-\eps} \nabla)
- \varphi O(\Psisc^{2,0;\rho}(\Mbar)) \varphi - O(\Psisc^{1,0;\rho}(\Mbar)).$$
Multiplying on left and right by $(1+H)^{-1/4}$, and then adding a large multiple of $(1+H)^{1/4} A_0 (1+H)^{1/4}$, we
can find $A_4 \in \Psisc^{0,0;\rho}(\Mbar)$ such that
$$ i[H,A_4] \geq (1+H)^{1/4} \langle z \rangle^{-1-2\eps} (1 + H)^{1/4} 
- O(\Psisc^{0,0;\rho}(\Mbar))$$
which gives part of \eqref{aha} (with $\eps := \rho/2$).  Finally, 
from \cite[Lemma 10.5]{HTW} we can obtain $A_5 \in \Psisc^{0,0;\rho}(\Mbar)$
such that
$$ i[H,A_5] \geq B^* B - O(\Psisc^{1,1+2\eps;\rho}(\Mbar)) - 
O(\Psisc^{0,0;\rho}(\Mbar)))$$
for some $\eps > 0$.  By adding a large multiple of $A_4$ to $A_5$ we can eliminate the 
$O(\Psisc^{1,1+2\eps;\rho}(\Mbar))$ error and obtain \eqref{aha} as desired.
\end{proof}

\emph{Proof of Proposition~\ref{parametrix-smoothing} (completion).} We choose $A$ as in Lemma~\ref{A-ok} and use this in \eqref{heisenberg}. Using \eqref{aha} and \eqref{posit}, we see that \eqref{ls-variant} holds, which completes the proof.
\end{proof}


\section{Proof of Theorem~\ref{main}}\label{sec:proof}

We can now deduce Theorem \ref{main} from the machinery developed in the previous sections (the proofs of some of which are deferred to Appendices).  
We fix $M$, $q$, $r$, $u$, and allow all constants $C$ to depend on $M$, $n$, $q$, $r$.
Write $f := u(0)$, thus $u := e^{-itH} f$. 

It suffices by \eqref{partition} and the triangle inequality to establish the estimate
$$ \| \psi^{(0)}_\alpha u(t) \|_{L^q([0,1]; L^r(M))} \leq C_{n,q,r,M} \| f \|_{L^2(M)} $$
for each $\alpha \in A$.
By \eqref{identity}, \eqref{nls} and the fundamental theorem of calculus we have the Duhamel identity
\begin{align*}
\psi^{(0)}_\alpha u(t) &= U_\alpha(t) \psi^{(0)}_\alpha u(0) +
\int_{0}^t \frac{d}{ds} (U_\alpha(t-s) \psi^{(0)}_\alpha u(s))\ ds\\
&= U_\alpha(t) \psi^{(0)}_\alpha f +
\int_{0}^t -U'_\alpha(t-s) \psi^{(0)}_\alpha u(s) - U_\alpha(t-s) \psi^{(0)}_\alpha iH u(s)\ ds \\
= U_\alpha(t) &\psi^{(0)}_\alpha f +
\int_{0}^t -(U'_\alpha + iU_\alpha H)(t-s) \psi^{(0)}_\alpha u(s) 
+ U_\alpha(t-s) [iH,\psi^{(0)}_\alpha] u(s)\ ds.
\end{align*}
We can multiply both sides by $\psi^{(1)}_\alpha$ and use \eqref{nesting} to obtain
\begin{multline*}
\psi^{(0)}_\alpha u(t) = 
\psi^{(1)}_\alpha U_\alpha(t) \psi^{(0)}_\alpha f + \\
\int_{0}^t - \psi^{(1)}_\alpha (U'_\alpha + iU_\alpha H)(t-s) \psi^{(0)}_\alpha u(s) 
+ \psi^{(1)}_\alpha U_\alpha(t-s) [iH,\psi^{(0)}_\alpha] u(s)\ ds.
\end{multline*}
Now we take $L^q_t L^r_x$ norms of both sides and use the triangle inequality.  The first term
will be acceptable from \eqref{S-endpoint}.  From Lemma \ref{LSIO-operator} we see
that $\psi^{(1)}_\alpha (U'_\alpha + iU_\alpha H)(t-s) \psi^{(0)}_\alpha$ is an LSIO, and
so the second term will also be acceptable from \eqref{S-endpoint} and Minkowski's inequality (noting that
$u(s)$ has $L^2$ bounded by a multiple of that of $u(0)=f$).  It thus will suffice to show that
$$ \| \int_{0}^t \psi^{(1)}_\alpha U_\alpha(t-s) [iH,\psi^{(0)}_\alpha] u(s)\ ds \|_{L^q_t L^r_x([0,1] \times M)}
\leq C \| f \|_{L^2(M)}.$$
From Lemma \ref{ls-lemma} it will thus suffice to show that
$$ \| \int_{0}^t \psi^{(1)}_\alpha U_\alpha(t-s) [iH,\psi^{(0)}_\alpha] v(s)\ ds \|_{L^q_t L^r_x([0,1] \times M)}
\leq C \| v \|_{L^2_t([0,1];X)}$$
for all Schwartz functions $v$.  
But by the Christ-Kiselev lemma (see Lemma \ref{ck-lemma}
in Appendix \S \ref{sec:lemmas}) and the non-endpoint hypothesis $q>2$, it will suffice to show that
\begin{equation}\label{bigmess}
\|\int_0^{t_0} \psi^{(1)}_\alpha U_\alpha(t-s) i[H,\psi^{(0)}_\alpha] v(s)\ ds \|_{L^q([t_0,1]; L^r(M))} \leq C \| v \|_{L^2_t([0,1]; X)}
\end{equation}
for all $0 \leq t_0 \leq 1$.

Fix $t_0$.  We shall take advantage of the approximate semigroup properties of $U_\alpha$ to approximately
factorize the expression in \eqref{bigmess}, thus decoupling the hybrid Strichartz/local smoothing estimate
\eqref{bigmess} into a Strichartz estimate component for the parametrix (which we can deal with by
\eqref{S-endpoint}) and a local smoothing component for the parametrix (which we can deal with
by Lemma \ref{tech}).  More precisely, we decompose
$$ U_\alpha(t-s) = U_\alpha(t-t_0) U_\alpha(t_0-s) + 
(U_\alpha(t-s) - U_\alpha(t-t_0) U_\alpha(t_0-s)).$$
Let us first consider the contribution of the error term $U_\alpha(t-s) - U_\alpha(t-t_0) U_\alpha(t_0-s)$.
From the definition of $U_\alpha$ and
Theorem \ref{sfio-concatenate} we see that this error term is an LSIO of the form $S_c(t-s)$, where $c = c_{t-t_0,t_0-s}$ is
defined by
$$ c(z,z') := \psi^{(2)}_\alpha(z) \psi^{(2)}_\alpha(z') (1 -
\psi^{(2)}_\alpha(w)^2) - \frac{t-s}{\langle z \rangle + \langle z'
  \rangle} e(z,z')=: c_1 + c_2$$ for some product symbol $e = e_{t-t_0,t_0-s}$ of order
$(0,0)$, and $w := \gamma_{z \to z'}(\frac{t-t_0}{t-s})$.  From Lemma
\ref{psilemma} we know that $\psi_\alpha^{(1)} c_1  i [H,\psi_\alpha^{(0)}]
=0.$  By \eqref{S-def} we thus have
\begin{align*}
& \psi^{(1)}_\alpha S_c(t-s) 
i[H,\psi^{(0)}_\alpha] f(z) \\
&:= \frac{1}{(2\pi i (t-s))^{n/2}} \int_M e^{i \frac{\Phi(z,z')}{t-s}} 
\frac{t-s}{\langle z \rangle + \langle z' \rangle} \tilde e(z,z') 
i[H,\psi^{(0)}_\alpha] f(z')\ dg(z').
\end{align*}
If we now integrate the first-order operator $i[H,\psi^{(0)}_\alpha]$ by parts, we thus see that
the operator $\psi^{(1)}_\alpha S_c(t-s) i[H,\psi^{(0)}_\alpha] = S_d(t-s)$ for some product symbol
$d = d_{t-s}$ of order $(0,0)$ (the damping factor $\frac{t}{\langle z \rangle + \langle z' \rangle}$ can 
absorb any term that arises when the $z'$ derivative hits $\Phi(z,z')/t$, thanks for instance to Lemma~\ref{metric-smooth}).  
Furthermore the symbol bounds of $d$ are uniform in $s,t$ in the range $0 \leq s \leq t \leq 1$.
In particular, we see from \eqref{S-endpoint} that for any $0 < s < 1$ we have
$$ \| \psi^{(1)}_\alpha S_{c_{t-t_0,t_0-s}}(t-s) i[H,\psi^{(0)}_\alpha] v(s) \|_{L^q_t L^r_x([0,1] \times M)}
\leq C \| v(s) \|_{L^2(M)}$$
The contribution of the error term to \eqref{bigmess} is thus acceptable from Minkowski's inequality since
the $X$ norm controls the $L^2(M)$ norm by \eqref{X-def}.

It remains to consider the contribution of the main term $U_\alpha(t-t_0) U_\alpha(t_0-s)$ to \eqref{bigmess}.
We can factorize this as
$$ \| \psi^{(1)}_\alpha U_\alpha(t-t_0) \int_0^{t_0} U_\alpha(t_0-s) i[H,\psi^{(0)}_\alpha] v(s)\ ds \|_{L^q([t_0,1]; L^r(M))} \leq C \| v \|_{L^2_t([0,1]; X)}$$
for all $0 \leq t_0 \leq 1$.  From \eqref{S-endpoint} it will suffice to show that
$$ \| \int_0^{t_0} U_\alpha(t_0-s) i[H,\psi^{(0)}_\alpha] v(s)\ ds \|_{L^2(M)} \leq C \| v \|_{L^2_t([0,1]; X)}$$
We test this inequality against an arbitrary test function $\tilde f$, and see that it suffices by duality and
\eqref{symmetry} to show that
$$ \Big |\int_0^{t_0} \langle U_\alpha(s-t_0) \tilde f, i[H,\psi^{(0)}_\alpha] v(s) \rangle\ ds \Big|
\leq C \| v \|_{L^2_t([0,1];X)} \| \tilde f \|_{L^2(M)}.$$
This follows from Proposition~\ref{parametrix-smoothing} and Lemma~\ref{X-split}. This completes the proof of Theorem \ref{main}.
\qed



\section{Remarks}\label{remarks:sec}

We comment briefly on some standard and automatic extensions of the
estimate in Theorem \ref{main}.  The time interval $[0,1]$ can easily be
replaced with any other finite interval by time translation invariance and
concatenating several time intervals together (and exploiting the fact that
free propagators $e^{-itH}$ preserve the $L^2(M)$ norm), although at present we
are unable to prevent the constants in the estimates from growing with the length of
the interval and so cannot prove a global-in-time Strichartz estimate (in analogy with \eqref{strich})
though we conjecture that such estimates are true.  Also one can add
a bounded lower order potential term $V$ to $H$ without affecting this
estimate by using Duhamel's formula and exploiting the time localization.
Also one can extend the Strichartz estimate to allow for the presence of a
forcing term by use of $TT^*$ arguments and the Christ-Kiselev lemma (see
e.g. \cite{tao}).    The condition that $u$ and $F$ are Schwartz
can be removed by the usual limiting arguments.  Again, we can use Sobolev
embedding and commute the equation with fractional powers of $H$ to extend
these results to other Sobolev spaces.


Somewhat trickier would be to obtain the endpoint $q=2$ (in dimensions $n \geq 3$, of course).  The Christ-Kiselev
lemma (Lemma \ref{ck-lemma}) is no longer applicable, so one would have to prove the hybrid local smoothing/Strichartz
estimate \eqref{bigmess} by other means.  We will pursue this matter elsewhere.

As is well known, Strichartz estimates can be used to obtain local well-posedness results for semilinear Schr\"odinger equations such as
$$ u_t = -iHu + \lambda |u|^{p-1} u; \quad u(0) = f$$ where $p>1$ and
$\lambda = \pm 1$ is fixed.  Indeed, the Strichartz estimates given above
imply that this equation is locally well-posed in $H^s(M)$ whenever $s \geq
\max(0, \frac{n}{2} - \frac{2}{p-1})$, and (in order to ensure $|u|^{p-1}
u$ is sufficiently regular) the exponent $p$ is either odd, or greater than
$\lfloor s \rfloor + 1$; this can be achieved by applying the arguments of
Cazenave and Weissler \cite{cwI} with Euclidean space replaced by the
manifold $M$.  (For $p$ not an odd integer, one needs to develop a
fractional Liebniz rule for Sobolev spaces on manifolds, but this can be
done by working in local co-ordinates, in which the Sobolev norms are
comparable to their Euclidean counterparts).  We omit the details.  For
similar reasons one can adapt the arguments in \cite{kato} to obtain global
well-posedness for this equation in the energy norm $\|u(t)\|_{H^1}$ when
$p < 1+4/(n-2)$ (in the defocusing case $\lambda = -1$) or when $p < 1+4/n$
(in the focusing case $\lambda = +1$), and to obtain the respective
endpoint values of $p$ when the $H^1$ norm is sufficiently small. Note that
none of these results require the endpoint Strichartz estimate (which was
only obtained later, in \cite{tao:keel}).  On the other hand, these
Strichartz estimates are not strong enough to obtain a scattering theory
for any of these equations (as is done for instance in \cite{caz}), because
they are only local in time.  To obtain global-in-time Strichartz estimates
would require one to extend local smoothing estimates such as those in
Lemma \ref{ls-lemma} to also be global in time, and would also require one
to treat lower order terms in the above arguments more carefully.



\section{Appendix: Christ-Kiselev lemma}\label{sec:lemmas}

In this section we prove a minor variant of the Christ-Kiselev lemma \cite{ck}.

\begin{lemma}[Christ-Kiselev Lemma \cite{ck}]\label{ck-lemma}  Let $X$, $Y$ be Banach spaces, and for all $s,t \in \R$
let $K(s,t): X \to Y$ be an operator-valued kernel from $X$ to $Y$.  Suppose we have the estimate
\begin{equation}\label{ck-hyp}
 \| \int_{s < t_0} K(s,t) f(s)\ ds \|_{L^q([t_0,\infty);Y)} \leq A \| f \|_{L^p(\R;X)}
\end{equation}
for some $A > 0$ and $1 \leq p < q \leq \infty$, and all $t_0 \in \R$ and
$f \in L^p((-\infty,t_0);X)$ and $t_0 \in \R$.  Then we have
$$ \| \int_{s < t} K(s,t) f(s)\ ds \|_{L^q(\R;Y)} \leq C_{p,q} A \| f \|_{L^p(\R;X)}.$$
\end{lemma}

\begin{remarks}
Note that the hypothesis \eqref{ck-hyp} is in particular satisfied if we have the global estimate
\begin{equation}\label{Apq}
 \| \int_\R K(s,t) f(s)\ ds \|_{L^q(\R;Y)} \leq A \| f \|_{L^p(\R;X)},
\end{equation}
since one can simply restrict $f$ to $(-\infty,t_0)$ and restrict the quantity on the left-hand side
to the region $t \geq t_0$.  For sake of completeness we shall prove the above variant of the Christ-Kiselev lemma
in Appendix \S \ref{sec:lemmas}.  This Lemma is very convenient in the non-endpoint setting for eliminating
the causality constraint $s < t$ in the operators that one wishes to bound, and has been utilized in 
previous work on the Schr\"odinger equation, see \cite{tao}, \cite{st}, \cite{burq}.  Clearly we may 
replace the time domain $\R$ by any smaller interval such as $[0,1]$ and obtain a similar result in which 
all time variables $s, t, t_0$ are restricted to this interval.  The estimate unfortunately fails at the endpoint
$p=q$; see \cite{tao} for an instance of this failure in a context of obtaining Strichartz estimates for the Schr\"odinger equation.
\end{remarks}

\begin{proof}[Proof of Lemma~\ref{ck-lemma}] 
This is a minor modification of the proof of \cite{tao}, Lemma 3.1 or \cite{ss}, Lemma 3.1,
but for sake of completeness we give the full proof (as our hypotheses are slightly different).  By a limiting argument
we may assume that the cumulative distribution function $F(t_0) := \int_{-\infty}^{t_0} \| f(t)\|_X^p\ dt$ is continuous
and strictly increasing in $t_0$.  We may also normalize $\|f\|_{L^p(\R;X)} = 1$, hence $F$ is a bijection from $\R$ to
$[0,1]$.

We impose the usual dyadic grid on $[0,1]$.  Given two dyadic intervals\footnote{We will ignore sets of measure zero and so we will not care whether these intervals are open or closed.} $I, J$ strictly contained in $[0,1]$, we
define the relation $I \sim J$ if $I$ and $J$ have the same length, $I$ is to the left of $J$, and that $I$ and $J$ are not
adjacent but have adjacent parents; let us write $S_n$ for the set of pairs $(I,J)$ with $|I| = |J| = 2^{-n}$ and $I \sim J$.  Also, for almost every $0 < x < y < 1$
there is exactly one pair of intervals $I, J$ such that $x \in I$, $y \in J$, and $I \sim J$.  Applying $F^{-1}$ to
this fact, we thus have (cf. \cite{tao}, \cite{ss})
$$ \int_{s < t} K(s,t) f(s)\ ds = \sum_{I,J: I \sim J} 1_{F^{-1}(J)}
\int_{F^{-1}(I)} K(s,t) f(s)\ ds.$$ We group the intervals based on their
dyadic length, take $L^q(\R;Y)$ norms, and apply the triangle inequality to
obtain
$$ \| \int_{s < t} K(s,t) f(s)\ ds \|_{L^q(\R;Y)} \leq \sum_{n=1}^\infty \|
\sum_{(I,J) \in S_n} 1_{F^{-1}(J)} \int_{F^{-1}(I)} K(s,t) f(s)\ ds
\|_{L^q(\R;Y)}.
$$ Observe that for fixed $n$ and $J$, there are at most two intervals $I$
for which $I \sim J$.
Thus the inner summands have finitely overlapping supports, and we can estimate
$$ \| \int\limits_{s < t} K(s,t) f(s)\ ds \|_{L^q(\R;Y)} \leq C_q \sum_{n=1}^\infty 
\Big( \!\! \sum_{(I,J) \in S_n} \!\!\!  \big\|   1_{F^{-1}(J)} \!\!\!\!\! \int\limits_{F^{-1}(I)} \!\!\!\!\! K(s,t) f(s)\ ds \big\|_{L^q(\R;Y)}^q \Big)^{1/q} \! .$$
Applying the hypothesis \eqref{ck-hyp} (with $t_0$ chosen to be some intermediate time between $F^{-1}(I)$ and $F^{-1}(J)$)
we obtain
$$ \| \int_{s < t} K(s,t) f(s)\ ds \|_{L^q(\R;Y)} \leq C_q A \sum_{n=1}^\infty 
(\sum_{(I,J) \in S_n} \| f(s) \|_{L^p(F^{-1}(J);X)}^q)^{1/q}.$$
From the definition of $F$ we have $\| f(s) \|_{L^p(F^{-1}(J);X)} = |J|^{1/p} = 2^{-n/p}$, thus we have
$$ \| \int_{s < t} K(s,t) f(s)\ ds \|_{L^q(\R;Y)} \leq C A \sum_{n=1}^\infty 
(\sum_{(I,J) \in S_n} 2^{-nq/p})^{1/q}
= C_q A \sum_{n=1}^\infty (2^n 2^{-nq/p})^{1/q} $$
which is less than infinity, so the claim follows.
\end{proof}

\section{Appendix: Phase nondegeneracy estimates}\label{sec:geometry}

In this section we provide the proof of Lemmas \ref{glc}, \ref{grad1} and
\ref{geodesic-lemma}.  

\subsection{Proof of Lemma \ref{glc}}

We begin with the proof of Lemma \ref{glc}.  If we restrict attention to any fixed compact subset $K_0$ of $M$,
then the claim about the existence of a unique geodesic $\gamma_{z \to z'}$ follows from local Riemannian geometry
by taking $\epsilon$ (and hence all the $o(1)$ factors) 
sufficiently small.  Thus by letting $K_0$ grow very slowly in $\epsilon$, it suffices to work in the asymptotic 
region $\{ x \leq o(1)\}$.  Now let us consider what happens in this region in the perfectly conic case $g = g_\conic$.

Consider any metric of the form $f(r) dr^2 + g(r) h(y, dy)$ on $M$, where $h$ is a metric on a cross-section manifold $N$. 
Consider a geodesic arc $\theta \mapsto \gamma(\theta)$ in the cross-section $(N, h)$.  
It is not difficult to check that the two-dimensional cone $\Cgamma = \{ (r,\gamma(\theta)) \} \subset M$ over $\gamma$ is totally geodesic in $M$. Let us write $d\theta^2$ for the metric on the geodesic $\gamma$. For the conic metric, $dr^2 + r^2 h$, the induced metric on $\Cgamma$ is exactly equal to the Euclidean metric $dr^2 +
r^2 d\theta^2$ on the plane $\R^2$ in polar co-ordinates.  Thus
the geometry of geodesics on the cone $\Cgamma$ is
identical to that on the Euclidean plane $\R^2$.  Conversely, every
geodesic $\lambda$ in $M$ is contained in a cone of the form $\Cgamma$, since one can easily locate such a cone for which
$\lambda$ is tangent to at least second order, and then the claim follows
from uniqueness of the second-order geodesic flow ODE.  Note that we can
now conclude uniqueness of geodesics and the cosine rule
\eqref{sharp-cosine} for perfectly conic manifolds, since this rule is
known to already hold in the Euclidean plane $\R^2$.  For similar reasons
we also obtain \eqref{cycle} in the perfectly conic case.

Similarly, the cone $\Cgamma$ is totally geodesic for the incomplete metric $\overline{g}$, and on $\Cgamma$ the  metric 
induced by $\overline{g}$ takes the form 
$$ 
\overline{g} = \frac{dx^2}{(1 + x^2)^2} + \frac{d\theta^2}{1+x^2}
= \frac{dr^2}{(1+r^2)^2} + \frac{r^2}{1+r^2} d\theta^2.$$

Following \cite{melrose}, we let $S^2_+$ denote the upper hemisphere of
the standard unit sphere in $\RR^3$ and $\RC: \RR^2 \to S^2_+$ denote the
``radial compactification'' (or ``inverse gnomonic projection'')
$$
z \mapsto \Big( \frac{z}{\sqrt{1+\abs{z}^2}}, \frac{1}{\sqrt{1+\abs{z}^2}}\Big)
$$
(see figure).
Letting $g_{S^2_+}$ denote the standard metric on the hemisphere, we then compute
$$
\RC^*(g_{S^2_+}) = \overline{g}
$$ The map $\RC$ identifies lines on the plane with great semi-circles on
 the hemisphere, and thus maps geodesics to geodesics.  From this we
 conclude that $g_\conic$ and $\overline{g}$ do indeed give the same
 geodesics as claimed.
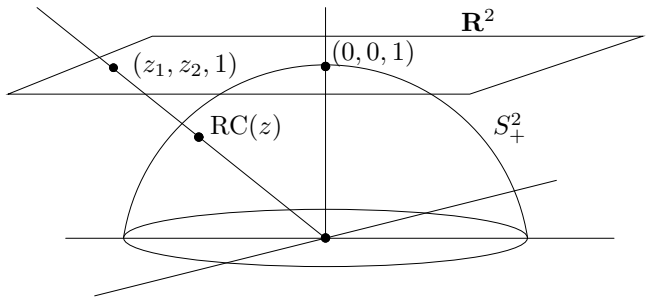
\begin{figure}
\setlength{\unitlength}{0.0005in}
\begingroup\makeatletter\ifx\SetFigFont\undefined%
\gdef\SetFigFont#1#2#3#4#5{%
  \reset@font\fontsize{#1}{#2pt}%
  \fontfamily{#3}\fontseries{#4}\fontshape{#5}%
  \selectfont}%
\fi\endgroup%
{\renewcommand{\dashlinestretch}{30}
\begin{picture}(6624,3039)(0,-10)
\put(3312.000,287.000){\arc{4250.000}{3.2951}{6.1296}}
\put(3312,612){\ellipse{4200}{600}}
\put(1107,2382){\blacken\ellipse{84}{84}}
\put(1107,2382){\ellipse{84}{84}}
\put(3312,2397){\blacken\ellipse{90}{90}}
\put(3312,2397){\ellipse{90}{90}}
\put(3312,612){\blacken\ellipse{90}{90}}
\put(3312,612){\ellipse{90}{90}}
\put(1992,1662){\blacken\ellipse{90}{90}}
\put(1992,1662){\ellipse{90}{90}}
\path(1512,2712)(6612,2712)(4812,2112)
        (12,2112)(1512,2712)
\path(6312,612)(612,612)
\path(5712,1212)(912,12)
\path(3312,612)(312,3012)
\path(3312,3012)(3312,612)
\put(1302,2337){$(z_1,z_2,1)$}
\put(2112,1692){$\mathrm{RC}(z)$}
\put(3372,2460){$(0,0,1)$}
\put(4722,2772){$\RR^2$}
\put(5052,1692){$S^2_+$}
\end{picture}
}
\caption{The radial compactification map}
\end{figure}


It remains to prove the claims of Lemma in the asymptotically conic case in the asymptotic region 
where $x,x' \leq o(1)$ and $d_{\overline M}(z,z') = o(1)$.  We begin with the claim that the geodesic $\gamma_{z \to z'}$ is unique and has diameter $o(1)$ in the $d_{\overline M}$ metric.  We repeat the arguments from
\cite[Lemma 9.1]{HTW}.  We observe the crude estimate
\begin{equation}\label{dm-crude} d_M(z,z') = |r-r'| + o( \min(r,r') )
\end{equation}
for $(z,z')$ in the above region; the lower bound follows from the trivial
comparison $g \geq dr^2$, and the upper bound follows (when $r \leq r'$,
say) by considering the path from $z'$ to $z$ formed by taking the radial
ray from $z' = (r',y')$ to $(r,y')$ and then the angular arc connecting
$(r,y')$ to $(r,y)$.  This inequality, combined with the triangle
inequality, already shows that any length-minimizing geodesic in $M$
connecting $z$ to $z'$ must stay within the asymptotic region $\{ x \leq
o(1) \}$, as a path which strays too close to the compact region $K_0$ will
necessarily have longer length than $|r-r'| + o(\min(r,r'))$.  In this
asymptotic region we have the comparision $g \geq \tilde g_\conic$, where
$\tilde g_\conic$ is the same perfectly conic metric as $g_\conic$ but with
$h$ replaced by $(1-o(1)) h$.  Thus $d_M(z,z') \geq \tilde d_\conic(z,z')$
in this region.  From this we also see that the $d_M$-length-minimizing
geodesic must have diameter at most $o(1)$ in the $d_{\overline M}$ metric,
since again a path which wanders further than $o(1)$ away from $z$ and $z'$
in this metric will require too much length in $\tilde d_\conic$ (as can be
seen from \eqref{dm-crude} and \eqref{sharp-cosine}, and a quantitative
version of the triangle inequality in the Euclidean plane $\R^2$; see
\cite[Lemma 9.1]{HTW}) and hence in $d_M$.

To prove the uniqueness of length-minimizing geodesics, we must show that given two points $z, z'$ whose distance is $o(1)$ in the 
compactified metric $d_{\overline M}$, there is only one geodesic between them which remains at distance $o(1)$ from $z$.  For short-range metrics this
is part of \cite[Proposition 9.2]{HTW}; for long-range metrics, we prove this in Corollary~\ref{uniq}.

Now we prove \eqref{cycle} for asymptotically conic metrics. For this we use the comparison $g = (1 + o(1)) g_\conic$ in the asymptotic region to obtain the estimate\footnote{This is 
still quite crude.  For a more precise estimate in the short-range case, see \cite[Proposition 9.4]{HTW}.}
\begin{equation}\label{dm-compare} d_M(z,z') = (1 + o(1)) d_\conic(z,z').
\end{equation}
Now set $Z'' := \gamma_{z \to z'}(\theta)$, and $Z''_\conic
:= \gamma_{z \to z',\conic}(\theta)$, where $\gamma_{z \to z',\conic}$ is the geodesic from $z$ to $z'$ 
in the perfectly conic metric $g_\conic$.  From \eqref{dm-compare} we see that
$$ d_\conic(z, Z'') = (1 + o(1)) \theta d_\conic(z,z''); \quad d_\conic(z', Z'') = (1 + o(1)) (1-\theta) d_\conic(z,z''),$$
while of course
$$ d_\conic(z, Z''_\conic) = \theta d_\conic(z,z''); d_\conic(z', Z''_\conic) = (1-\theta) d_\conic(z,z'').$$
From the triangle inequality we then see that
$$ d_\conic(Z'', Z''_\conic) = o(1) \min(\theta,1-\theta) d_\conic(z,z'')
= o( \ang{Z''_\conic} )$$
(by \eqref{cycle} in the perfectly conic case), and we can now extend \eqref{cycle} to the asymptotically conic case
by the triangle inequality. 
\endprf

\subsection{Proof of Lemma \ref{grad1}}

Now we prove Lemma \ref{grad1}.  We begin by assuming that $z, z', z''$ all lie in a fixed compact set $K_0 \subset M$.  Then the estimates \eqref{det-est}, \eqref{grad}
are immediate from normal coordinates around $z''$
(with constants depending on $K_0$ of course). 
Indeed, for $z, z'$ close enough to $z''$, $z, z', z'' \in
K_0$, the function $\Phi(z,z'')$ agrees with the Euclidean function
$|z-z''|^2/2$ up to terms of third order in $z-z''$; similarly if we
replace $z$ by $z'$.  Also in these normal co-ordinates $|z-z'|$ is
comparable to $d_M(z,z')$.  Since the claims \eqref{det-est}, \eqref{grad} are easily verified in the Euclidean case, 
the claim follows if $z, z', z''$ are close enough.  The estimate \eqref{gradw} is proven similarly, but using
Fermi normal co-ordinates around the geodesic $\gamma_{z \to z'}$ instead of normal co-ordinates around $z''$;
we omit the details.

In light of the above statements, we see that in order to prove Lemma \ref{grad1}, 
we can now fix attention to the asymptotic region $x = o(1)$
(by choosing $K_0$ to be growing slowly as $\epsilon \to 0$.  As in the proof of Lemma \ref{glc}, we begin this
task by considering the case of perfectly conic metrics $g = g_\conic$; the general case turns out to be treatable
by an easy perturbation argument (using Lemma \ref{metric-smooth}) which we give at the end of this argument.

First we prove \eqref{det-est} for perfectly conic metrics. 
From \eqref{sharp-cosine} (and Gauss's Lemma for $\partial M$) 
we see that the matrix $-\nabla_{z''} \nabla_z \Phi(z,z'')$ 
of partial derivatives in \eqref{det-est} (with respect to a standard orthonormal frame associated to the co-ordinates $x, y$) takes the form 
\begin{equation}
-\nabla_{z''} \nabla_z \Phi(z,z'') = 
\begin{pmatrix}
\cos d_{\partial M}(y, y'') & \dots & \sin d_{\partial M}(y,y'') \hat \mu & \dots
\ \ \\ \vdots & & & \\ \sin d_{\partial M}(y,y'') \hat \mu'' & & \nabla_y
\nabla_{y''} \cos d_{\partial M}(y,y'') & \\ \vdots & & &
\end{pmatrix}
\label{zz''derivs}\end{equation}
where $\hat \mu$ is the unit vector at $y''$ parallel to the geodesic from
$y''$ to $y$, and similarly for $\hat \mu'$.
The $(n-1) \times (n-1)$ lower right submatrix is close to the identity
since $\cos d_{\partial M} = 1 - d_{\partial M}^2/2 + O(d_{\partial M}^4)$,
and $d_{\partial M} = d_{\partial M}(y,y'')$ is assumed to be small, namely
$o(1)$. Thus this matrix is the identity plus $o(1)$, and its
determinant is $1 + o(1)$, and hence we have \eqref{det-est} as desired.

Next, we prove \eqref{grad} for perfectly conic metrics. We first do this assuming that $z$ and $z'$ are relatively close
in the uncompactified metric $d_M$, say $d_M(z, z') < o(\ang{z})$ where the $o()$ decay is as slow as we please. In this case, we consider a curve $s \mapsto \gamma(s)$ between $z = (x,y)$ and $z' = (x',y')$ to be chosen shortly, and compute
\begin{equation}
\nabla_{z''} \big( \Phi_\conic(z'',z) - \Phi_\conic(z'',z') \big) = \int_0^1 \gamma'(s) \cdot \nabla^z \nabla_{z''} \Phi_\conic(z'', \gamma(s)) \, ds.
\label{sss}\end{equation}
The vector $\gamma'(s)$ can be taken to be $a x^2 \partial_{x} + b x \partial_{y_1}$ in a suitable local coordinate system by choosing $\gamma$ appropriately. A short calculation using the exact formula \eqref{sharp-cosine} shows that the length of this curve is comparable to $d_\conic(z,z')$, or more precisely that $a^2 + b^2 = (1 + o(1)) d_\conic(z,z')^2$. Since the matrix of partial derivatives in \eqref{zz''derivs} is equal $I + o(1)$, we thus have
$$
\gamma' \cdot \nabla_z \nabla_{z''} \Phi_\conic(z'', \gamma(s)) = a (x'')^2 \partial_{x''} + b x'' \partial_{y_1} + o(1) \nabla_{x',y'}$$
Thus the length of this vector is equal to $\sqrt{a^2 + b^2} + o(1)$, the integral \eqref{sss} is equal to
$(1 + o(1)) d_\conic(z,z')$. 

To complete the proof of \eqref{grad} for perfectly conic metrics
we have to consider the case when $d_\conic(z,z') > o(\ang{z})$, we may choose
the decay in the $o()$ factor to be as slow as we please.
By \eqref{dm-crude} we thus have $|r-r'| \geq o(\ang{z})$. By \eqref{sharp-cosine} we have
\begin{equation*}\begin{gathered}
|\nabla^{z''} \big(\Phi_\conic(z'',z) - \Phi_\conic(z'',z') \big)|_{g(z'')}^2 \\ = \Big( \frac{\cos d_{\partial M}(y,y'')}{x} - \frac{\cos d_{\partial M}(y',y'')}{x'}  \Big)^2  
 + \Big| \frac{\sin d_{\partial M}(y,y'') \hat \mu}{x}  - \frac{\sin d_{\partial M}(y',y'') \hat \mu'}{x'} \Big|_{h(y'')}^2   \\
=  d_\conic(z,z')^2 + 2rr'T, 
\end{gathered}
\end{equation*}
where $T$ is the quantity
$$
T := \cos d_{\partial M}(y,y') - \cos d_{\partial M}(y,y'') \cos d_{\partial M}(y',y'')   - \sin d_{\partial M}(y,y'') \sin d_{\partial M}(y',y'') \hat \mu \cdot \hat \mu' .
$$
But since $z,z',z''$ are within $o(1)$ of each other in the compactified metric $d_{\overline M}$, we have
$T = o(1)$.  Since $|r-r'| \geq o(\ang{z})$ for some slowly decaying $o()$, we thus have
$rr' T = o(d_\conic(z,z')^2)$, and \eqref{grad} follows.

We now we verify \eqref{gradw} for perfectly conic metrics.  Fix $z, z', \theta$, and let $\Psi(z'')$ denote
the function
$$ \Psi(z'') := (1-\theta) \Phi_\conic(z'',z) + \theta \Phi_\conic(z'',z')$$
appearing in \eqref{gradw}, and set $Z'' := \gamma_{z \to z'}(\theta)$; our task is to show that
\begin{equation}\label{nabla-grow}
 |\nabla_{z''} \Psi(z'')| \geq (1 + o(1)) d(z'', Z'').
\end{equation}
Note that this is already true in the case $z''=Z''$ thanks to \eqref{id}.

Next, let us study the Hessian
$\nabla_{z''} \nabla_{z''} \Phi_{\conic}(z,z'')$, or equivalently let us study the second derivatives
$\frac{d^2}{d\tau^2} \Phi_\conic(z,z''(\tau))$ when $z'' = z''(\tau)$ moves along a geodesic for $g_{\conic}$ 
at unit speed.  We first observe that
$$ \frac{d^2}{d\tau^2} \frac{1}{2} r''(\tau)^2 = 1$$
since the geodesic $z''()$ can be contained inside a flat totally geodesic sector ${\mathcal C} \gamma$,
at which point the claim follows from the corresponding claim for the Euclidean plane $\R^2$.  From 
the cosine rule \eqref{sharp-cosine}
we thus see that
$$ \frac{d^2}{d\tau^2} \Phi_\conic(z,z''(\tau)) = 1- r \frac{d^2}{d\tau^2}
(r''(\tau) \cos d_{\partial M}(y,y''(\tau))).$$ In flat space, $(r''(\tau)
\cos d_{\partial M}(y,y''(\tau)))$ is a linear function of $\tau$ (it is
the $y$ component of $z''(\tau)$) and so the second derivative vanishes.
For perfectly conic manifolds, this quantity is homogeneous of order 1 in
$z''(\tau)$ and so by working in normal co-ordinates in $\partial M$ around
$y$ we obtain the bound
$$
\frac{d^2}{d\tau^2} \Phi_\conic(z,z''(\tau)) = 1 + o(\frac{r}{r''(\tau)}).
$$
and hence we can control the Hessian as
\begin{equation}\label{nabla-phi-conic}
 \nabla_{z''} \nabla_{z''} \Phi_\conic(z,z'') = I + o(\frac{r}{r''})
\end{equation}
where $I$ is the bilinear form $I := g(z'')^{-1}$ on $T_{z''}^* M$ (i.e.\ the identity matrix, if measured in an
orthonormal frame).  Hence
$$
\nabla_{z''} \nabla_{z''} \Psi(z'') = I + o( \frac{\theta r + (1-\theta) r'}{r''} ) 
= I + o( \frac{R''}{r''} )$$
where the latter estimate follows from \eqref{cycle}.
This allows us to establish \eqref{nabla-grow} in the region $r'' \geq o(R'')$ (where the
$o()$ factor decays sufficiently slowly in $\epsilon$), by integrating
this estimate along the geodesic $\gamma_{Z'' \to z'}$ and using the fact that \eqref{nabla-grow} is already
true at $Z''$.  Note that we can control the radial variable of this geodesic without difficulty
since this geodesic lies in a flat totally geodesic sector ${\mathcal C} \gamma$, and thus behaves like a Euclidean
geodesic.

Now consider what happens if $r'' = o(R'')$, so that $d(z'',Z'') = (1 + o(1)) R''$
by \eqref{dm-crude}. To prove \eqref{nabla-grow} it will suffice to bound
the radial gradient, i.e.\ to prove
\begin{equation}\label{nabla-grow-radial}
 |\partial_{r''} \Psi(z'')| \geq (1 + o(1)) R''.
\end{equation}
Using the cosine rule \eqref{sharp-cosine} again, we have 
\begin{equation}\label{phi-conic-form}
\partial_{r''} \Phi_\conic(z,z'') = r'' - r \cos d_{\partial M}(y,y'') = o(R'') - (1 + o(1)) r
\end{equation}
and hence
$$
\partial_{r''} \Psi(z'') = o(R'') - (1-\theta) (1 + o(1)) r - \theta (1 + o(1)) r'.$$
But by \eqref{cycle} we have
$$(1-\theta) (1 + o(1)) r + \theta (1 + o(1)) r' = (1 + o(1)) R'',$$
and \eqref{nabla-grow-radial} follows.  This completes the proof of \eqref{gradw}, and hence Lemma \ref{grad1} has 
been established for perfectly conic metrics.

Finally, we  show (by perturbative arguments) that the estimates \eqref{det-est}, \eqref{grad}, \eqref{gradw}
for the perfectly conic metric will imply the corresponding estimates for the original metric in the asymptotic region $x = o(1)$.  We recall the error function $e := \Phi - \Phi_\conic$ used in Lemma \ref{metric-smooth}.  
Since \eqref{grad} was already proven for the perfectly conic function $\Phi_\conic$, it suffices by the triangle
inequality to prove that
$$ \big| \nabla_{z''} e(z'',z) - \nabla_{z''} e(z'',z') \big| = o(d_M(z,z'))$$
whenever $d_{\overline M}(z, z'), d_{\overline M}(z,z''), d_{\overline M}(z',z'') = o(1)$
and $x, x', x'' = o(1)$.  But this follows by integrating \eqref{e-mixedderiv} on $\gamma_{z \to z'}$.
Note that a similar (and simpler) argument also gives \eqref{det-est} in the asymptotically conic case.


We next prove \eqref{gradw} for asymptotically conic manifolds in the asymptotic region $x = o(1)$. 
We repeat the argument in the perfectly conic case, except now with some additional error terms that can be controlled
with Lemma \ref{metric-smooth}.  As before, it suffices to prove \eqref{nabla-grow}, which is already true when $z''=Z''$
thanks to \eqref{id}.  By Lemma \ref{metric-smooth}, the replacement of $\Phi_\conic$ with $\Phi$
introduces errors of $o(\frac{r}{r''}) + o(1)$ to \eqref{nabla-phi-conic}, 
and errors of $o(r) + o(r'')$ to \eqref{phi-conic-form},
and none of these errors significantly affect the argument.  
This completes the proof of \eqref{gradw}, and thus of Lemma \ref{grad1}.
\endprf

\subsection{Proof of Lemma \ref{geodesic-lemma}}

We now prove Lemma \ref{geodesic-lemma}.  This lemma is trivially true by compactness arguments if $z, z', z'', z_0$
are restricted to any given compact subset $K_0$ of $M$, so we shall assume instead that we are in the asymptotic region
$x = o(1)$.

The claim \eqref{comparable} already follows from \eqref{cycle}, so we turn to 
the symbol estimates \eqref{symbol}. From \eqref{id} and \eqref{gradw} we already know that
$Z'' := \gamma_{z \to z'}(\theta)$ can be defined implicitly by the equation
\begin{equation} \nabla_{z''} \big( (1-\theta)\Phi(Z'',z) + \theta \Phi(Z'',z') \big) = 0.
\label{id-2}\end{equation}
Then \eqref{symbol} follows from \eqref{comparable} and the implicit
function theorem, or more precisely by applying several powers of $\langle
z \rangle \nabla_z$ and $\langle z' \rangle \nabla_{z'}$ to \eqref{id-2}
and solving recursively for $(\langle z \rangle \nabla_z)^\alpha (\langle
z' \rangle \nabla_{z'})^\beta \gamma_{z \to z'}(\theta)$; we omit the
details.

Finally, the bound \eqref{lipschitz} follows directly from \eqref{grad} and Gauss's Lemma \eqref{gauss-lemma}.
This completes the proof of Lemma \ref{geodesic-lemma}.
\endprf

\section{Appendix: The geometry of long range metrics}\label{longrange}

In this appendix we prove Lemma \ref{metric-smooth} for long-range metrics.  In the course of doing so, we shall also
establish uniqueness of length-minimizing geodesics between points which are close in the compactified metric
$d_{\overline M}$. We may work entirely in the scattering region $\{x = o(1)\}$ since in the near region $K_0$ there
is no distinction between long-range and short-range metrics (and Lemma \ref{metric-smooth} is easy to prove from
normal co-ordinates and a compactness argument in that case anyway).

Recall that we call $g$ a long range scattering metric if it takes the form in some smooth coordinates $(x,y)$ near infinity, with $r = x^{-1}$,
\begin{equation}
g = dr^2 + r^2 (h(y, dy) + r^{-\delta} k_{ij} dy^i dy^j),
\label{lr}\end{equation}
where the $k_{00}$, $k_{0j}$, $k_{ij}$ are symbols of order zero on $\Mbar$, and $\delta > 0$.   We can assume that
$\delta$ is small, and allow all constants (for instance in the $o(1)$ notation) to depend on $\delta$.  
In particular, the principal symbol $\sigma(H)$ of the Hamiltonian $H$ in the scattering co-ordinates
$(x,y,\nu,\mu)$ of the scattering cotangent bundle introduced in Section \ref{sec:doi} becomes
\begin{equation}\label{lr-symbol2}
\sigma(H) = \half \Big( \nu^2 + (h^{ij} + x^\delta k^{ij}) \mu_i \mu_j  \Big).
\end{equation}

Before we prove Lemma \ref{metric-smooth}, we need a number of geometric preliminaries.

\subsection{The blown up product manifold $\MMb$.}

It will be convenient to work on the compactified product manifold-with-boundary
$\overline{M} \times \overline{M} := \{ (z,z'): z,z' \in \overline{M} \}$, which contains in particular
the compactified diagonal $\overline{\Delta_0} := \{ (z,z): z \in \overline{M} \}$.  
We also need to introduce the blown-up double space $M^2_b$ given by
$$
M^2_b := [M^2; (\partial M)^2],
$$ i.e.\ it is the double space $M^2$ with the corner blown up. The corner
is thus removed and replaced with a new hypersurface consisting of the
inward-pointing spherical normal bundle of $(\partial M)^2$. This new
hypersurface will be denoted $\bface$ and is a quarter-circle bundle over
$(\partial M)^2$. Coordinates on $\bface$ are $y, y'$ and $x/(x + x')$,
which varies between $0$ at the left boundary $\lb$ of $M^2_b$ and $1$ at
the right boundary $\rb$ of $M^2_b$ (see figure). 
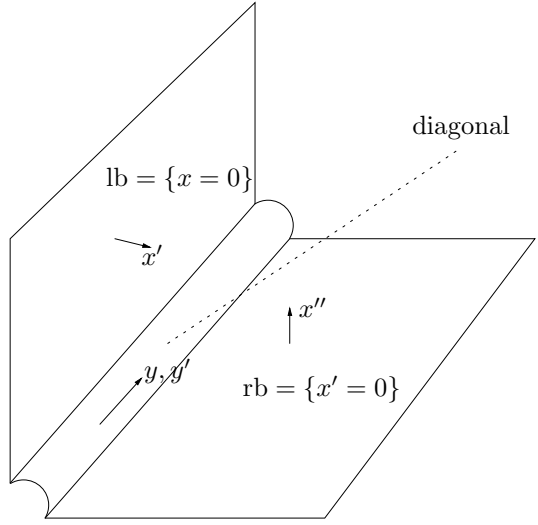
\begin{figure}
\setlength{\unitlength}{0.0003in}
\begingroup\makeatletter\ifx\SetFigFont\undefined%
\gdef\SetFigFont#1#2#3#4#5{%
  \reset@font\fontsize{#1}{#2pt}%
  \fontfamily{#3}\fontseries{#4}\fontshape{#5}%
  \selectfont}%
\fi\endgroup%
{\renewcommand{\dashlinestretch}{30}
\begin{picture}(9024,8889)(0,-10)
\put(237.000,237.000){\arc{874.643}{4.1720}{6.8236}}
\put(4437.000,5037.000){\arc{874.643}{4.1720}{6.8236}}
\path(12,612)(4212,5412)
\path(612,12)(4812,4812)
\path(612,12)(5412,12)(9012,4812)(4812,4812)
\path(12,612)(12,4812)(4212,8862)(4212,5412)
\dashline{60.000}(2712,3012)(7662,6312)
\path(4812,3012)(4812,3612)
\blacken\path(4842.000,3492.000)(4812.000,3612.000)(4782.000,3492.000)(4842.000,3492.000)
\path(1812,4812)(2412,4662)
\blacken\path(2288.307,4662.000)(2412.000,4662.000)(2302.859,4720.209)(2288.307,4662.000)
\path(1554,1625)(2262,2412)
\blacken\path(2204.046,2302.724)(2262.000,2412.000)(2159.440,2342.852)(2204.046,2302.724)
\put(6912,6612){diagonal}
\put(1662,5712){$\lb = \{ x = 0 \}$}
\put(4002,2112){$\rb = \{ x' = 0 \}$}
\put(3762,4512){$\bf$}
\put(4962,3462){$x''$}
\put(2262,4362){$x'$}
\put(2312,2412){$y,y'$}
\end{picture}
}
\caption{The blown up space $\overline{M}^2_b$}\label{fig:m2b}
\end{figure}

\subsection{Algebras of conormal functions}

Since long-range metrics do not smoothly extend to the boundary $\Mbar$, we need
the algebra $\CI(\Mbar)$ by the larger algebra of \emph{conormal functions} $\CId(\Mbar)$. 

In the case of a manifold with boundary, $\CId(\Mbar)$ is defined as a set
to be the sum $\CI(\Mbar) + \Ad(\Mbar)$, where $\Ad(\Mbar)$ is the
Fr\'echet space of conormal functions of order $\delta$ \cite{cocdmc}:
\begin{equation}\label{conormality}
\Ad(\Mbar) = \{ g \in L^\infty(\Mbar) \mid (x\partial_x)^k \partial_y^\alpha g \in 
x^\delta L^\infty(\Mbar) \text{ for all } k, \alpha \}.
\end{equation}
Given a continuous extension operator $E : \CI(\partial \Mbar) \to \CI(\Mbar)$, we may equivalently say that $f \in \CId(\Mbar)$ if $f | \partial \Mbar \in \CI(\partial \Mbar)$ and $f - E(f | \partial \Mbar) \in \Ad(\Mbar)$. This naturally gives rise to a countable family of seminorms on $\CId(\Mbar)$ which makes $\CId(\Mbar)$ into a Fr\'echet space. 

To define the space of conormal functions on the b-double space $\MMb$, we must
we consider a smooth manifold $X$ with corners of codimension two.  To
define $\CId(X)$ we require first that $f \in \CId(X)$ restricts to each
boundary hypersurface $H$ to be in $\CId(H)$. As a consequence the
restriction of $f$ to the codimension two corner is smooth. We then use
smooth extension operators $E_K, E_{H_i}$ from each hypersurface $H$ and
the corner $K$, and require that
$$
f - E_K(f |_{K}) - \sum_i E_{H_i}\big( (f - E_k(F |_{K})) \big) \in \mathcal{A}_{(\delta,\delta)}(X),
$$
where $\mathcal{A}_{(\delta,\delta)}$ is defined by modifying
\eqref{conormality}, with iterated regularity under $x_1\pa_{x_1},$
$x_2 \pa_{x_2},$ and $\pa_y$ leaving the distribution in $x_1^\delta
x_2^\delta L^\infty,$ where $x_i$ are the definining functions for the
boundary hypersurfaces, and $y$ are tangential coordinates.
We give $\CId(X)$ the natural Fr\'echet topology. 

Let $X$ be a manifold with corners of codimension at most two. 
It is not hard to check that $\CId(X)$ is an algebra, and that if $f \in \CId(X)$ is nowhere zero, then also $1/f \in \CId(X)$. When $X = \MMb$, it is straightforward to check that
\begin{equation}
f \in \CId(\MMb) \implies f |_{M^2} \text{ is a local product symbol of order } (0,0).
\label{product-prop}\end{equation}

From now on in this section, by a ``conormal'' function  we mean one in in space $\CId(X)$ where $X$ is the manifold under discussion. By a conormal change of variables on a manifold with boundary, we mean a change of variables of the form $\tilde x(x,y), \tilde y(x,y)$ where $\tilde x(x,y) = x f(x,y)$ and $f$ and $\tilde y_i$ are in $\CId(\Mbar)$; this definition extends in an obvious way to manifolds with corners.

\subsection{Proof of Lemma \ref{metric-smooth} in the long-range case}

We can now prove Lemma \ref{metric-smooth} in the long-range case, by checking
that the proof of \cite[Proposition 9.2]{HTW} extends to the Hamilton vector field associated
to Hamiltonian symbols of the form \eqref{lrsymbol-2}. We assume some familiarity with our arguments in our
previous paper.

As in \cite[Proposition 9.4]{HTW}, we parameterize the cotangent bundle of $\MMb$ by $(\rho,y',x'',y''; \nu', \mu', \nu'',\nu'')$, where $(x',y',\nu',\mu')$ and $(x'',y'',\nu'',\mu'')$ parameterize the two factor cotangent bundles
$\Tscstar \Mbar$, and $\rho := x'/x''$.  We divide the Hamilton vector field of $\sigma(H)$ by $x'$ (i.e.\ we 
reparameterize the time variable $s$ by $ds' := x' ds$, and observe that the vector field now takes 
the form (compare with equations (9.10) and (9.12) of \cite{HTW})
\begin{equation}\begin{split}
\frac{d}{ds'} \rho &= - \rho \nu - \rho^3  (x'')^2 c(\mu, \nu) \\ 
\frac{d}{ds'} {(y')}^i &=  (h^{ij} + (\rho x'')^\delta k^{ij}) \mu'_j + \rho^2 (x'')^2c(\nu) \\ 
\frac{d}{ds'} \nu' &=  (h^{ij} + (\rho x'')^\delta k^{ij}) \mu'_i \mu'_j + \frac{\rho}{2} \frac{\partial (h^{ij} + (\rho x'')^\delta k^{ij})}{\partial \rho} \mu'_i \mu'_j  + \rho^2 (x'')^2 c(\mu \nu, \nu^2)\\ 
\frac{d}{ds'} \mu'_k &= - \mu'_k \nu' - \frac{1}{2}  \frac{\partial (h^{ij}+ (\rho x'')^\delta k^{ij})}{\partial (y')^k} \mu'_i \mu'_j + \rho^2 (x'')^2 c(\nu^2, \nu \mu, \mu^2) \\
\frac{d}{ds'} x'' &=
\frac{d}{ds'} y'' = 
\frac{d}{ds'} \nu'' = 
\frac{d}{ds'} \mu'' = 0 
\end{split} 
\label{hvf}\end{equation}
Here we write $c(\mu)$, $c(\mu \nu)$, etc, for any function which is linear in $\mu$, $\mu \nu$, etc, with coefficients in $\mathcal{A}_0(\Mbar)$. The 
evolution of the function $\kappa := x' \Psi$, where $\Psi(z',z'') := d_M(z',z'')$ is the distance function, is given by 
$$
\frac{d}{ds'} \kappa = 1 - \kappa \nu'  + \kappa \rho^2 (x'')^2 c(\nu, \mu).
$$
Following \cite{HTW}, we introduce $M := \mu/\rho$ and $K := (\kappa - 1)/\rho$ and divide the flow by $\rho$ 
(i.e.\ reparameterize the time variable as $ds'' = \rho ds'$) to obtain
\begin{equation}
\begin{split}
\frac{d}{ds''} \rho &= - \nu + \rho (x'')^2 c(M, \nu) \\
\frac{d}{ds''} {y}^i &= (h^{ij} + (\rho x'')^\delta k^{ij}) M_j + \rho (x'')^2c(\nu)   \\ 
\frac{d}{ds''} \nu &= \rho h^{ij} M_i M_j + \frac{\rho^3}{2}
\frac{\partial (h^{ij} + (\rho x'')^\delta k^{ij})}{\partial \rho} M_i M_j  +  \rho (x'')^2 c(M \nu, \nu^2) \\
\frac{d}{ds''} M_k &=- \frac{1}{2}  \frac{\partial (h^{ij} + (\rho x'')^\delta k^{ij})}{\partial y^k} M_i M_j + x'' c(M^2, M \nu)\\
\frac{d}{ds''} K &= \frac{(h^{ij} + (\rho x'')^\delta k^{ij})M_i M_j}{1 + \nu} + (x'')^2 c(\mu, \nu).\end{split}\label{hvf4}\end{equation}
This is a vector field with conormal coefficients, which is tangent to bf and transverse to $\{ \rho = 0 \}$. We need to prove that the flowout from the set $\{ y' = y'', \rho = 1, \mu' = - \mu'', \nu' = -\nu'' \}$ has conormal regularity. This will follow from

\begin{lemma}\label{conormal-vf} Suppose that $V$ is a conormal, nonvanishing vector field on a manifold with codimension two corners.

(i) Suppose that near the boundary hypersurface $H$, with local coordinates $(x,y)$, $V$ takes the form
$$
V = \sum_i f_i \frac{\partial}{\partial y^i}, \quad f_1 = 1.
$$
Then there is a conormal change of coordinates in which $V$ takes the form 
$$
V = \frac{\partial}{\partial y^1}.
$$

(ii) Suppose that, near a codimension two corner, with local coordinates $(x^1, x^2, y)$, $V$ takes the form 
$$
V = \frac{\partial}{\partial x^1} + \sum_i f_i \frac{\partial}{\partial y^i}.
$$
Then there is a conormal change of coordinates  in which $V$ takes the form 
$$
V = \frac{\partial}{\partial x^1}.
$$
\end{lemma}

\begin{proof} (i)
The flow of this vector field, starting say from the hypersurface $\{ y^1 = 0 \}$, is given by
\begin{equation}
(t, x, y) \mapsto (x, Y_1(t), \dots, Y_m(t)), \quad Y_1(t) = t
\end{equation}
where $Y' = (Y_j)$, $j \geq 2$ satisfies the ODE
\begin{equation}
\frac{dY'}{dy_1}(x, t, y') = f_k(x, t, Y'(x, t, y'))) \quad Y'(x, 0, y') = y'.
\end{equation}
Here we have written $y' = (y_2, \dots, y_m)$. Consequently $Y'$ solves the integral equation
\begin{equation}
Y'(x, y_1, y') = y' + \int_0^{y_1} f_k(x, t, Y'(x, t, y'))) \, dt.
\end{equation}
It follows from this equation that $Y'_j$ are conormal functions. Indeed, by iteratively differentiating the ODE with respect to $x \partial_x$ or $\partial_y$, we see that $Y' \in \mathcal{A}_0(\Mbar)$ locally. To show membership in $\CId(X)$ we derive a linear ODE for the function $x^{-\delta}(Y'(x, y_1, y') - Y'(0, y_1, y'))$ (whose coefficients depend on $Y'$)
and then interatively differentiate the ODE  with respect to $x \partial_x$ or $\partial_y$, as above; we omit the details of this argument which follows standard lines. 

 Then we change to coordinates $(x, y_1, z')$, where\begin{equation}z'(x, y_1, Y'(x, y_1, y')) = y'.\end{equation}Thus the value of $z'$ at $(x, y_1, Y')$ is the value of $y'$ at $y_1 = 0$ along the integral curve passing through $(x, y_1, Y')$. It follows that the Jacobian determinant of $z'$ as a function of $Y'$ is conormal (using the algebra and inverse properties of $\CId(X)$). 
Using a similar argument as above, one checks that $z'$ is also conormal. Finally $z'$ is constant along the flow, so $V$ takes the required form in the coordinates $(x, y_1, z')$. 

The proof of (ii) is similar.
\end{proof}

Continuing the proof of Lemma~\ref{metric-smooth}, we consider the flowout $G$
from the set $\Ginit := \{ y' = y'', \rho = 1, \mu' = -\mu'', \nu' = -
\nu'', g' = g'' = 1 \}$ (with $g',$ $g''$ denoting the metric function
$2\sigma(H)$ in each factor); $G$ is the graph of $d\Psi$, where $\Psi$ is
the distance function $d(z',z'')$, and the value of $\Psi$ is given by
$(\rho x'')^{-1} \kappa = (\rho x'')^{-1} + (x'')^{-1} K$ on $G$. By
symmetry of the distance function, we may restrict to the set where $x' <
2x''$, that is $\rho < 2$. (In the complementary region, $x'' < 2x'$, we
would look at the flow in the doubly-primed variables instead.) Then the
integral curve joining such a point to $\Ginit$ lies wholly within the set
$x' < Cx''$ for some $C < \infty$ by Lemma~\ref{glc}, so we restrict
attention to this region. This is away from the boundary hypersurface $\rb
= \{ x''/x' = 0 \}$, so we may use boundary defining functions $x''$ for
$\bfc$ and $\rho$ for $\lb$ throughout the discussion.

Near any given point we will have either $(d/ds'') \rho \neq 0$ or
$(d/ds'') (y')^k \neq 0$ for some $k$; for definiteness let us assume that
$(d/ds'') \rho \neq 0$. Then by dividing the flow equations by the value of
$(d/ds'') \rho$ (which does not change the flowout from $\Ginit$), we may
assume that $\dot \rho = \pm1$ locally. Applying the lemma above, we find
that $y'$ and $\kappa$ are conormal functions of $y'', \mu'', x'', \rho$ on
$G$, provided that $\rho > 0$. Similar reasoning shows that $y'$ and $K$
are conormal functions of these variables near $\rho = 0$. Moreover, the
flow restricted to $\bfc$ is the same as the exact conic flow, and since
the distance function has the explicit form \eqref{sharp-cosine}, we see
that $(y', y'', \rho)$ are coordinates on $G$ restricted to $\{ x'' = 0
\}$, provided that $y'$ and $y''$ are close in the $h$-metric. Also, we
have $dx'' \neq 0$ at the intersection of $\Ginit \cap \{ x'' = 0 \}$, and
the Lie derivative of $dx''$ vanishes under the flow, so
$dx'' \neq 0$ on $G \cap \{ x'' = 0 \}$ in the region under
consideration. Thus  $(y', y'', \rho, x'')$ are coordinates on
$G$ on a small neighbourhood of $\bfc$, when $y'$ and $y''$ are close in
the $h$-metric, which is equivalent to saying that 
\begin{equation}
\text{ $(y', y'', \rho, x'')$
are coordinates on $G$ on $\overline{\Delta_\epsilon} \subset \MMb$ provided $\epsilon$ is small. }
\label{proj}\end{equation}
 It follows that $\kappa$ (or $K$ when $\rho$
is small) is a conormal function of these variables.  Hence,
$$
\Phi = \frac1{2} \Psi^2  = \frac1{2} \Big( \frac{1 + \rho K}{x'} \Big)^2 = \frac{1 + 2\rho K + \rho^2 K^2}{(x')^2}
$$
near $\rho = 0$, where $K$ is conormal. Therefore 
\begin{equation}
\frac1{\ang{z'}\ang{z''}} \Big( \Phi - \frac{ \ang{z'}^2 + \ang{z''}^2}{2} \Big) = \frac1{2} \Big( 2K + \rho K^2 - \rho \Big)
\label{Phi-symbol}\end{equation}
is a conormal function on $\overline{\Delta_\epsilon} \subset \MMb$. By \eqref{product-prop}, the
expression in \eqref{Phi-symbol} is thus a local product symbol of order (0,0).   
Now we prove the error estimates. We have
$$
e = \Phi - \Phi_{\conic} =  \frac{\ang{z'} \ang{z''}}{2} \Big( 2(K - K_{\conic}) + \rho (K^2 - K^2_{\conic}) \Big).
$$
Since $K$ and $K_{\conic}$ agree at $\bfc$, we see from this that $e / \langle z \rangle \langle z' \rangle$
extends smoothly to $\overline{\Delta_\epsilon}$ and vanishes to order $\delta$ at
$\bfc$. This implies all the claims \eqref{e-onederiv}, \eqref{e-mixedderiv}, \eqref{e-zzderiv}, \eqref{e-zzpderiv},
while \eqref{e-zeroderiv} follows from \eqref{dm-compare}.  This proves Lemma \ref{metric-smooth} in the
long-range case.
\endprf

\begin{corollary}\label{uniq}
Given two points $z, z'$ whose distance is $o(1)$ in the 
compactified metric $d_{\overline M}$, there is only one geodesic between them which remains at distance $o(1)$ from $z$.
\end{corollary}

\begin{proof} 
Any geodesic from $z$ to $z'$ that stays in a $o(1)$ neighbourhood of $z$ in the compactified metric $\overline{g}$ is represented by an integral curve of the flow \eqref{hvf4} which is contained in $G$ and remains uniformly close to $\bfc$. If there were more than one, this would contradict \eqref{proj}.
\end{proof}

\subsection{Long-range metrics not in normal form}\label{normal-form}
In this section we consider more general long range metrics, as in Remark~\ref{general-lr-metrics}, and justify the statement made there that such metrics can be brought to the normal form \eqref{ac-metric} where $h_{ij}(y)$ has the form \eqref{ac-lr-metric}.

Consider a more general long range metric which in the asymptotic region takes the form
$$
g = dr^2 + r^2 h_{ij}(y) dy^i dy^j + r^{-\delta} \Big( k_{00} dr^2 + r k_{0j} dr dy^j + r^2 k_{ij} dy^i dy^j \Big),
$$
for some symbols $k_{00}, k_{0j}, k_{ij}$ of order 0.  In particular, the principal symbol of the Hamiltonian $H$
takes the form
\begin{equation}
\sigma(H) = \half \Big((1 + x^\delta k^{00}) \nu^2 + 2 x^\delta k^{0j} \nu \mu_j + (h^{ij} + x^\delta k^{ij}) \mu_i \mu_j \Big).
\label{lrsymbol-1}\end{equation}
Let us sketch a proof that there is a conormal change of variables bringing such a metric to normal form. 
This fact was essentially proved by Joshi and S\'a Barreto \cite{jsb}; their proof is only written for short range metrics, but we now sketch why the argument extends to the long range setting.  By decreasing $\delta$ if necessary
we may assume that $1/\delta$ is not an integer.  We first establish a preliminary
result, which is that given a symbol of the form 
\begin{equation}
\half \Big( (1 + x^{n \delta} k^{00}) \nu^2 + 2 x^{n \delta} k^{0j} \nu \mu_j + (h^{ij} + x^\delta k^{ij}) \mu_i \mu_j  \Big)
\label{lrsymbol-3}\end{equation}
for some integer $n \geq 0$,
that there is a conormal change of coordinates such that in the new coordinates, the symbol becomes
\begin{equation}
\half \Big( (1 + x^{(n+1) \delta} k^{00}) \nu^2 + 2 x^{(n+1) \delta} k^{0j} \nu \mu_j + (h^{ij} + x^\delta k^{ij}) \mu_i \mu_j  \Big)
\label{lrsymbol-4}\end{equation}
where the functions $k^{00}$, $k^{0j}$, $k^{ij} \in \mathcal{A}_0(\Mbar)$ may of course change from 
\eqref{lrsymbol-3} to \eqref{lrsymbol-4}. 

Given a change of coordinates
\begin{equation*}
\tilde x = x \alpha(y,x), \quad \tilde y = \beta(y,x),
\end{equation*}
the variables $\nu$ and $\mu$ change according to
\begin{equation}
\nu = \frac{\tilde \nu}{\alpha} + \frac{\tilde \nu}{\alpha^2} x \partial_x \alpha +  \tilde \mu_i x \partial_x \beta^i, \quad 
\mu_i = -\frac{\tilde \nu}{\alpha^2} \partial_{y^i} \alpha + \tilde \mu_j \partial_{y^i} \beta^j.
\end{equation}
We try to choose $\alpha$ and $\beta$ in the form 
\begin{equation}
\alpha = 1 + x^{n \delta} a, \quad \beta^i = y^i + x^{n \delta} b^i,
\quad a, b^i \in \mathcal{A}_0(\Mbar),
\end{equation}
hence $\alpha, \beta \in \CId(\Mbar).$
Substituting these expressions into \eqref{lrsymbol-3}, we find that we can arrange \eqref{lrsymbol-4} provided that 
\begin{equation}\begin{gathered}
2 \big( x \partial_x a   + (-1 + n\delta) a \big) = k^{00}, \\
x \partial_x b^i + n\delta b^i = k^{0i} - h^{ij} \partial_{y^i} a.
\end{gathered}\end{equation}
These have solutions $a, b^i \in \mathcal{A}_0(\Mbar)$. A finite number of applications of this result allows us to
represent the symbol $\sigma(H)$ of $H$ in the form
there is a conormal change of co-ordinates under which the symbol $\sigma(H)$ takes the form
\begin{equation}
\sigma(H) = \half \Big( (1 + x^2 k^{00}) \nu^2 + 2 x^2 k^{0j} \nu \mu_j + (h^{ij} + x^\delta k^{ij}) \mu_i \mu_j  \Big).
\label{lrsymbol-2}\end{equation}
(for instance).  Now one applies the arguments of Joshi-S\'a
Barreto involves starting at the submanifold $S_\epsilon = \{ x = \epsilon
\}$ and showing that Fermi normal coordinates with respect to $S_\epsilon$
extend all the way to infinity. This argument goes through in our
setting; effectively we use a  
single space version of the argument in Lemma~\ref{metric-smooth}, introducing
the blowup coordinate $M = \mu/x$ and dividing by a power of $x$ to make
the flow transverse to the boundary. It is then a simple matter to show
that the new $y$ coordinates generated by the flow (i.e.\ the $Y$
coordinates of the proof of Lemma~\ref{conormal-vf}) remain nonsingular to
$x=0$, provided that $\epsilon$ is sufficiently small. We omit the details.

\section{Appendix: Basic estimates for LSIOs}\label{FIO-sec}

In this section we establish the basic estimates for LSIOs that were needed
in the paper, or more specifically Theorem \ref{sfio-l2}, Theorem
\ref{sfio-concatenate}, and Theorem \ref{sfio-infty}.  In our arguments 
we will not bother to explicitly establish the fact that our constants
depends only on finitely many of the product symbol constants for the symbols $b$ or $b'$, but such an estimate can 
be extracted from the argument we give.  Henceforth we allow all our constants to depend on $b$ and $b'$.

\subsection{$L^2$ boundedness.} \label{l2boundedness}

We begin with the proof of Theorem \ref{sfio-l2}.  

\begin{proof}
We first observe that it suffices to prove this result assuming that
$b(z,z',t)$ is supported in the region where $\ang{z} > \ang{z'}/2$ (using
the symmetry of the distance function).  Under this assumption, we use the
$T^* T$ method and show that $V(t)^* V(t)$ is bounded on $L^2$. (In the
region where $\ang{z'} > \ang{z}/2$, we would simply look at $V(t) V(t)^*$
instead.)  Note further that, given that we will show $L^2$ bounds depends
only upon a finite number of symbol seminorms, we may in fact drop the
time-dependence of the symbol.

Consider the kernel of $V(t)^* V(t)$. This is given by an 
integral\footnote{This is an oscillatory integral, not a convergent 
integral, and is given meaning by inserting a cutoff 
$\chi(\ang{z''}/R)$ into the integral, and showing that the limit as 
$R \to \infty$ exists as a distribution.
We omit the standard details.}
\begin{equation}
t^{-n} \int e^{-i\Phi(z'',z)/t} e^{i\Phi(z'',z')/t} 
\overline{b(z,z'')}  b(z'',z') a(z,z'') a(z'',z) \, dg_{z''}.
\label{comp-int}\end{equation}
Because of the $\chi$ cutoff in the kernel, we may suppose (by Lemma \ref{glc} and a 
partition of unity) that $z, z'$ and $z''$ lie in a small ball $B$ in $\Mbar$
of some small radius $\epsilon$ in the compactified metric $d_{\overline M}$. If 
$B$ is contained compactly within the interior of $M$, 
then the result follows by the standard H\"ormander-Eskin $L^2$ bound \cite{hormander}, \cite{eskin} 
for non-degenerate oscillatory integrals (Fourier integral 
operators) (see e.g. \cite{stein:large} for a proof).  Hence we may assume that $B$ is a ``thin cone''; i.e., $y, 
y'$ and $y''$ lie in a ball of radius $o(1)$ in $\partial \Mbar$, and $\ang{z}$, 
$\ang{z'}$ and $\ang{z''}$ are larger than $1/o(1)$. Under this 
assumption, we shall show that $V^*(t) V(t)$  is a semiclassical 
pseudodifferential operator of order zero (with $t$ playing the role 
of small parameter, usually denoted $h$), which proves the lemma 
since it is a standard result that such operators are uniformly 
bounded on $L^2$ \cite{melrose}, \cite{VZ}. 

From Lemma \ref{metric-smooth} we see
that $\nabla_{z,z''} \Phi(z,z'')$ is a local product symbol of order (0,0).
%
In particular on the support of $b$, where $\ang{z''} > \ang{z}/2$, we have
\begin{equation}
\big| \nabla_z \nabla_{z''}^\beta \nabla_{z,z''} \Phi(z,z'') \big| 
\leq C_{\beta}  \ang{z''}^{1-\beta}.
\label{d-symbol}\end{equation}
Let us fix a co-ordinate chart on $B$ (e.g. by using $x$ and $y^j$),
recalling that the measure $a(z,z'') a(z'',z) \, dg_{z''}$ is invariant under
changes of co-ordinates.
With $z(s)$ the point given in coordinates by $z(s) = (1-s) z + s z'$, we 
define a function $\zeta = \zeta(z, z', z'')$ implicitly by setting
\begin{equation*}
\Phi(z'',z) - \Phi(z'',z')  = (z - z') \cdot \zeta, \qquad
\zeta = \int_0^1   \frac{\partial}{\partial z_j} \Phi(z'', 
z(s))  \, ds,.
\end{equation*}
This allows us to write the phase in pseudodifferential form,
$(z-z') \cdot \zeta  / t$. From \eqref{d-symbol} and \eqref{det-est} we have
\begin{equation}
\frac{\partial \zeta_i}{\partial z''_j} = I + o(1)
\label{mixt-derivs}\end{equation}
where $I$ denotes the identity operator on $\RR^n$. Thus we may invert and express $z''$ 
as a function of $(z, z', \zeta)$, $z'' = Z''(z,z',\zeta)$.

We must now show that the function $s(z,z',\zeta)$ defined by
$$
s(z, z', \zeta) := \overline{ab(z,Z''(z, z',\zeta))} ab(z', 
Z''(z,z',\zeta)) \chi(z,Z''(z, z',\zeta)) \chi(z',Z''(z, z',\zeta))
$$ 
is a classical pseudodifferential symbol of order zero 
in $\zeta$; explicitly, we must show that
\begin{equation}
\Big| \big( \ang{z} \frac{\partial}{\partial z} \big)^\alpha \big( 
\ang{z'} \frac{\partial}{\partial z'} \big)^\beta \big( \ang{\zeta} 
\frac{\partial}{\partial \zeta} \big)^\gamma s(z, z', \zeta) \Big| 
\leq C_{\alpha, \beta, \gamma}.
\end{equation}
Since
$Z''(z, z', \zeta(z, z',z'')) = z''$, we see that $\partial_z Z'' = 
\partial_{z'} Z'' = 0$. Thus $z$ or $z'$ derivatives only hit the 
first argument of $ab$ and symbolic estimates for these derivatives 
follow directly from the symbolic estimates on $ab$.

As for the $\zeta$ partial derivatives, we first need good bounds on partial $\zeta$-derivatives of $Z''$. It is straightforward to show 
by induction that
\begin{equation}
\frac{\partial^{k+1} Z''}{\partial \zeta^{k+1}} = \sum_{(\alpha_1, 
\dots, \alpha_l)}^{\alpha_1 + \dots \alpha_l = k} \prod_{j=1}^l \big( 
\frac{\partial Z''}{\partial \zeta} \big)^{b_j} 
\frac{\partial^{\alpha_j+1} \zeta}{(\partial z'')^{\alpha_j+1}}.
\label{Z''-derivs}\end{equation}
Here the product on the right hand side is shorthand for a sum of 
terms  each of which is a monomial of degree $b_j$ in the components 
of $\partial Z'' / \partial \zeta$ and linear in the components of
$\partial^{\alpha_j + 1} \zeta / (\partial z'')^{\alpha_j + 1}$.
By \eqref{mixt-derivs}, the matrix $\partial Z''/\partial \zeta$ is 
bounded between $1/2$ and $2$, while from the definition of $\zeta$ 
and \eqref{d-symbol}, $\partial \zeta / \partial z''$ is 
a symbol of order zero in $z''$, uniformly in $z, z' \in U$. It 
follows from these facts and \eqref{Z''-derivs} that
\begin{equation*}
\Big| \frac{\partial^{k+1} Z''}{\partial \zeta^{k+1}} \Big| \leq C_k 
\ang{z''}^{-k}.
\end{equation*}
Recall now that the integrand in \eqref{comp-int} is supported  where $\ang{z''}
> \max(\ang{z}/2,\ang{z'}/2)$. 

We now note that from \eqref{ac-metric} (and \eqref{ac-lr-metric}) that
the geodesic evolution of the radial co-ordinate $r$ is
\begin{equation*}
\ddot r = r^{-1} \big( h^{ij} + o(1) \big) \mu_i \mu_j;  \end{equation*}
in particular, since $h^{ij}$ is positive definite, we see that $r(s)$ is a convex function of $s$
in the asymptotic region.
Using this convexity, and by comparing the metric with a larger
conic metric and using the exact formula \eqref{sharp-cosine} for conic
distance, we have
\begin{equation*}
|\zeta|^2 \leq \sup_s d^2(z'', z(s)) \leq \ang{z''}^2 + \big(\ang{z} 
+ \ang{z'} \big)^2 + \eta \ang{z''}\big(\ang{z} + \ang{z'} \big) \leq 
C \ang{z''}^2.
\end{equation*}
Hence we find that
\begin{equation*}
\Big| \frac{\partial^{k+1} Z''}{\partial \zeta^{k+1}} \Big| \leq C_k 
\ang{\zeta}^{-k},
\end{equation*}
and using this it is simple to check that $s$ is symbolic of order 
zero in $\zeta$. This completes the proof of the lemma.
\end{proof}

\subsection{The composition law.}

We now prove Theorem \ref{sfio-concatenate}.  Fix $b, b'$.  From \eqref{S-def} we have
$S_b(t) S_{b'}(t') = S_{\tilde c_{t,t'}}(t+t')$, where $\tilde c_{t,t'}$ is the function
\begin{multline*} \tilde c_{t,t'}(z,z'') := \\ \frac{(2\pi i (t +
    t'))^{n/2}}{a(z,z'')(2\pi i t)^{n/2} (2\pi i t')^{n/2}} 
\int\limits_M e^{i( \frac{\Phi(z,z')}{t} + \frac{\Phi(z',z'')}{t'} - \frac{\Phi(z,z'')}{t+t'}} 
ab(z,z',t) ab'(z',z'',t')\ dg(z').\end{multline*}
Since $b$, $b'$ are localized to a region of the form $\Delta_\epsilon$ for some small $\epsilon$,
so is $\tilde c_{t,t'}$ (with a slightly larger value of $\eps$. It thus suffices to show that $\tilde c_{t,t'}$
has the expansion $\tilde c_{t,t'} = c + \frac{t+t'}{\langle z \rangle + \langle z' \rangle}  e_{t,t'}$ 
for some local product symbol $e_{t,t'}$ of order (0,0), and $c$ given by \eqref{b-weight}; observe from \eqref{symbol} that this would in fact ensure that $\tilde c_{t,t'}$ is a product local symbol of order (0,0) as claimed.

Let us make the ``semi-classical'' change of variables $\hbar :=
\frac{tt'}{t+t'}$ and $\theta := \frac{t}{t+t'},$ and set
$$f_{\hbar,\theta}=a(z,z'')\tilde c_{t,t'}$$ so that we have
$$ f_{\hbar,\theta}(z,z'') = \frac{1}{(2\pi i \hbar)^{n/2}} 
\int_M e^{\frac{i}{\hbar} \Psi_\theta(z,z',z'')} 
ab(z,z',t) ab'(z',z'',t')\ dg(z')$$
where $\Psi_\theta$ is the phase
$$ \Psi_\theta(z,z',z'') := (1-\theta) \Phi(z,z') + \theta\Phi(z',z'') -
\theta(1-\theta)\Phi(z,z'').$$  Here we have elected not to apply the
change of variables to the $t,t'$ in the symbols (one may think of these as
auxiliary parameters) so as not to require any time-regularity of the
symbols later on.

Our task is to show that
\begin{equation}\label{c-compose}
 f_{\hbar,\theta}(z,z'') = a(z,z'') b(z,z'_0,t) b'(z'_0,z'',t') + \frac{\hbar}{\theta (1-\theta) (\langle z \rangle + \langle z'' \rangle)} e_{\hbar,\theta}(z,z'')
\end{equation}
where $e_{\hbar,\theta}$ is a local product symbol of order (0,0) which
depends on $\theta$ and $\hbar$ but obeys local product symbol bounds which
are uniform with respect to these parameters, and $z'_0 := \gamma_{z \to
  z''}(\theta)$.  (We have absorbed the $a$ factor into $e_{\hbar,\theta}$
since $a$ and $a^{-1}$ are both local product symbols of order $(0,0)$ on
the support of $f_{\hbar,\theta},$ thanks to Lemma~\ref{a-smooth}.)

We first verify the claim in the classical limit $\hbar = 0$, or more precisely
that
\begin{equation}\label{limit} 
\begin{split}
\lim_{\hbar \to 0}
\frac{1}{(2\pi i \hbar)^{n/2}} &
\int_M e^{\frac{i}{\hbar} \Psi_\theta(z,z',z'')} 
ab(z,z',t) ab'(z',z'',t')\ dg(z')
\\
& = a(z,z'') b(z,z'_0,t) b'(z'_0,z'',t').
\end{split}
\end{equation}
From Lemma \ref{grad1} we have the phase oscillation estimate
\begin{equation}\label{phase-osc} |\nabla_{z'} \frac{\Psi_\theta(z,z',z'')}{\hbar}|
\geq c \frac{1}{\hbar} d(z', z'_0).
\end{equation}
In particular the only stationary point of the integral occurs at $z' = z'_0$ (cf. \eqref{id}).  At this point, the phase is equal to
\begin{equation}\label{phase}
 \Psi_\theta(z,z'_0,z'') =
(1-\theta)(\theta d(z,z''))^2/2 + \theta ((1-\theta) d(z,z''))^2/2 - \theta(1-\theta) d(z,z'')^2/2 = 0.
\end{equation}
Applying the principle of stationary phase (see \cite{stein:large}), it
thus suffices to verify the identity
$$ \det( \nabla_{z'}^2 \Psi_\theta(z,z',z'')|_{z' = z'_0})^{-1/2} =
\frac{a(z,z'')}{a(z,z'_0) a(z'_0,z'')\ dg(z'_0)}$$ or (by \eqref{a-def})
that
\begin{equation}\label{phase-magic}
 \det( \nabla_{z'}^2 \Psi_\theta(z,z',z'') |_{z' = z'_0})
= \frac{ \det( - \nabla_z \nabla_{z'} \Phi(z,z') )  
\det( - \nabla_{z'} \nabla_{z''} \Phi(z',z'') )  
} { \det( - \nabla_z \nabla_{z''} \Phi(z,z'') )  }|_{z' = z'_0}.
\end{equation}
To see this, we now let $z$ and $z''$ vary (keeping $\theta$ fixed), so that $z'_0$ is now thought of as a function
of $z$ and $z''$.  We rewrite \eqref{phase} as
$$ \theta(1-\theta) \Phi(z,z'') = [(1-\theta) \Phi(z,z') + \theta \Phi(z',z'')] |_{z' = z_0}$$
and observe from \eqref{id} that the right-hand side is stationary in $z'$ at $z_0$:
\begin{equation}\label{gauss}
\nabla_{z'} \Psi(z,z',z'')|_{z' = z'_0}
= \nabla_{z'} [(1-\theta) \Phi(z,z') + \theta \Phi(z',z'')] |_{z' = z'_0} = 0,
\end{equation}
Differentiating both
sides with respect to $z''$ (and using the previous observation to ignore the $z'$ variation), we thus obtain
$$ \theta(1-\theta) \nabla_{z''} \Phi(z,z'') = \theta \nabla_{z''} \Phi(z'_0, z'').$$
Differentiating again, this time in $z$, and using the chain rule, we obtain
$$ \theta(1-\theta) \nabla_z \nabla_{z''} \Phi(z,z'') = \theta
\nabla_{z'} \nabla_{z''} \Phi(z', z'')|_{z' = z'_0}
\cdot \nabla_z z'_0(z,z'').$$
On the other hand, if we differentiate \eqref{gauss} in $z$ using the chain rule, we obtain
$$ ((1-\theta) \nabla_{z} \nabla_{z'} \Phi(z,z')
+ \nabla_{z'}^2 \Psi(z,z',z'') |_{z' = z'_0})
\cdot \nabla_z z'_0(z,z'') = 0.$$
Combining this with the previous identity, we obtain after some algebra
\begin{align*} -\nabla_z \nabla_{z''} \Phi(z,z'') =
&[-\nabla_{z'} \nabla_{z''} \Phi(z', z'')\\
&\cdot \nabla_{z'}^2 \Psi(z,z',z'')^{-1}\\
&\cdot (-\nabla_{z} \nabla_{z'} \Phi(z,z'))]|_{z' = z'_0}.
\end{align*}
Taking determinants of both sides we obtain \eqref{phase-magic} as
desired\footnote{Note that we did not really use
any Riemannian geometry to prove this lemma, instead relying exclusively on
the fact that the phase \eqref{phase} vanished and was stationary at
$z'_0$.  Indeed this computation is really one in symplectic geometry
rather than Riemannian geometry, and reflects the multiplicative law for
the symbol of Fourier integral operators under composition, and the
untruncated parametrix $S_1(t)$ has symbol $1$ (see Remark~\ref{unitary}).}. This
proves \eqref{limit}.

In light of \eqref{limit}, it will now suffice to obtain a two-term expansion of the form
$$
 f_{\hbar,\theta}(z,z'') = f_{0,\theta}(z,z'') + \frac{\hbar}{\theta (1-\theta) (\langle z \rangle + \langle z'' \rangle)} e_{\hbar,\theta}(z,z'')
$$
where $c_{0,\theta}(z,z'')$ is some
function independent of $\hbar$, since by taking $\hbar \to 0$ we then see that $f_{0,\theta}(z,z'')$ equals
the expression in \eqref{limit}.   In fact by \eqref{comparable}, \eqref{symbol} it will suffice to obtain
a representation of the form
\begin{equation}\label{asymptotic}
 f_{\hbar,\theta}(z,z'') = f_{0,\theta}(z,z'') + \frac{\hbar}{\langle z'_0 \rangle^2} e_{\hbar,\theta}(z,z'')
\end{equation}
for some product symbol $e_{\hbar,\theta}$ of order (0,0), with bounds uniform in $\hbar$ and $\theta$.

Introduce a smooth cutoff $\varphi_{z,z''}(z')$
which equals 1 when $d(z',z'_0) \rangle \leq C^{-1} \langle z'_0 \rangle$ and is supported in the region
$d(z',z'_0) \leq 2C^{-1} \langle z'_0 \rangle$ for some large $C \gg 1$.
Consider first the non-local contribution to $f_{\hbar,\theta}(z,z''):$
\begin{multline*} f_{\hbar,\theta}(z,z'') - \tilde f_{\hbar,\theta}(z,z'')
  \\ :=
\frac{1}{(2\pi i \hbar)^{n/2}} \int_M e^{\frac{i}{\hbar}
\Psi_\theta(z,z',z'')} ab(z,z',t) ab'(z',z'',t') (1 - \varphi_{z,z''}(z'))\
dg(z').\end{multline*}  From \eqref{phase-osc} and repeated
integration by parts we obtain the bounds
$$ |(f_{\hbar,\theta} - \tilde f_{\hbar,\theta})(z,z'')| \leq C_N \hbar^{N} \langle z'_0 \rangle^{-N}$$
for any $N \geq 0$, and in fact we have the additional estimates
$$ |(\langle z \rangle \partial_z)^\alpha (\langle z'' \rangle \partial_{z''})^\beta (f_{\hbar,\theta} - \tilde f_{\hbar,\theta})(z,z'')| \leq C_{N,\alpha,\beta} \hbar^{N} \langle z'_0 \rangle^{-N}$$
for any $\alpha, \beta \geq 0$ (mostly thanks to \eqref{symbol}).  Thus this contribution to \eqref{asymptotic} is
certainly acceptable.

It remains to consider the local contribution
$$ \tilde f_{\hbar,\theta}(z,z'') := \frac{1}{(2\pi i \hbar)^{n/2}} \int_M
e^{\frac{i}{\hbar} \Psi_\theta(z,z',z'')} ab(z,z',t) ab'(z',z'',t')
\varphi_{z,z''}(z')\ dg(z')$$ 
But this contribution will be acceptable, thanks to the standard result in
stationary phase that a oscillatory integral with a single non-degenerate
critical point and with phase and amplitude smooth and smoothly depending
on certain parameters will enjoy an asymptotic series, with all the main
terms and error terms in the series depending smoothly on the
parameters---see Theorem 7.7.6 of \cite{Ho}. 
This concludes the proof of Theorem \ref{sfio-concatenate}.

\qed

\subsection{The dispersive inequality.}

We now prove Theorem \ref{sfio-infty}.  

When one of $s$ or $t$ is zero, the claim follows from \eqref{S0-def} and
\eqref{easy-dispersive}, so assume $s,t$ are non-zero.  When $t$ and $s$
have opposite sign then the claim follows quickly from Lemma
\ref{sfio-basic}, Theorem \ref{sfio-concatenate}, and
\eqref{easy-dispersive}.  Thus we may assume that $t$ and $s$ have the same
sign\footnote{These cases are not symmetric; the case of the equal signs is
both more difficult and more important.  The situation is analogous to that
of the $L^2$ or $L^p$ theory of the Hilbert transform kernel
$\frac{1}{s-t}$; the portion of this kernel when $s$ and $t$ have opposing
signs is easy to treat by Hardy's inequality (even when absolute values are
placed on the kernel), but the portion when $s$ and $t$ have the same signs
is more delicate.}.  By duality and Lemma \ref{sfio-basic} we may in fact
assume that $0 < s < t$.  By expanding out the kernel of $S_b(t)^*
S_{b'}(s)$ using \eqref{S-def}, it thus suffices to show that
\begin{equation}\label{z-race}
 |\int_M e^{i (\frac{\Phi(z_1,z_2)}{s} - \frac{\Phi(z_1,z_0)}{t})} \overline{ab(z_1,z_0,t)} ab'(z_1,z_2,s)\ dg(z_1)|
\leq C \frac{t^{n/2} s^{n/2}}{|t-s|^{n/2}}
\end{equation}
for all $z_0, z_2 \in M$.  Clearly we may restrict (via smooth truncation of $b$ or $b'$) to the region where
\begin{equation}\label{z1-large}
\langle z_1 \rangle \geq C \big(\frac{ts}{|t-s|}\big)^{1/2}
\end{equation}
since the portion of the integral where \eqref{z1-large} fails can be estimated by taking absolute values everywhere.

The idea is to use the composition law (Theorem \ref{sfio-concatenate}) 
to factorize $S_b(t)^* S_{b'}(s)$ into something resembling $S_{b''}(t-s)^* S_{b'''}(s)^* S_{b'}(s)$.  
Indeed, from Theorem \ref{sfio-concatenate} (with $b, b' := \chi$ for some suitable cutoff $\chi$ to a region
$\Delta_\epsilon$) we have
\begin{multline*} \frac{(2\pi i t)^{\frac n 2}}{(2\pi i s)^{\frac n 2} (2\pi i (t-s))^{\frac n 2}} 
\int_M e^{i( \frac{\Phi(z,z')}{s} + \frac{\Phi(z',z'')}{t-s}} 
a(z,z') a(z',z'') \chi(z,z') \chi(z',z'') \ dg(z') \\
= e^{i\Phi(z,z'')/t} (a(z,z'') \chi(z,z') \chi(z',z'') + \frac{t}{\langle z \rangle + \langle z'' \rangle} e(z,z''))
\end{multline*}
for some local product symbol $e$ of order (0,0).  We may arrange matters so that $\chi(z,z')
\chi(z',z'') = 1$ on the support of $b$.  When $t$ is small, or $\langle z
\rangle$ or $\langle z'' \rangle$ is large, the error term
$\frac{t}{\langle z \rangle + \langle z'' \rangle} e(z,z'')$ is dominated
by the main term, and in particular the term in parentheses is bounded away
from zero. When $t$ is large and $\langle z \rangle$ or $\langle z''
\rangle$ is small, then this is not necessarily the case, but this is
easily fixed in this case by modifying the amplitude $a(z,z') a(z',z'')$ on
the left-hand side slightly in the region $\langle z \rangle, \langle z'
\rangle, \langle z'' \rangle = O(1)$, replacing it with a local product symbol
$\tilde a(z, z', z'')$ of order $(0,0,0)$ in the three variables.  We can
then multiply both sides by a suitable local product symbol of order (0,0) in $z$ and $z''$ to
obtain a representation of the form
$$ e^{i\frac{\Phi(z,z'')}{t}} ab(z, z'',t) = \frac{(2\pi i t)^{\frac{n}{2}}}{(2\pi i
s)^{\frac{n}{2}} (2\pi i (t-s))^{\frac{n}{2}}} \int\limits_M e^{i( \frac{\Phi(z,z')}{s} +
\frac{\Phi(z',z'')}{t-s})} \tilde b(z,z',z'',s,t)\ dg(z')$$ for some
local product symbol $\tilde b$ of order $(0,0,0)$, localized so that $z,z',z''$
are all close to each other in the compactified metric $d_{\overline{M}}$.
Note that while $\tilde b$ will depend on both $s$ and $t$, the local product symbol
bounds are uniform in these parameters.  Using this representation, we can
rewrite \eqref{z-race} as\footnote{Semiclassically, what is going on is the
following.  The expression in \eqref{z-race} is sort of a quantum
mechanical version of a trajectory consisting of the concatenation of the
geodesics $\gamma_{z_2 \to z_1}$ and $\gamma_{z_1 \to z_0}$, where
$\gamma'_{z_1 \to z_0}(0) = - \frac{t}{s} \gamma'_{z_2 \to z_1}(0)$ (so the
path $\gamma_{z_1 \to z_0}$ overshoots past $z_2$).  The point $z'$ then
corresponds to the point $\gamma_{z_1 \to z_0}(s/t) = z_2$, and so this
factorization corresponds to a splitting of this path into a loop from
$z_2$ to $z_1$ and back to $z' = z_2$, followed by a geodesic from $z'=z_2$
to $z_0$.}
$$
s^{-n} |\int_M \int_M e^{i (\frac{\Phi(z_1,z_2)}{s} - \frac{\Phi(z_1,z')}{s} - \frac{\Phi(z_0,z')}{t-s})} 
\overline{\tilde b(z_1,z',z_0,s,t)} b'(z_1,z_2,s)\ dg(z_1) dg(z')|
\leq C.$$
Note that all four points $z_0, z', z_1, z_2$ are constrained to lie close to each other in $d_{\overline{M}}$.  

Let $m \geq 0$ be an integer, and denote by $Q_m$ the quantity
$$s^{-n} \Big|\int\limits_{M^2} e^{i (\frac{\Phi(z_1,z_2)}{s} - \frac{\Phi(z_1,z')}{s} - \frac{\Phi(z_0,z')}{t-s})} 
\overline{\tilde b(z_1,z',z_0,s,t)} b'(z_1,z_2,s) \varphi_m(z_1)\ dg(z_1) dg(z') \Big|$$
where $\varphi_m$ is a suitable smooth cutoff to the region where $\langle z_1 \rangle$ is comparable to $2^m$.
From \eqref{z1-large} we may assume that $2^{2m} \geq st/(t-s)$.
It thus suffices to show that 
\begin{equation}\label{qm-bound}
\sum_{m \geq 0: 2^{2m} \geq st/(t-s)} Q_m \leq C.
\end{equation}
Let us first observe that each $Q_m$ is individually bounded.
Indeed, from \eqref{grad} we observe the phase oscillation estimate
$$ \Big|\nabla_{z_1} (\frac{\Phi(z_1,z_2)}{s} - \frac{\Phi(z_1,z')}{s} -
\frac{\Phi(z_0,z')}{t-s})\Big| \geq C^{-1} d(z_2,z')/s,$$
while for higher derivatives we see from \eqref{symbol} and Lemma \ref{metric-smooth} that
$$ \Big|\nabla^k_{z_1} (\frac{\Phi(z_1,z_2)}{s} - \frac{\Phi(z_1,z')}{s} -
\frac{\Phi(z_0,z')}{t-s})\Big| \geq C_k \langle z_1 \rangle^{1-k}/s \hbox{ for } k \geq 2.$$
These estimates imply
(after repeated integration by parts in $z_1$, and exploiting the local product symbol
bounds for $\tilde b$, $b'$ as well as the localization of $\varphi_m$)
that for fixed $z' \in M$ we have
\begin{equation}\label{z1-int}
 \Big|\int_M e^{i (\frac{\Phi(z_1,z_2)}{s} - \frac{\Phi(z_1,z')}{s} - \frac{\Phi(z_0,z')}{t-s})} 
\overline{\tilde b} \,  b' \, \varphi_m 
\ dg(z_1) \Big| \leq C_N 2^{nm} (1 + 2^m d(z_2,z')/s)^{-N}
\end{equation}
for any $N > 0$; taking $N = 100n$ (for instance) and integrating this estimate over $z'$, we obtain
the bound $Q_m = O(1)$ as claimed.  A similar argument shows that the contribution of the integral to $Q_m$ coming
from the region where $\langle z' \rangle$ is not comparable to $\langle z_2 \rangle$ is extremely small (e.g.
it is bounded by $O(2^{-100m})$ so we shall discard this portion and smoothly truncate to the region where $\langle z' \rangle$ is comparable to $\langle z_2 \rangle$.

To improve upon these bounds and obtain \eqref{qm-bound}
we have to take advantage of oscillation in the $z'$ variable as well as the $z_1$ 
variable.  
From \eqref{gradw} we have
\begin{equation}\label{gw}
 \Big|\nabla_{z'} (\frac{\Phi(z_1,z_2)}{s} - \frac{\Phi(z_1,z')}{s} - \frac{\Phi(z_0,z')}{t-s})\Big|
\geq C^{-1} d_M(z',\gamma_{z_1 \to z_0}(\frac{s}{t})) \frac{t}{s(t-s)}.
\end{equation}
Since we have already localized $z'$ to be close to $z_2$, it is now natural to introduce the quantity
$$ \rho_m := \inf_{z_1 \in \supp \varphi_m} d_M(z_2,\gamma_{z_1 \to z_0}(\frac{s}{t})).$$
This quantity is usually quite large:

\begin{lemma}  We have $\rho_m \geq C^{-1} \frac{t-s}{t} 2^m$ for all but $O(1)$ values of $m$.
\end{lemma}

\begin{proof}  Let $m \geq \tilde m$ be any two integers such that $\rho_m \leq C_0^{-1} \frac{t-s}{t} 2^m$ and
$\rho_{\tilde m} \leq C_0^{-1} \frac{t-s}{t} 2^{\tilde m}$ for some large constant $C_0 > 1$ to be chosen later.  By definition of $\rho_m$ and the triangle inequality, we can thus find $z_1 \in \supp \varphi_m$, $\tilde z_1 \in \supp \varphi_{\tilde m}$ such that
$$ d_M( \gamma_{z_1 \to z_0}(\frac{s}{t}), \gamma_{\tilde z_1 \to z_0}(\frac{s}{t}) ) \leq 2 C_0^{-1}\frac{t-s}{t} 2^m.$$
But by \eqref{lipschitz} , the left-hand side is comparable to
$$ \frac{t-s}{t} |\gamma'_{z_0 \to z_1}(0) - \gamma'_{\tilde z_0 \to z_1}(0)|$$
which by \eqref{lipschitz} again is comparable to
$$ \frac{t-s}{t} d_M( z_1, \tilde z_1 ).$$
Thus $d_M(z_1, \tilde z_1) \leq C C_0^{-1} 2^m$, which by the localization of $z_1, \tilde z_1$ forces $m = m' + O(1)$
if $C_0$ is chosen large enough.  The claim follows.
\end{proof}

We may assume that 
\begin{equation}\label{rhom}
\rho_m \geq C^{-1} \frac{t-s}{t} 2^m,
\end{equation}
since the contribution of the $O(1)$ exceptional values of $m$ to \eqref{qm-bound} can be dealt
with by the bound $Q_m = O(1)$ obtained previously.  Since $2^{2m} \geq st/(t-s)$, this implies that
\begin{equation}\label{rhomba}
\rho_m \geq C^{-1} s / 2^m. 
\end{equation}
Now we smoothly truncate the $z'$ integral into the
region where $d_M(z',z_2) \leq \min(1,\rho_m)/2$ and where $d_M(z',z_2) \geq \min(1,\rho_m)/4$.  The contribution of
the latter case to $Q_m$ is at most $C_N (2^m \min(1,\rho_m) / s)^{-N}$ for any $N \geq 0$ thanks to \eqref{z1-int}
and \eqref{rhomba}, which then gives an acceptable contribution to \eqref{qm-bound}
thanks to \eqref{rhom}.  For the former case, we observe from \eqref{gw} and the triangle inequality that
$$ \Big|\nabla_{z'} (\frac{\Phi(z_1,z_2)}{s} - \frac{\Phi(z_1,z')}{s} -
\frac{\Phi(z_0,z')}{t-s})\Big| \geq C^{-1} \frac{\rho_m}{2}
\frac{t}{s(t-s)},$$ 
while higher derivatives are controlled using Lemma \ref{metric-smooth} and \eqref{symbol} by
the somewhat crude estmate
$$ \Big|\nabla^k_{z'} (\frac{\Phi(z_1,z_2)}{s} - \frac{\Phi(z_1,z')}{s} -
\frac{\Phi(z_0,z')}{t-s})\Big| \leq C_k \frac{t}{s(t-s)} \hbox{ for all } k \geq 2$$ 
so by repeated integration by parts
in $z'$ (noting that every derivative in $z'$ of 
a cutoff function or of phase derivatives costs us at most $O( 1 / \min(1,\rho_m)$) we see that the
contribution of this case to $Q_m$ is at most $$C_N s^{-n} 2^{nm}
\min(1,\rho_m)^n (C \rho_m \min(1,\rho_m) \frac{t}{s(t-s)} )^{-N}$$ for any
$N \geq 0$, which is in turn bounded by $C_N (2^m \min(1,\rho_m) / s
)^{-N + O(1)}$ thanks to \eqref{rhomba}, and so as before
this gives an acceptable contribution to \eqref{qm-bound}.  This completes
the proof of Theorem \ref{sfio-infty}.  \qed


\end{document}